\documentclass[opre,nonblindrev]{informs4_arxiv}

\usepackage{amsmath,amssymb,epsfig,graphicx,float,url,tabularx,multirow}
\usepackage{algorithm,algorithmicx,algpseudocode}
\usepackage[hypertexnames=false,hyperfootnotes=false]{hyperref}
\usepackage[hmargin=1in,vmargin=1in]{geometry}
\usepackage{float}

\usepackage[sort&compress]{natbib}
 \bibpunct[, ]{[}{]}{,}{n}{}{,}%
 \def\bibfont{\small}%

\usepackage[dvipsnames]{xcolor}
\usepackage{subcaption}
\usepackage{wrapfig}
\usepackage{tabularx}

\usepackage{mathabx}  

\OneAndAHalfSpacedXI 

\newcommand{\opt}{^{\star}}
\newcommand{\tb}{\textcolor{black}}

\usepackage[colorinlistoftodos,bordercolor=orange,backgroundcolor=orange!20,line
color=orange,textsize=scriptsize]{todonotes}
\catcode`\@=11
\catcode`\@=12

\DeclareMathOperator{\kl}{KL}

\newtheorem{theorem}{Theorem}
\newtheorem{proposition}[theorem]{Proposition}

\newtheorem{remark}[theorem]{Remark}
\newtheorem{example}[theorem]{Example}
\newtheorem{corollary}[theorem]{Corollary}
\newtheorem{definition}[theorem]{Definition}
\newtheorem{lemma}[theorem]{Lemma}
\newtheorem{assumption}[theorem]{Assumption}

\newcommand{\ep}{\varepsilon}

\newcommand{\E}{\mathbb{E}}

\newcommand{\cH}{\mathcal{H}}

\newcommand{\tr}{^{\top}}

\newcommand{\X}{\mathcal{S}}

\newcommand{\A}{\mathcal{A}}
\newcommand{\U}{\mathcal{U}}
\newcommand{\N}{\mathbb{N}}
\newcommand{\R}{\mathbb{R}}
\newcommand{\gammab}{\gamma_{{\sf bw}}}
\newcommand{\gammaa}{\gamma_{{\sf avg}}}

\newcommand{\tbm}{{\em The Big Match}}

\newcommand{\Ravg}{R_{{\sf avg}}}  
\newcommand{\hRavg}{\hat{R}_{{\sf avg}}}

\newcommand{\PiH}{\Pi_{\sf H}}
\newcommand{\PiS}{\Pi_{\sf S}}
\newcommand{\PiSD}{\Pi_{\sf SD}}

\newcommand{\PiM}{\Pi_{\sf M}}

\newcommand{\err}{{\sf error}}
 
\usepackage{bm}
\numberwithin{equation}{section}
\numberwithin{theorem}{section}

\ECRepeatTheorems
\EquationsNumberedBySection

\hypersetup{colorlinks=true,citecolor=Blue,linkcolor=MidnightBlue,urlcolor=Green}

\RUNAUTHOR{Grand-Cl{\'e}ment, Petrik and Vieille}

\RUNTITLE{Robust MDPs with average and Blackwell optimality}

\TITLE{Beyond discounted returns: Robust Markov decision processes with average and Blackwell optimality}

\ARTICLEAUTHORS{%
\AUTHOR{Julien Grand-Cl{\'e}ment}
\AFF{Information Systems and Operations Management Department, HEC Paris, \EMAIL{grand-clement@hec.fr}}
\AUTHOR{Marek Petrik}
\AFF{Department of Computer Science, University of New Hampshire, \EMAIL{mpetrik@cs.unh.edu}}
\AUTHOR{Nicolas Vieille}
\AFF{Economics and Decision Sciences Department, HEC Paris, \EMAIL{vieille@hec.fr}}
} 

\ABSTRACT{%
  Robust Markov Decision Processes (RMDPs) are a widely used framework for sequential decision-making under parameter uncertainty. RMDPs have been extensively studied when the objective is to maximize the discounted return, but little is known for average optimality (optimizing the long-run average of the rewards obtained over time) and Blackwell optimality (remaining discount optimal for all discount factors sufficiently close to $1$). In this paper, we prove several foundational results for RMDPs beyond the discounted return. We show that average optimal policies can be chosen stationary and deterministic for sa-rectangular RMDPs \tb{but, perhaps surprisingly, we show that for s-rectangular RMDPs average optimal policies may not exist, and if they exist, may need to be history-dependent (Markovian)}. We also study Blackwell optimality for sa-rectangular RMDPs, where we show that $\epsilon$-Blackwell optimal policies always exist, although Blackwell optimal policies may not exist. We also provide a sufficient condition for their existence, which encompasses virtually any examples from the literature. We then discuss the connection between average and Blackwell optimality, and we describe several algorithms to compute the optimal average return. Interestingly, our approach leverages the connections between RMDPs and stochastic games. \tb{Overall, our paper emphasizes the superior practical properties of distance-based sa-rectangular models over s-rectangular models for average and Blackwell optimality.}
  }%


\begin{document}
\maketitle
\section{Introduction}
Markov decision process (MDP)~\citep{puterman2014markov} is one of the most popular frameworks to model sequential decision-making by a single agent, with recent applications to game solving and reinforcement learning~\citep{sutton2018reinforcement,mnih2013playing}, healthcare decision-making~\citep{bennett2013artificial,steimle2017markov} and finance~\citep{bauerle2011markov}. Robust Markov decision processes (RMDPs)~\citep{iyengar2005robust,wiesemann2013robust} are a generalization of MDPs, dating as far back as \cite{satia1973markov} and extensively studied more recently after the seminal papers of \cite{iyengar2005robust} and \cite{nilim2005robust}. Robust MDPs consider a single agent repeatedly interacting with an environment with unknown instantaneous rewards and/or transition probabilities. The unknown parameters are assumed to be chosen adversarially from an {\em uncertainty set}, modeling the set of all plausible values of the parameters. 

Most of the robust MDP literature has focused on computational and algorithmic considerations for the {\em discounted return}, where the future instantaneous rewards are discounted with a discount factor $\gamma \in (0,1)$. In this case, RMDPs become tractable under certain {\em rectangularity} assumptions, with models such as sa-rectangularity~\citep{iyengar2005robust,nilim2005robust}, s-rectangularity~\citep{le2007robust,wiesemann2013robust}, k-rectangularity~\citep{mannor2016robust}, and r-rectangularity~\citep{goh2018data,goyal2022robust}. For rectangular RMDPs, iterative algorithms can be efficiently implemented for various distance- and divergence-based uncertainty sets~\citep{givan1997bounded,iyengar2005robust, ho2021partial,grand2021first,behzadian2021fast}, as well as gradient-based algorithms~\citep{wang2022convergence,kumar2023policy,li2022first,li2023policy}. 

Robust MDPs have found some applications in healthcare~\citep{zhang2017robust,goh2018data,grand2023robustness}, where the set of states represents the potential health conditions of the patients, actions represent the medical interventions, and it is critical to account for the potential errors in the estimated transition probabilities representing the evolution of the patient's health, and in inverse reinforcement learning~\citep{viano2021robust,chae2022robust} for imitations that are robust to shifts between the learner's and the experts' dynamics.
In some applications, introducing a discount factor can be seen as a modeling choice, e.g. in finance~\citep{deng2016deep}. In other applications like game-solving~\cite{mnih2013playing,brockman2016openai}, the discount factor is merely introduced for algorithmic purposes, typically to ensure the convergence of iterative algorithms or to reduce variance~\citep{baxter2001infinite}, and it may not have any natural interpretation. Additionally, large discount factors may slow the convergence rate of the algorithms.
 {\em Average optimality}, i.e., optimizing the limit of the average of the instantaneous rewards obtained over an infinite horizon, and {\em Blackwell optimality}, i.e., finding policies that remain discount optimal for all discount factors sufficiently close to $1$, provide useful objective criteria for environments with no natural discount factor, unknown discount factor, or large discount factor. To the best of our knowledge, the literature on robust MDPs that address these optimality notions is scarce. The authors in \cite{tewari2007bounded,wang2023robust} study RMDPs with the average return criterion but only focus on computing the optimal worst-case average returns among {\em stationary} policies without any guarantee that average optimal policies may be chosen in this class of policies. Similarly, some previous papers show the existence of Blackwell optimal policies for sa-rectangular RMDPs, but they require some restrictive assumptions, such as polyhedral uncertainty sets~\cite{tewari2007bounded,goyal2022robust} or unichain RMDPs with a unique average optimal policy~\citep{wang2023robust}. 
 
\vspace{2mm}
\noindent
{\bf Our goal} is here to study RMDPs with average and Blackwell optimality in the general case. Our {\bf main contributions} can be summarized as follows.
\begin{enumerate}
    \item {\bf Average optimality.} For sa-rectangular robust MDPs with compact uncertainty sets, we show the existence of a stationary deterministic policy that is average optimal, thus closing an important gap in the literature. Additionally, we describe several strong duality results and we highlight that the worst-case transition probabilities need not exist, a fact that has been overlooked by previous work. We also discuss the case of robustness against history-dependent worst-case transition probabilities.  \tb{In addition, we show that for s-rectangular RMDPs, surprisingly, average-optimal policies may not exist, and when they exist, stationary policies may be suboptimal. That is, history-dependent (Markovian) policies may be necessary to achieve the optimal average return. }
    \item {\bf Blackwell optimality.} We provide an extensive treatment of Blackwell optimality for sa-rectangular RMDPs. In this case, we provide a counterexample where Blackwell optimal policies fail to exist, and we identify the pathological oscillating behaviors of the robust value functions as the key factor for this non-existence. We show however that policies that remain $\epsilon$-Blackwell optimal for all values of $\epsilon>0$ always exist, \tb{and, in fact, that this property entirely characterizes average optimal policies. This result provides an interesting link between the three objective criteria studied in this paper (discount optimality, Blackwell optimality and average optimality)}. Finally, we introduce the notion of {\em definable} uncertainty sets, which are sets built on simple enough functions to ensure that the discounted value functions are ``well-behaved''. The definability assumption captures virtually all uncertainty sets used in practice for sa-rectangular RMDPs \tb{and it can be seen as a methodological choice of the decision-maker}. We show that for sa-rectangular RMDPs with definable uncertainty sets, a stationary deterministic policy always exists that is Blackwell optimal. We also highlight the connections between Blackwell optimal policies and average optimal policies. For s-rectangular RMDPs, we provide a simple example where Blackwell optimal policies do not exist, even in the simple case of polyhedral uncertainty. 
    \item {\bf Algorithms and numerical experiments.} We introduce three algorithms that converge to the optimal average return, in the case of definable sa-rectangular uncertainty sets. 
    We note that efficiently computing an average optimal policy for sa-rectangular RMDPs is a long-standing open problem in algorithmic game theory, and we present numerical experiments on three RMDP instances. 
    \item \tb{{\bf Methodological takeaways.} 
    Overall, our paper uncovers critical distinctions between sa-rectangular and s-rectangular frameworks. For sa-rectangular models, most of the insights from the widely studied discounted case carry over: the optimal average return and policies align with their discounted return counterparts as discount factors increase. Additionally, uncertainty sets based on common distances (e.g., Kullback-Leibler, $\ell_{p}$, and Wasserstein) yield Blackwell optimal policies and stable discounted value functions. In contrast, s-rectangular models lack guarantees for the existence of average optimal policies, which, if they exist, can even require complex, non-stationary structures in simple cases. This contrast positions sa-rectangular uncertainty as highly advantageous for practical RMDPs, ensuring implementable (stationary deterministic) optimal policies and simplifying policy interpretation.
}

\end{enumerate}

\vspace{2mm}
\noindent 
Table~\ref{tab:results-summary} summarizes our main results and compares them with prior results for robust MDPs and nominal MDPs. Our results are highlighted in {\bf bold}.

\begin{table}[htb]
    \centering
        \caption{Properties of optimal policies for different objective criteria. The set of states and the set of actions are finite. Our results are in bold. }

    \begin{tabular}{|m{4cm} |m{3.3cm}|m{3.5cm}| m{5.7cm}|} 
     \hline
    Uncertainty set $\U$ & Discount optimality & Average optimality & Blackwell optimality \\ [0.5ex] 
     \hline \hline
    Singleton (MDPs) & stationary,
    
    deterministic & stationary,
    
    deterministic & stationary,
    deterministic \\ 
     \hline
     sa-rectangular, compact & stationary,
     
     deterministic &  {\bf stationary,}
    
    {\bf deterministic}  &
    $\;$
     \begin{itemize}
         \item {\bf may not exist}
         \item {\bf  \boldmath$\exists \; \pi$ stationary deterministic,} 

        {\bf \boldmath$\pi \; \epsilon$-Blackwell optimal, \boldmath$\forall \epsilon>0$ }

         \item {\bf \boldmath$\pi$ also average optimal}
     \end{itemize}
     \\
      \hline
     sa-rectangular, compact,
     
     definable & stationary, 
     
     deterministic & {\bf stationary,}
    
    {\bf deterministic}  & 
    $\;$
     \begin{itemize}
         \item {\bf stationary, deterministic}
         \item {\bf \boldmath$\pi$ also average optimal}
     \end{itemize}
       \\
      \hline
     s-rectangular,
     
     compact convex & stationary,
     
     randomized &   
         $\;$
    \tb{\begin{itemize}
         \item {\bf may not exist}
         \item {\bf may be history-dependent, randomized}
     \end{itemize}}
     
       & {\bf may not exist} \\
      \hline
    \end{tabular}
    \label{tab:results-summary}
\end{table}  
\paragraph{Outline of the paper.} The rest of the paper is organized as follows. 
We introduce robust MDPs in Section \ref{sec:RMDPs} and provide our literature review there. We study average optimality in Section \ref{sec:average-optimality} and Blackwell optimality in Section \ref{sec:blackwell-optimality}. We study algorithms for computing average optimal policies for sa-rectangular RMDPs in Section \ref{sec:avg-alg}. \tb{Our main take-aways are summarized in Section \ref{sec:discussion}.}
\subsection{Notations.}
Given a finite set $\Omega$, we write $\Delta(\Omega)$ for the simplex of probability distributions over $\Omega$, and $P(\Omega)$ for the power set of $\Omega$. We allow ourselves to conflate vectors in  $\R^{|\Omega|}$ and functions on $\R^{\Omega}$.
For any two vectors $\bm{u},\bm{v} \in \R^{\Omega}$, the inequality $\bm{u} \leq \bm{v}$ is understood componentwise: $u_{s} \leq v_{s}, \forall \; s \in \Omega$, and for $\epsilon \geq 0$, we overload the notation $\bm{u} \leq \bm{v} + \epsilon$ to signify $u_{s} \leq v_{s} + \epsilon, \forall \; s \in \Omega$.
\section{Robust Markov decision processes}\label{sec:RMDPs}
In this section, we introduce discounted robust MDPs, and discuss the connections with stochastic games which prove useful in the rest of the paper.
\subsection{Discount optimality}\label{sec:RMDPs-discounted}
A robust Markov decision process is defined by a tuple $\left(\X,\A,\bm{r},\U,\bm{p}_{0}\right)$ where $\X$ is the finite set of states, $\A$ is the finite set of actions available to the decision maker (DM),  and $\bm{r} \in \R^{\X \times \A \times \X}$ is the instantaneous reward function. The set $\U \subset \left(\Delta(\X)\right)^{\X \times \A}$, usually called the {\em uncertainty set}, models the set of all possible transition probabilities $\bm{P} = \left(\bm{p}_{sa}\right)_{sa} \in \Delta(\X)^{\X \times \A}$,  and $\bm{p}_{0} \in \Delta(\X)$ is the initial distribution over the set of states $\X$. We emphasize that we assume throughout $\X$ and $\A$ to be finite. We also assume that $\U$ is a non-empty compact set. This is well-motivated from a practical standpoint, where $\U$ is typically built from some data based on some statistical distances, see the end of this section for some examples.

We write $\Pi_{\sf H}$ for the set of history-dependent, randomized policies that possibly depend on the entire past history $(s_{0},a_{0},...,s_{t})$ to choose\tb{, possibly randomly,} an action at period $t$. We write $\Pi_{\sf M}$ for the set of \emph{Markovian} policies, i.e., for policies $\pi \in \Pi_{\sf H}$ that only depend on $t$ and the current state $s_{t}$ for the choice of $a_t$. The set $\Pi_{\sf S}\subset \Pi_M$ of \emph{stationary} policies consists of the Markovian policies that do not depend on time, and $\Pi_{\sf SD}\subset \Pi_{\sf S}$ are those policies that are furthermore deterministic~\citep{puterman2014markov}. 
Given a policy $\pi \in \Pi_{\sf H}$ and transition probabilities $\bm{P}\in \U$, we denote by $\E_{\pi,\bm{P}}$ the expectation with respect to the distribution of the sequence $(s_t,a_t)_{t\geq 0}$ of states and actions induced by $\pi$ and $\bm{P}$ (as a function of the initial state).
Given $\gamma <1$, the value function induced by  $(\pi,\bm{P}) \in \Pi_{\sf H} \times \U$ is 
\[v^{\pi,\bm{P}}_{\gamma,s} = \E_{\pi,\bm{P}} \left[ \sum_{t=0}^{+\infty} \gamma^{t}r_{s_{t}a_{t}s_{t+1}} \; | \; s_{0} = s \right], \forall \; s \in \X,\]
and the discounted return is $R_{\gamma}(\pi,\bm{P}) = \bm{p}_{0}\tr\bm{v}^{\pi,\bm{P}}_{\gamma}$. 
A \emph{discount optimal policy} is an optimal solution to the optimization problem:
\begin{equation}\label{eq:robust-mdp-discounted}
    \sup_{\pi \in \Pi_{\sf H}} \min_{\bm{P} \in \U} \, R_{\gamma}(\pi,\bm{P}).
\end{equation}
\vspace{2mm}

\noindent {\bf Rectangularity, Bellman operator and adversarial MDP.} The discounted robust MDP problem~\eqref{eq:robust-mdp-discounted} becomes tractable when the uncertainty set $\U$ satisfies the following {\em s-rectangularity} condition:
\begin{equation}\label{eq:s-rectangularity}
    \U = \times_{s \in \X} \, \U_{s}, \quad \U_{s} \subseteq \Delta(\X)^{\A}, \quad \forall \; s \in \X,
\end{equation}
with the interpretation that $\U_s$ is the set of probability transitions from state $s \in \X$, as a function of the action $a \in \A$.
The authors in \cite{wiesemann2013robust} show that if the uncertainty set $\U$ is s-rectangular and compact convex, a discount optimal policy $\pi\opt$ can be chosen stationary (but it may be randomized). Additionally, $\pi\opt$ can be computed by first solving the equation: $\bm{v}\opt=T(\bm{v}\opt)$  where $T:\R^{\X} \rightarrow \R^{\X}$ is the Bellman operator, defined as
\begin{equation}\label{eq:operator-bellman}
    T_{s}(\bm{v}) = \max_{\pi \in \Pi_{\sf S}} \min_{\bm{P}_{s} \in \U_{s}} \; \sum_{a \in \A} \pi_{sa} \bm{p}_{sa}\tr\left( \bm{r}_{sa} + \gamma \bm{v}\right), \forall \; \bm{v} \in \R^{\X},\forall \; s \in \X,
\end{equation}
then returning the policy $\pi\opt$ that attains the maximum on each component of $T(\bm{v}\opt)$.

Let us now fix a stationary policy $\pi \in \Pi_{\sf S}$, and a s-rectangular RMDP instance.  The problem of computing the worst-case return of $\pi$, defined as
\begin{equation}\label{eq:adversarial MDP}
    \min_{\bm{P} \in \U} R_{\gamma}(\pi,\bm{P})
\end{equation}
can be reformulated as an MDP instance, called the {\em adversarial MDP}. In the adversarial MDP, the state set is the finite set $\X$, the set of actions available at $s \in \X$ is the compact set $\U_{s}$, and the rewards and transitions are continuous, see~\citep{ho2021partial,goyal2022robust} for more details. In particular, the {\em robust value function} $\bm{v}^{\pi,\U}_{\gamma} \in \R^{\X}$ of $\pi \in \Pi_{\sf S}$, defined as
$v^{\pi,\U}_{\gamma,s}= \min _{\bm{P}\in \U} v^{\pi,\bm{P}}_{\gamma,s}, \forall \; s \in \X$, 
satisfies the fixed-point equation:
\[
  v^{\pi,\U}_{\gamma,s}
  \; =\;
  \min_{\bm{P}_{s} \in \U_{s}} \, \sum_{a \in \A} \pi_{sa} \bm{p}_{sa}\tr\left( \bm{r}_{sa} + \gamma \bm{v}^{\pi,\U}_{\gamma}\right),
  \quad\forall \; s \in \X.
\] 
Crucially, the worst-case transition probabilities for a given $\pi \in \Pi_{\sf S}$ can be chosen in the set of extreme points of $\U$.

Under the stronger {\em sa-rectangularity} assumption:
\begin{equation}\label{eq:sa-rectangularity}
  \U = \bigtimes_{(s,a) \in \X \times \A} \, \U_{sa}, \quad
  \U_{sa} \subseteq \Delta(\X), \quad
  \forall \; (s,a) \in \X \times \A,
\end{equation}
a discount optimal policy may be chosen stationary and deterministic when $\U$ is compact~\citep{iyengar2005robust,nilim2005robust}. 
Typical examples of sa-rectangular uncertainty sets are based on $\ell_{p}$-distance with $p \in \{1,2,+\infty\}$ , with $\U_{sa}$ defined as
$\U_{sa} = \{\bm{p} \in \Delta(\X) \; | \; \| \bm{p} - \bm{p}_{sa} \|_{p} \leq \alpha_{sa}\}$ or based on Kullback-Leibler divergence: $\U_{sa} = \{\bm{p} \in \Delta(\X) \; | \; \kl\left(\bm{p},\bm{p}_{sa}\right)\leq \alpha_{sa}\}$ for $(s,a) \in \X \times \A$ and some radius $\alpha_{sa} \geq 0$~\citep{iyengar2005robust,nilim2005robust,panaganti2022sample}.
Typical s-rectangular uncertainty sets are constructed analogously to sa-rectangular ambiguity sets~\citep{wiesemann2013robust,grand2021first}.
\subsection{Connections with stochastic games}\label{sec:RMDPs-SGs}
 Two-player zero-sum stochastic games (abbreviated SGs in the rest of this paper) were introduced in the seminal work of Shapley~\cite{shapley1953stochastic} and model repeated interactions between two agents with opposite interests. 
 In each period, the players choose each an action, which influences the evolution of the (publicly observed) state. Agents aim at optimizing their gain, which depends on states and actions. The literature on stochastic games is extensive, we refer the reader to \cite{neyman2003stochastic,mertens2015repeated} for classical textbooks and to \cite{laraki2015advances,renault2019tutorial} for recent short surveys.
 
Stochastic games and robust MDPs share many similarities despite historically distinct communities, motivations, terminology, and lines of research. In this paper, we leverage some existing results in SGs to prove results for RMDPs. This section briefly describes the similarities and differences between the two fields. 
\vspace{2mm}

\noindent 
{\bf Similarities between SGs and RMDPs.}
It has been noted in the RMDP literature that rectangular RMDPs as in~\eqref{eq:robust-mdp-discounted} can be reformulated as stochastic games. The first player in the game is the decision-maker in the RMDP, who chooses actions in $\A$, while the adversary chooses the transition probabilities $\bm{P}$ in $\U$. The uncertainty set thus corresponds to the action set of the adversary. This connection has been alluded to multiple times in the robust MDP literature, e.g., section~5 of \cite{iyengar2005robust}, the last paragraph of page 7 of \cite{xu2010distributionally}, the last paragraph of section 2 in \cite{wiesemann2013robust}, and the introduction of \cite{grand2024convex}. \tb{For the sake of completeness, we make these connections precise in Appendix~\ref{app:more-detailed-sgs-rmdps}.}

It is crucial to note that polyhedral uncertainty sets in s-rectangular RMDPs can be modeled by a finite action set for the second player in the associated SG. This is because the worst-case transition probabilities may be chosen in the set of extreme points of $\U$, and polyhedra have finitely many extreme points. With this in mind, extreme points of the s-rectangular set $\U$ correspond to stationary {\em deterministic} strategies for the second player of the SG.

Interestingly, sa-rectangular RMDPs can be reformulated as a special case of {\em perfect information} stochastic games (see section 5 in~\citep{iyengar2005robust}). Perfect information SGs are SGs with the property that each state is controlled (in terms of rewards and transitions) by only one player, see  Chapter 4 in \cite{neyman2003stochastic} for an introduction. In the perfect information SG reformulation of sa-rectangular RMDPs, the set of states is the union of $\X$ and of $\X\times \A$, the first player controls states $s \in \X$, and the second player chooses the transition probabilities $\bm{p}_{sa}$ in states $(s,a) \in \X \times \A$. Instantaneous rewards are only obtained at the states of the form $(s,a)$. For the sake of completeness, a detailed construction is provided in Appendix~\ref{app:more-detailed-sgs-rmdps-sarec}. 

\vspace{2mm}

\noindent 
{\bf A fundamental distinction between SGs and RMDPs.} 
We now describe a crucial distinction between RMDPs and SGs, which has received limited attention in the RMDP literature. 
In the classical stochastic game framework, the first player {\em and} the second player can choose history-dependent strategies. 
In contrast, in robust MDPs as in~\eqref{eq:robust-mdp-discounted}, it is common to assume that the decision-maker chooses a history-dependent policy ($\pi \in \PiH$), but that the adversary is restricted to {\em stationary} policies, i.e., to choose $\bm{P} \in \U$, instead of $\bm{P} \in \U_{\sf H}$ with $\U_{\sf H}$ is the set of (randomized) history-dependent policies for the adversary~\cite{wiesemann2013robust, goyal2022robust}.  This fundamental difference between SGs and RMDPs is one of the main reasons why existing results for SGs do not readily extend to results for RMDPs. The focus of this paper is on stationary adversaries, as classical for RMDPs, and for clarity we make it explicit in the statements of all our results.

\tb{
\begin{remark}[History-dependent adversaries]
Note that we overload notations for the definition of history-dependent adversaries, as compared to the definition of stationary adversaries. In particular, when we write ``$\bm{P} \in \U$", we have that $\bm{P}$ is a collection transition probabilities from $\X \times \A$ to $\X$, i.e. $\bm{P}$ belongs to $\Delta(\X)^{\X \times \A}$. In contrast, when we write $\bm{P} \in \U_{\sf H}$, we have that $\bm{P}$ is a map from the set of all finite histories $\cH$ to $\U$, with $\cH = \cup_{t \in \N} \cH_{t}$, with $\cH_{t}$ the set of all possible histories up to period $t \in \N$: \[\cH_{t} = \{ \left(s_{0},a_{0},\bm{P}_{0},s_{1},a_{1},\bm{P}_{1},...,s_{t-1},a_{t-1},\bm{P}_{t-1},s_{t}\right) \; | \;  \left(s_i\right)_{i \leq t} \in \X^{t},\left(a_i\right)_{i \leq t-1} \in \A^{t-1},\left(\bm{P}_i\right)_{i \leq t-1} \in \U^{t-1}\}.\]
In fact, the distinction between stationary adversaries and history-dependent adversaries is irrelevant in the discounted return case, as we show in Proposition~\ref{prop:discounted-return-min-inf} below.
This proposition follows from~\cite{shapley1953stochastic} and we provide a concise proof in Appendix \ref{app:history dependent equivalent stationary - discount}.
Equality~\eqref{eq:rmdp-discounted-pi-hr-2} shows that facing a stationary adversary or a history-dependent adversary is equivalent for RMDPs with discounted return (and a decision-maker that can choose history-dependent policies). This proposition answers some of the discussions about stationary vs. non-stationary adversaries from the seminal paper on discounted sa-rectangular RMDPs~\citep{iyengar2005robust}.
\end{remark}
}
\begin{proposition}\label{prop:discounted-return-min-inf}
Let $\U$ be a convex compact s-rectangular uncertainty set. Then
    \begin{equation}\label{eq:rmdp-discounted-pi-hr-2}
     \sup_{\pi \in \Pi_{\sf H}} \inf_{\bm{P} \in \U_{\sf H}} R_{\gamma}(\pi,\bm{P}) = \sup_{\pi \in \Pi_{\sf H}} \min_{\bm{P} \in \U} R_{\gamma}(\pi,\bm{P}).
 \end{equation}
\end{proposition}

\section{Robust MDPs with average optimality}\label{sec:average-optimality}
In this section, we study the average optimality criterion for rectangular RMDPs. We start by introducing average optimality, and we discuss the gaps in the existing literature in Section \ref{sec:RMDPs-average}. We then show in Section \ref{sec:sa-rec-avg} that stationary deterministic average optimal policies exist for sa-rectangular compact uncertainty sets. However and perhaps surprisingly, we show in Section \ref{sec:s-rec-avg} for s-rectangular RMDPs that  Markovian policies may strictly improve upon stationary policies. \tb{We note that the main results in this section are entirely based on a direct approach 
to average optimality, i.e. we obtain our results without studying the connections between average optimality and discount/Blackwell optimality. We will make these connections explicitly in Section \ref{sec:blackwell-optimality}.}
\subsection{Average optimality and gaps in the existing literature}\label{sec:RMDPs-average}
{\em Average optimality} is a fundamental objective criterion extensively studied in nominal MDPs and in reinforcement learning, see Chapters 8 and~9 in \cite{puterman2014markov} and the survey~\cite{dewanto2020average}. Average optimality alleviates the need for introducing a discount factor by directly optimizing a long-run average $\Ravg(\pi,\bm{P})$ of the instantaneous rewards received over time. Intuitively, we would want $\Ravg(\pi,\bm{P})$ to capture the limit behavior of the average payoff $\E_{\pi,\bm{P}} \left[ \frac{1}{T+1}  \sum_{t=0}^{T} r_{s_{t}a_{t}s_{t+1}} \; | \; s_{0} \sim \bm{p}_{0} \right]$ over the first $T$ periods, for $T \in \N$. A well-known issue however, see for instance Example 8.1.1 in \cite{puterman2014markov}, is that these average payoffs need not have a limit as $T \rightarrow + \infty$. 

 \begin{example}\label{ex:no-average-limit}
 Consider a MDP with a single state and two actions $a_{0}$ and $a_{1}$, with payoffs 0 and 1 respectively. In this setup, a (deterministic) Markovian strategy is a sequence $(a_t)_{t \in \N}$ in $\{0,1\}^{\N}$, and the average payoff over the first $T$ periods is the frequency of action $a_{1}$ in these periods. Therefore, if the sequence $(a_t)_{t \in \N}$ is such that these frequencies do not converge, the average payoff up to $T$ does not converge either as $T\to +\infty$. We provide a detailed example in Appendix \ref{app:detail-example-no-average-limit}.
\end{example}

In this paper, we focus on the following definition of the average return $\Ravg(\pi,\bm{P})$ of a pair $\left(\pi,\bm{P}\right) \in \Pi_{\sf H} \times \U$:
\begin{equation}\label{eq:definition-average-return-exp-limsup}
  \Ravg(\pi,\bm{P}) =  \E_{\pi,\bm{P}} \left[ \limsup_{T \rightarrow + \infty}\frac{1}{T+1} \sum_{t=0}^{T} r_{s_{t}a_{t}s_{t+1}} \; | \; s_{0} \sim \bm{p}_{0} \right].
\end{equation}
Other natural definitions of the average return are possible, such as using the $\liminf$ instead of the $\limsup$, or taking the expectation before the $\limsup$. We will show in the next section that our main theorems still hold for these other definitions (see Corollary \ref{cor:avg-reward-main-other-ravg}).

Given a {\em stationary} policy $\pi \in \Pi_{\sf S}$ and for $\bm{P} \in \U$, our definition of the average return coincides with the usual definition of average return:
\[\Ravg(\pi,\bm{P}) =  \lim_{T \rightarrow + \infty} \E_{\pi,\bm{P}} \left[ \frac{1}{T+1} \sum_{t=0}^{T} r_{s_{t}a_{t}s_{t+1}} \; | \; s_{0} \sim \bm{p}_{0} \right], \forall \; (\pi,\bm{P}) \in \Pi_{\sf S} \times \U.\]
\tb{
Since the discount factor $\gamma \in [0,1)$ captures the importance of future rewards in the discounted return, one can expect that the average return coincides with the discounted return as $\gamma$ approaches $1$. This is indeed the case, and the following equation shows the connection between average return and discounted return of a fixed pair of stationary policy and transition probabilities:
 \[\lim_{\gamma \rightarrow 1} (1-\gamma)R_{\gamma}(\pi,\bm{P}) = \Ravg(\pi,\bm{P}), \forall \;(\pi,\bm{P}) \in \Pi_{\sf S} \times \U.\]
In fact, for {\em finite MDPs}, it is well-known that the same connection can be made {\em at optimality}, i.e. that for a fixed $\bm{P} \in \U$, we have
 \[\lim_{\gamma \rightarrow 1} \max_{\pi \in \PiS} (1-\gamma)R_{\gamma}(\pi,\bm{P}) =  \max_{\pi \in \PiS} \lim_{\gamma \rightarrow 1} (1-\gamma)R_{\gamma}(\pi,\bm{P}) =\max_{\pi \in \PiS} \Ravg(\pi,\bm{P}).\]
We refer to Section 2.4 in \cite{feinberg2012handbook} for more details on this connection. We will provide an analogous connection between the optimal {\em worst-case} discounted return and the optimal {\em worst-case} average return in Section \ref{sec:sa-rec-relation-blackwell-opt-avg-opt}.
}

One of our main goals in this section is to study the robust MDP problem with the average return criterion:
\begin{equation}\label{eq:avg-return-rmdps-general-0}
    \sup_{\pi \in \Pi_{\sf H}} \inf_{\bm{P} \in \U} \Ravg(\pi,\bm{P})
\end{equation} 
A solution $\pi \in \Pi_{\sf H}$ to \eqref{eq:avg-return-rmdps-general-0}, {\em when it exists}, is called an {\em average optimal policy}. 
The definition of the optimization problem \eqref{eq:avg-return-rmdps-general-0} warrants two important comments.

{\bf Properties of average optimal policies.}
First, we stress that \eqref{eq:avg-return-rmdps-general-0} optimizes the worst-case average return over \emph{history-dependent} policies. In contrast, prior work on average optimal RMDPs~\citep{tewari2007bounded, le2007robust,wang2023robust} have solely considered 
the optimization problem $ \sup_{\pi \in \Pi_{\sf S}} \inf_{\bm{P} \in \U} \Ravg(\pi,\bm{P})$, thus restricting the decision maker to {\em stationary policies}. At this point, it is not clear if stationary policies are optimal in \eqref{eq:avg-return-rmdps-general-0}, and this issue has been overlooked in the existing literature. One of our main contributions is to answer this question by the positive for sa-rectangular RMDPs and by the negative for s-rectangular RMDPs.

\begin{figure}[htb]
\begin{center}
    \begin{subfigure}{0.45\textwidth}
\includegraphics[width=\linewidth]{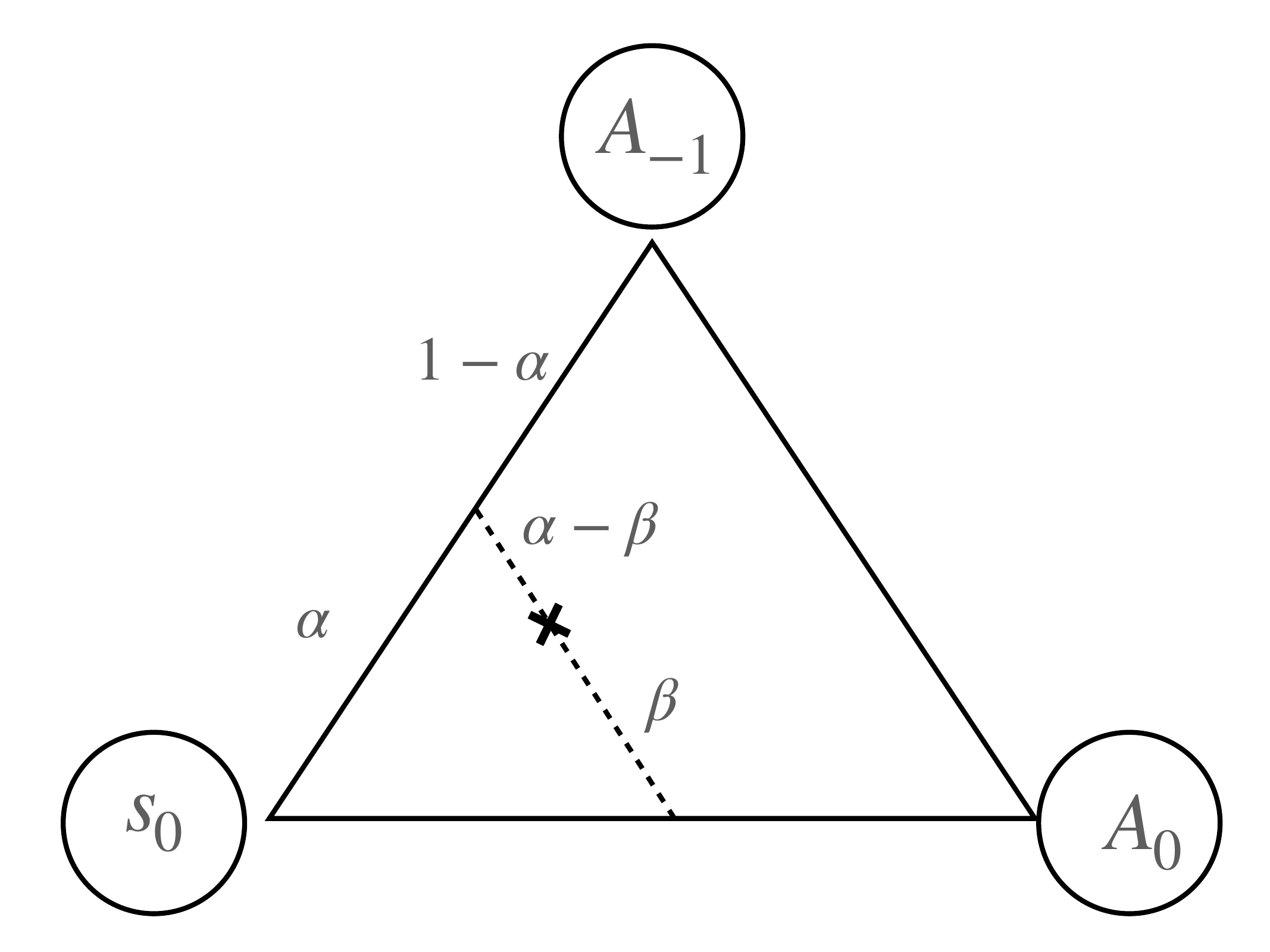}
    \caption{The point $(1-\alpha,\beta,\alpha-\beta)$ in $\Delta(\{s_{0},A_{-1},A_{0}\})$}
    \label{fig:simplex-alpha-beta}
    \end{subfigure}
    \begin{subfigure}{0.45\textwidth}
\includegraphics[width=\linewidth]{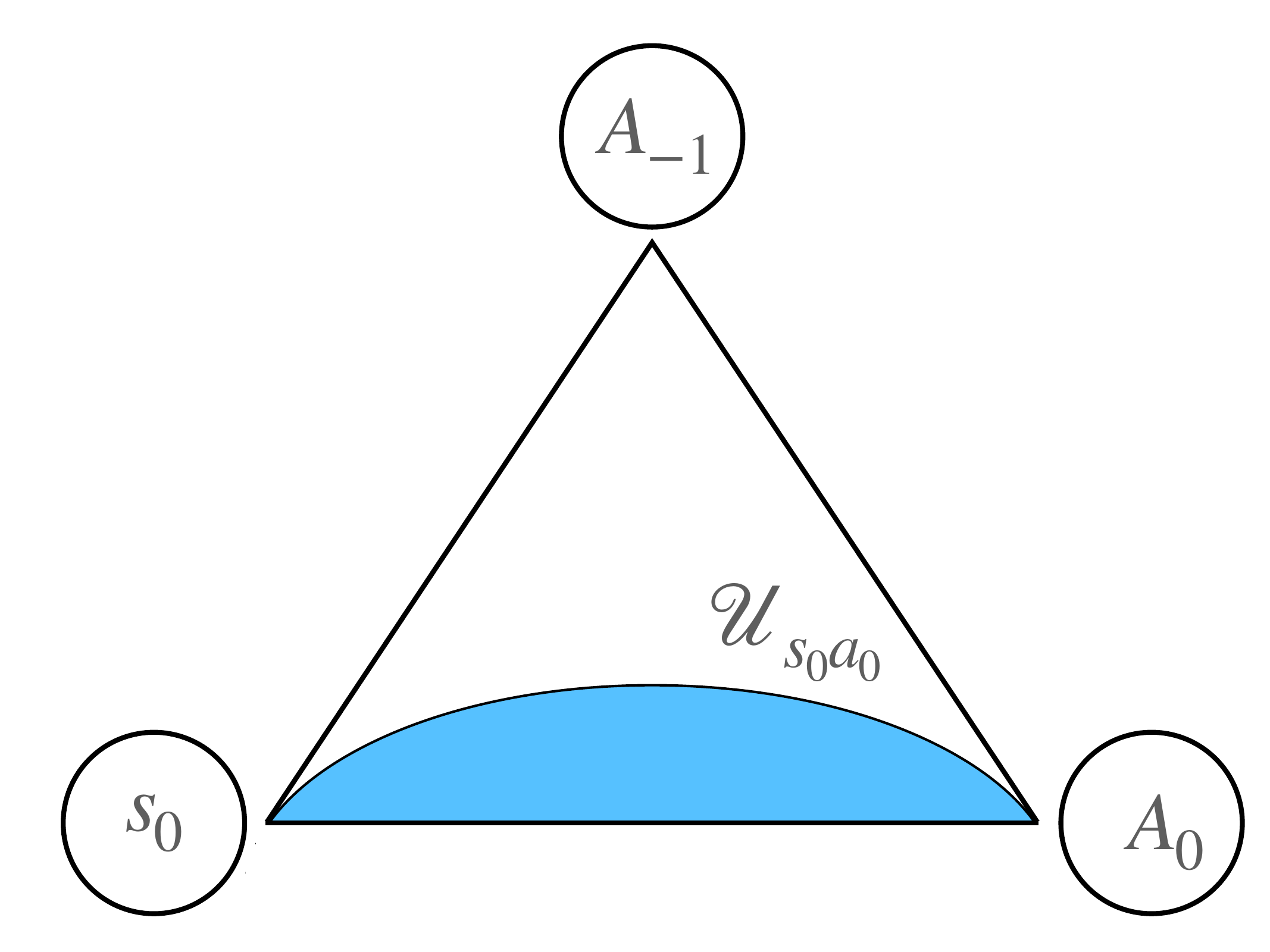}
    \caption{The uncertainty set $\U_{s_{0}a_{0}}$  (in blue)}
    \label{fig:D-D1}
    \end{subfigure}
    \end{center}
\caption{A simple robust MDP instance where $\inf_{\bm{P} \in \U} \Ravg(\pi,\bm{P})$ is not attained.} 
\label{fig:counter-example-inf-min}
\end{figure}

{\bf Non-existence of the worst-case transition probabilities.}
Second, we note that \eqref{eq:avg-return-rmdps-general-0} considers the {\em infimum} over $\bm{P} \in \U$ of the average return $\Ravg(\pi,\bm{P})$. The reason is that this infimum may not be attained. This issue has been neglected in prior work on robust MDPs with average return~\citep{tewari2007bounded,wang2023robust}, but recognized elsewhere repeatedly, see e.g. Section 1.4.4 in \cite{sorin2002first} or Example 4.10 in \cite{leizarowitz2003algorithm}. In the next proposition, we provide a simple counterexample.
\begin{proposition}\label{prop:inf-is-not-min-avg-rew}
There exists a robust MDP instance with an sa-rectangular compact convex uncertainty set, for which $\inf_{\bm{P} \in \U} \Ravg(\pi,\bm{P})$ is not attained for any $\pi \in \Pi_{\sf S}$.
\end{proposition}
The proof builds upon a relatively simple RMDP instance with three states and one action. We also reuse this RMDP instance for later important results in the next sections. For this reason, we describe it in detail below.
\proof{Proof of Proposition \ref{prop:inf-is-not-min-avg-rew}.}
We consider the robust MDP instance from Figure \ref{fig:counter-example-inf-min}.
There are three states $\{s_{0},A_{-1},A_{0}\}$, and only one action $a_0$, so that the RMDP reduces to a MDP with the adversary aiming at minimizing average returns. States  $A_{-1}$ and $A_{0}$ are  absorbing states with reward $-1$ and $0$ respectively. The initial state is $s_{0}$, and the instantaneous reward in $s_{0}$ is $0$. Elements of $\Delta(\X)$ are written $(p_0,p_1,p_2)$. The uncertainty set is
\begin{align*}
    \U_{s_{0}a_{0}} = \{ \left(1-\alpha,\beta,\alpha-\beta\right)\in \Delta(\X) \; | \; \alpha, \beta \in [0,1],\beta \leq \alpha(1-\alpha) \}, 
\end{align*}
see Figure \ref{fig:simplex-alpha-beta}. Given $\left(1-\alpha,\beta,\alpha-\beta\right)\in \U_{s_{0}a_{0}} $,  $\alpha$ is the probability to leave state $s_{0}$ and $\beta$ as the probability to go to state $A_{-1}$.  

The decision-maker has a single policy, which we denote $\pi= a_0$. We claim that $\inf_{\bm{P} \in \U_{s_{0}a_{0}}} \Ravg(\pi,\bm{P})= -1$, with 
$\Ravg(\pi,\bm{P})>-1$ for each $\bm{P} \in \U_{s_{0}a_{0}}$. Let $\bm{P}= (1-\alpha,\beta, \alpha-\beta)\in \U_{s_{0}a_{0}}$ be arbitrary. If $\alpha= 0$, the sequence of states does not leave $s_0$, and $\Ravg(\pi,\bm{P})=0$. If $\alpha>0$, the sequence of states eventually reaches $A_{-1}$ or $A_{0}$, and the  probability of reaching $A_{-1}$ (resp. $A_{0}$) is $\frac{\beta}{\alpha}$ (resp. $\frac{\alpha-\beta}{\alpha}$). Therefore, $\Ravg(\pi,\bm{P})=-\frac{\beta}{\alpha}$. Since $\alpha>0$, one has $\beta \leq \alpha(1-\alpha)< \alpha$, hence $\Ravg(\pi,\bm{P})>-1$. On the other hand, taking $\beta= \alpha (1-\alpha)$, one has $\Ravg(\pi,\bm{P})=-(1-\alpha)$ which can be arbitrarily close to $-1$, hence 
$\inf_{\bm{P} \in \U_{s_{0}a_{0}}} \Ravg(\pi,\bm{P})= -1$.
\hfill \halmos \endproof
\tb{The MDP instance from Figure \ref{fig:counter-example-inf-min} in the proof of Proposition \ref{prop:inf-is-not-min-avg-rew} highlights several properties that distinguish the average return from the discounted return:}
\begin{itemize}
    \item \tb{The discounted return $(\pi,\bm{P}) \mapsto R_{\gamma}(\pi,\bm{P})$ is always continuous in the entries of $\pi \in \PiS$ and $\bm{P} \in \U$, and this is why $\inf_{\bm{P} \in \U} R_{\gamma}(\pi,\bm{P})$ is attained when $\U$ is compact. In contrast, Figure \ref{fig:counter-example-inf-min} shows an instance where the average return is discontinuous in the input of the transition probabilities, and even for a sequence $\left(\bm{P}_{n}\right)_{n \in \N}$ of transition probabilities attaining $\inf_{\bm{P} \in \U} \Ravg(\pi,\bm{P})$ (at the limit $n \rightarrow +\infty$) we may have
    \[\lim_{n \rightarrow + \infty} \Ravg(\pi,\bm{P}_{n}) \neq \Ravg(\pi,\lim_{n \rightarrow + \infty}\bm{P}_{n}).\]
    This is the case for instance for 
    \[\bm{P}_{n} = (1-\alpha_n,\alpha_n(1-\alpha_n),\alpha_{n}^{2}),\alpha_n = 1/n, n \in \N,\]
    in which case $\lim_{n \rightarrow + \infty} \Ravg(\pi,\bm{P}_{n})=-1$ but $\lim_{n \rightarrow + \infty}\bm{P}_{n} = (1,0,0)$ and $\Ravg(\pi,\lim_{n \rightarrow + \infty}\bm{P}_{n}) = 0$. This shows that the limit point of a sequence $\left(\bm{P}_{n}\right)_{n \in \N}$ attaining the optimal average return (in the limit as $n \rightarrow + \infty$) may be vacuous and should be interpreted with caution.}
    \item \tb{Figure \ref{fig:counter-example-inf-min} also highlights the subtle interplay between discounted and average optimality. In particular, it is straightforward to verify that the worst-case transition probabilities for the discounted return is $\bm{P}\opt_{\gamma} = (1-\alpha\opt_{\gamma},\alpha\opt_{\gamma}(1-\alpha\opt_{\gamma}),\alpha^{\star 2}_{\gamma})$ with 
    $\alpha\opt(\gamma)= \frac{\sqrt{1-\gamma}-(1-\gamma)}{\gamma}$ for a given discount factor $\gamma \in [0,1)$ (see Appendix \ref{app:proof-th-sa-rec-no-Blackwell}). Since $\lim_{\gamma \rightarrow 1} \alpha\opt(\gamma)=0$, we have
    \[\lim_{\gamma \rightarrow 1} \bm{P}\opt_{\gamma} = (1,0,0).\]
    Similar to the previous point, this means that the limit of the sequence $\left(\bm{P}_{\gamma}\opt\right)_{\gamma \in [0,1)}$ of optimal worst-case transition probabilities for the discounted case may become vacuous and converge to some transition probabilities that achieve an average return very far from the worst-case average return.}
\end{itemize}
Prior works on RMDPs have focused on sa-rectangular RMDPs, either with polyhedral uncertainty set~\citep{le2007robust,tewari2007bounded} or under the assumption that the Markov chains induced by any pair of policy and transition probabilities are unichain~\citep{wang2023robust}. In these special cases, we prove in Appendix \ref{app:ravg-inf-equal-min} that $\inf_{\bm{P} \in \U} \Ravg(\pi,\bm{P}) = \min_{\bm{P} \in \U} \Ravg(\pi,\bm{P})$, i.e., worst-case transition probabilities indeed exist for any stationary policy $\pi \in \Pi_{\sf S}$. We note that this question is not discussed in prior work~\citep{tewari2007bounded,wang2023robust}. Since we don't make any assumption on the uncertainty set $\U$ beyond compactness, we will keep the notation $\inf_{\bm{P} \in \U} \Ravg(\pi,\bm{P})$ instead of $\min_{\bm{P} \in \U} \Ravg(\pi,\bm{P})$ in the rest of the paper.

We conclude this section with the following technical lemma, adapted from the literature on {\em nominal} MDP with a finite set of states and compact set of actions~\citep{bierth1987expected}. We briefly discuss this lemma in Appendix \ref{app:discussion-lemma-avg-rew-adversarial mdp}.
\begin{lemma}[Adapted from Theorem 2.5, \cite{bierth1987expected}]\label{lem:avg-rew-nom-eps-stationary}
    Let $\U$ be compact and s-rectangular. Let $\pi \in \Pi_{\sf S}$. Then 
    \[ \inf_{\bm{P} \in \U_{\sf H}} \Ravg(\pi,\bm{P}) = \inf_{\bm{P} \in \U} \Ravg(\pi,\bm{P}).\]
\end{lemma}
\subsection{The case of sa-rectangular robust MDPs}\label{sec:sa-rec-avg}
We first focus on the case of sa-rectangular RMDPs. 
Our main result in this section is that there always exist average optimal policies that are stationary and deterministic. 
\begin{theorem}\label{th:avg-reward-main}
    Consider an sa-rectangular robust MDP with a compact uncertainty set $\U$. There exists an average optimal policy that is stationary and deterministic:
    \[ 
    \sup_{\pi \in \Pi_{\sf H}} \inf_{\bm{P} \in \U} \Ravg(\pi,\bm{P}) \;=\; 
    \max_{\pi \in \Pi_{\sf SD}} \inf_{\bm{P} \in \U} \Ravg(\pi,\bm{P}).
    \]
\end{theorem}
The proof proceeds in two steps and leverages existing results from the literature on perfect information SGs~\cite{liggett1969stochastic}.
\proof{Proof of Theorem \ref{th:avg-reward-main}.}
In the first step of the proof, we show that Theorem \ref{th:avg-reward-main} is true in the special case where $\U$ is polyhedral. In this case, $\U$ only has a finite number of extreme points, and the robust MDP problem is equivalent to a perfect information stochastic game with finitely many actions for both players, as described in Section \ref{sec:RMDPs}. 
\tb{
From Theorem 1 in \cite{liggett1969stochastic}, one has
\[\inf_{\bm{P} \in \U} \max_{\pi \in \Pi_{\sf SD}}  \Ravg(\pi,\bm{P}) =  \max_{\pi \in \Pi_{\sf SD}} \inf_{\bm{P} \in \U} \Ravg(\pi,\bm{P}).\]
}

\tb{Now this implies that
\begin{equation}\label{eq:instead of gimbert-1a}
    \inf_{\bm{P} \in \U} \sup_{\pi \in \Pi_{\sf H}}  \Ravg(\pi,\bm{P}) = \inf_{\bm{P} \in \U} \max_{\pi \in \Pi_{\sf SD}}  \Ravg(\pi,\bm{P}) =  \max_{\pi \in \Pi_{\sf SD}} \inf_{\bm{P} \in \U} \Ravg(\pi,\bm{P}) \leq \sup_{\pi \in \Pi_{\sf H}} \inf_{\bm{P} \in \U} \Ravg(\pi,\bm{P})  
\end{equation}
where the first equality is from the inner supremum being an MDP since the adversary is stationary and the inequality is from $\Pi_{\sf SD} \subset \Pi_{\sf H}$.
}
From weak duality we always have $\sup_{\pi \in \Pi_{\sf H}} \inf_{\bm{P} \in \U} \Ravg(\pi,\bm{P}) \leq \inf_{\bm{P} \in \U} \sup_{\pi \in \Pi_{\sf H}}  \Ravg(\pi,\bm{P})$.
Therefore, \tb{all terms in \eqref{eq:instead of gimbert-1a} are equal, and we have}
\[\sup_{\pi \in \Pi_{\sf H}} \inf_{\bm{P} \in \U} \Ravg(\pi,\bm{P}) = \max_{\pi \in \Pi_{\sf SD}} \inf_{\bm{P} \in \U} \Ravg(\pi,\bm{P}).\]
In the second step of the proof, we show that Theorem \ref{th:avg-reward-main} holds for general compact uncertainty sets, without the assumption that $\U$ is polyhedral as in the first step of the proof. Let $\epsilon>0.$ For each $\pi \in \Pi_{\sf SD}$, consider $\bm{P}^{\pi}=(\bm{p}_{sa}^{\pi})_{s,a} \in \U$ such that
\[ \inf_{\bm{P} \in \U} \Ravg(\pi,\bm{P}) \leq \Ravg(\pi,\bm{P}^{\pi}) \leq \inf_{\bm{P} \in \U} \Ravg(\pi,\bm{P}) + \epsilon. \]
Note that  $\{ \bm{P}^{\pi} \; | \; \pi \in \Pi_{\sf SD}\}$ is a finite set since $\Pi_{\sf SD}$ is a finite set. \tb{Let us now consider the sa-rectangular uncertainty set $\U^{\sf f}= \times_{(s,a) \in \X \times \A} \U^{\sf f}_{sa}$ where $\U^{\sf f}_{sa} \subset \Delta(\X)$ is the convex hull of the finite set $\{ \bm{p} \in \Delta(\X) \; | \; \exists \; \pi \in \Pi_{\sf SD}, \bm{p} = \ \bm{p}_{sa}^{\pi}\}.$} By construction, $\U_{\sf f}$ is polyhedral and sa-rectangular. Additionally, $\U_{\sf f} \subset \U$, and for any $\pi \in \Pi_{\sf SD}$, we have
\begin{equation}\label{eq:ineq-Uf-U}
    \inf_{\bm{P} \in \U_{\sf f}} \Ravg(\pi,\bm{P}) \leq \Ravg(\pi,\bm{P}^{\pi}) \leq \inf_{\bm{P} \in \U} \Ravg(\pi,\bm{P}) + \epsilon. 
\end{equation}
Therefore, 
\begin{align*}
     \sup_{\pi \in \Pi_{\sf H}} \inf_{\bm{P} \in \U} \Ravg(\pi,\bm{P}) & \leq \sup_{\pi \in \Pi_{\sf H}} \inf_{\bm{P} \in \U_{\sf f}} \Ravg(\pi,\bm{P}) \\
     & = \max_{\pi \in \Pi_{\sf SD}} \inf_{\bm{P} \in \U_{\sf f}} \Ravg(\pi,\bm{P}) \\
     & \leq \max_{\pi \in \Pi_{\sf SD}} \inf_{\bm{P} \in \U} \Ravg(\pi,\bm{P}) + \epsilon.
\end{align*}
where the first inequality uses $\U_{\sf f} \subset \U$, the equality follows from the first step of the proof since $\U_{\sf f}$ is polyhedral, \tb{and the last equality follows from \eqref{eq:ineq-Uf-U}.}
Therefore, for all $\epsilon>0$, we have $\sup_{\pi \in \Pi_{\sf H}} \inf_{\bm{P} \in \U} \Ravg(\pi,\bm{P}) \leq \max_{\pi \in \Pi_{\sf SD}} \min_{\bm{P} \in \U} \Ravg(\pi,\bm{P}) + \epsilon$, and we can conclude that 
\[\sup_{\pi \in \Pi_{\sf H}} \inf_{\bm{P} \in \U} \Ravg(\pi,\bm{P}) \leq \max_{\pi \in \Pi_{\sf SD}} \inf_{\bm{P} \in \U} \Ravg(\pi,\bm{P})\]
which  implies 
\[\sup_{\pi \in \Pi_{\sf H}} \inf_{\bm{P} \in \U} \Ravg(\pi,\bm{P}) =  \max_{\pi \in \Pi_{\sf SD}} \inf_{\bm{P} \in \U} \Ravg(\pi,\bm{P}).\]
\hfill \halmos \endproof

\vspace{2mm}
\noindent
Theorem \ref{th:avg-reward-main} has several noteworthy consequences. 

First, to the best of our knowledge, we are the first to study average return RMDPs in all generality without constraining the problem to stationary policies and to show that average optimal policies for sa-rectangular RMDPs may be chosen stationary and deterministic. In this respect, Theorem \ref{th:avg-reward-main} closes a significant gap that has been overlooked in the existing literature on RMDPs with average return~\cite{tewari2007bounded,wang2023robust}. 

Another consequence of Theorem \ref{th:avg-reward-main} is the following strong duality statement. We provide the detailed proof in Appendix \ref{app:proof sa-rec strong duality}. 
\begin{theorem}\label{th:sa-rec-strong-duality}
    Consider an sa-rectangular robust MDP with a compact uncertainty set $\U$. Then the following strong duality results hold:
    \begin{align}
        \sup_{\pi \in \Pi_{\sf H}} \inf_{\bm{P} \in \U} \Ravg(\pi,\bm{P}) & =  \inf_{\bm{P} \in \U} \sup_{\pi \in \Pi_{\sf H}} \Ravg(\pi,\bm{P}) \label{eq:strong-duality-sup-inf}\\
        \max_{\pi \in \Pi_{\sf SD}} \inf_{\bm{P} \in \U} \Ravg(\pi,\bm{P}) & = \inf_{\bm{P} \in \U} \max_{\pi \in \Pi_{\sf SD}} \Ravg(\pi,\bm{P}) \label{eq:strong-duality-max-min}
    \end{align}
\end{theorem}
Theorem \ref{th:sa-rec-strong-duality} is akin to the strong duality results for discounted RMDPs~\citep{wiesemann2013robust,goyal2022robust}. It states that strong duality still holds for sa-rectangular RMDPs with average optimality. We note that in Equality \eqref{eq:strong-duality-max-min}, the maximum over $\pi \in \Pi_{\sf SD}$ is always attained (it is the supremum over a finite set) while the infimum may not be attained,  as we illustrate in Proposition~\ref{prop:inf-is-not-min-avg-rew}. 

Strong duality is crucial to study the case of {\em history-dependent adversaries}, as we now show. The case of a {\em non-stationary} adversary has also gathered interest and is discussed in \cite{iyengar2005robust,nilim2005robust}. Interestingly, we show below that this latter model is equivalent to the case of a stationary adversary. The proof is  concise and relies on Theorem \ref{th:avg-reward-main}, on Theorem \ref{th:sa-rec-strong-duality} and on the properties of MDPs with compact action sets. 
\begin{theorem}\label{th:avg-rew-history-dep-adversaries}
    Consider an sa-rectangular robust MDP with a compact uncertainty set $\U$. Then
    \[ 
    \sup_{\pi \in \Pi_{\sf H}} \inf_{\bm{P} \in \U_{\sf H}} \Ravg(\pi,\bm{P}) \;=\; 
    \sup_{\pi \in \Pi_{\sf H}} \inf_{\bm{P} \in \U} \Ravg(\pi,\bm{P}).
    \]\end{theorem}

Note that the same results also hold for the case of Markovian adversaries.
\proof{Proof of Theorem \ref{th:avg-rew-history-dep-adversaries}.}
We have
\begin{align*}
   \inf_{\bm{P} \in \U_{\sf H}} \sup_{\pi \in \Pi_{\sf H}}  \Ravg(\pi,\bm{P}) & \leq \inf_{\bm{P} \in \U} \sup_{\pi \in \Pi_{\sf H}}  \Ravg(\pi,\bm{P}) \\
   & =  \sup_{\pi \in \Pi_{\sf H}} \inf_{\bm{P} \in \U} \Ravg(\pi,\bm{P}) \\
   & = \max_{\pi \in \Pi_{\sf SD}} \inf_{\bm{P} \in \U} \Ravg(\pi,\bm{P}) \\
   & = \max_{\pi \in \Pi_{\sf SD}} \inf_{\bm{P} \in \U_{\sf H}} \Ravg(\pi,\bm{P})  \\
    & \leq \sup_{\pi \in \Pi_{\sf H}} \inf_{\bm{P} \in \U_{\sf H}} \Ravg(\pi,\bm{P})
\end{align*}
where the first inequality uses $\U \subset \U_{\sf H}$,  the three equalities follow respectively from  Theorem \ref{th:sa-rec-strong-duality}, Theorem \ref{th:avg-reward-main} and   Lemma \ref{lem:avg-rew-nom-eps-stationary}, and the last inequality from  the inclusion $\Pi_{\sf SD} \subset \Pi_{\sf H}$.
\hfill \halmos \endproof

\vspace{2mm}
\noindent
{\bf Other natural definitions of average optimality.}
Alternative definitions of the average return induced by  $(\pi,\bm{P}) \in \Pi_{\sf H} \times \U$ exist,  such as any of the following:
\begin{align}
    & \liminf_{T \rightarrow + \infty} \E_{\pi,\bm{P}} \left[ \frac{1}{T+1}  \sum_{t=0}^{T} r_{s_{t}a_{t}s_{t+1}} \; | \; s_{0} \sim \bm{p}_{0} \right], \label{eq:liminf-exp} \\
    & \limsup_{T \rightarrow + \infty} \E_{\pi,\bm{P}} \left[ \frac{1}{T+1}  \sum_{t=0}^{T} r_{s_{t}a_{t}s_{t+1}} \; | \; s_{0} \sim \bm{p}_{0} \right], \label{eq:limsup-exp} \\
    &   \E_{\pi,\bm{P}} \left[ \liminf_{T \rightarrow + \infty}\frac{1}{T+1} \sum_{t=0}^{T} r_{s_{t}a_{t}s_{t+1}} \; | \; s_{0} \sim \bm{p}_{0} \right]. \label{eq:exp-liminf}
\end{align}
 The following inequality holds   for every $ (\pi,\bm{P}) \in \Pi_{\sf H} \times \U$, (see Lemma 2.1 in \cite{feinberg2012handbook}):
  \begin{equation}\label{eq:ordering-ravg}
      \eqref{eq:exp-liminf} \leq~\eqref{eq:liminf-exp} \leq~\eqref{eq:limsup-exp} \leq \Ravg(\pi,\bm{P}),
  \end{equation}
with equality whenever  $(\pi,\bm{P}) \in \Pi_{\sf S} \times \U$.

At this point, the reader may wonder if the results that we proved in this section also hold for these other definitions of the average return and if these returns each need to be analyzed separately. Fortunately, based on Inequality~\eqref{eq:ordering-ravg} and the fact that all the definitions of the average return coincides over $\Pi_{\sf S} \times \U$, we can show in a straightforward manner that Theorem \ref{th:avg-reward-main}, Theorem~\ref{th:sa-rec-strong-duality} and Theorem \ref{th:avg-rew-history-dep-adversaries} still hold for the average return as in $\eqref{eq:exp-liminf},\eqref{eq:liminf-exp}$ and $\eqref{eq:limsup-exp}$. This shows that all these definitions of the average return are equivalent in the sense that there exists a stationary deterministic policy that is average optimal simultaneously for all these objective functions. In particular, we obtain the following important corollary. For the sake of conciseness, we provide the proof in Appendix \ref{app:proof-corollaries}. 
\begin{corollary}\label{cor:avg-reward-main-other-ravg}
    Consider an sa-rectangular robust MDP with a compact uncertainty set $\U$. Let $\hRavg$ \tb{be} as defined in \eqref{eq:liminf-exp}, \eqref{eq:limsup-exp}, or \eqref{eq:exp-liminf}.
    \begin{enumerate}
        \item[{\bf P0}.] There exists an average optimal policy that is stationary and deterministic, and which coincides with an average optimal policy for $\Ravg$ as in \eqref{eq:definition-average-return-exp-limsup}:
    \[ 
    \sup_{\pi \in \Pi_{\sf H}} \inf_{\bm{P} \in \U} \hRavg(\pi,\bm{P}) \;=\; 
    \max_{\pi \in \Pi_{\sf SD}} \inf_{\bm{P} \in \U} \hRavg(\pi,\bm{P}) = \max_{\pi \in \Pi_{\sf SD}} \inf_{\bm{P} \in \U} \Ravg(\pi,\bm{P}).
    \] 
    \item[{\bf P1}.] The strong duality results from Theorem \ref{th:sa-rec-strong-duality} still \tb{holds} when  replacing $\Ravg$ by $\hRavg$.
    \item[{\bf P2}.] The equivalence between stationary adversaries and history-dependent adversaries from Theorem \ref{th:avg-rew-history-dep-adversaries} still hold when replacing $\Ravg$ by $\hRavg$.
    \end{enumerate}
\end{corollary}
\tb{Overall, our results in this section highlight the advantages of sa-rectangular models for the practical implementations of RMDP models for real-world applications: stationary and deterministic policies are easier to interpret and deploy than history-dependent or randomized policies, optimal policies may consider stationary or non-stationary adversaries alike, and any reasonable definitions of the average returns yield the same average optimal policies. We will show that the situation is very different for s-rectangular models in the next section.} 
\subsection{The case of s-rectangular robust MDPs}\label{sec:s-rec-avg}
\tb{We now study s-rectangular robust MDPs with the average return criteria. We show two surprising results: average optimal policies may not exist, even for polyhedral uncertainty, and if average optimal policies exist they may need to be history-dependent (Markovian). These two results are in contrast with the good properties of {\em sa-rectangular} uncertainty sets when using the average return, as highlighted in the previous section. 
}

\tb{
  {\bf An example with no average optimal policy.}
We first show a simple RMDP instance with s-rectangular uncertainty, which lacks an optimal policy. Our proof is based on the RMDP instance represented in Figure \ref{fig:smex}, with threes states $s_{0},A_{0},A_{1}$ and two actions $T$ and $B$. Each arrow corresponds to a possible transition, and is labelled with its probability, and with the payoff obtained along this transition. The decision maker starts in $s_{0}$ while the states $A_{0}$ and $A_{1}$ are absorbing with rewards $0$ and $1$ respectively. The adversary chooses the scalar $p \in [0,1]$, corresponding to some probability transitions $\bm{P} \in \U$ as represented in Figure \ref{fig:smex}. In particular, we have the following proposition.}
\tb{
\begin{proposition}\label{prop:avg-s-rec no existence opt policies}
For the polyhedral s-rectangular robust MDP instance represented in Figure \ref{fig:smex}, the following statements hold: 
\begin{itemize}
    \item[{\bf P0}.]  $\sup_{\pi\in \Pi_{\sf H}} \inf_{\bm{P} \in \U} \Ravg(\pi,\bm{P}) = 1.$
    \item[{\bf P1}.] Average optimal policies do not exist: $\forall \; \pi \in \PiH, \inf_{\bm{P} \in \U} \Ravg(\pi,\bm{P}) < 1$. 
\end{itemize}
\end{proposition}
\proof{Proof.}
{\bf P0.} Note that all rewards are in $\{0,1\}$, and therefore $\sup_{\pi\in \Pi_{\sf H}}\inf_{\bm{P} \in \U} \Ravg(\pi,\bm{P}) \leq 1$. Now for $\pi$ the stationary policy $\pi$ that assigns the probability $\ep >0$ to the action $T$ we have $R_{\mathrm{avg}}(\pi,\bm{P})\geq 1-\ep$ for each $\bm{P} \in \U$. This shows that $\sup_{\pi\in \Pi_{\sf S}}\inf_{\bm{P} \in \U} \Ravg(\pi,\bm{P}) = 1$ so that $\sup_{\pi\in \Pi_{\sf H}}\inf_{\bm{P} \in \U} \Ravg(\pi,\bm{P}) = 1$.}
  \tb{  \begin{figure}
\begin{center}
    \begin{subfigure}{0.4\textwidth}
\includegraphics[width=\linewidth]{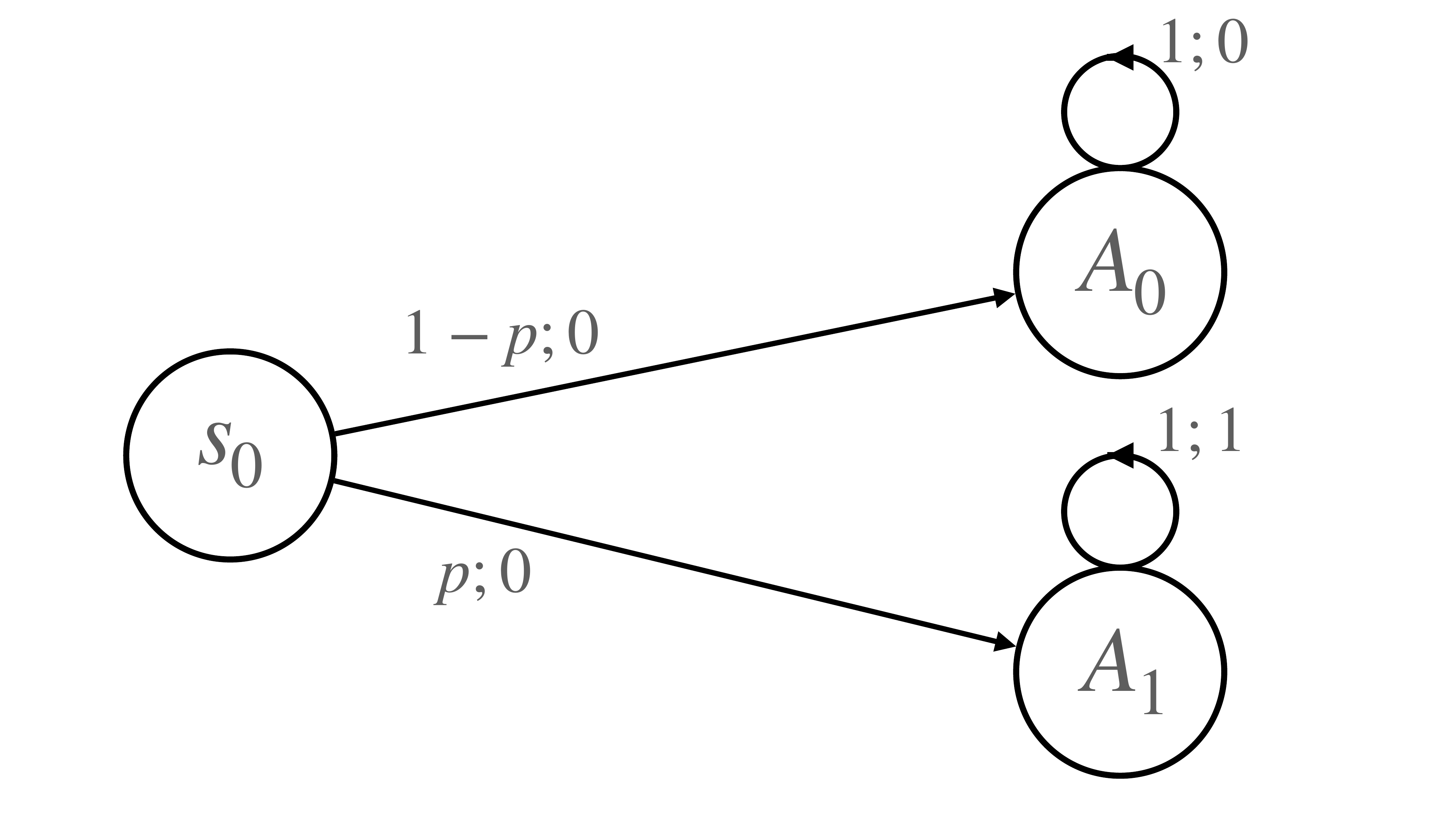}
    \caption{\tb{Transition for action $T$}}
    \label{fig:smex_T}
    \end{subfigure}
    \begin{subfigure}{0.4\textwidth}
\includegraphics[width=\linewidth]{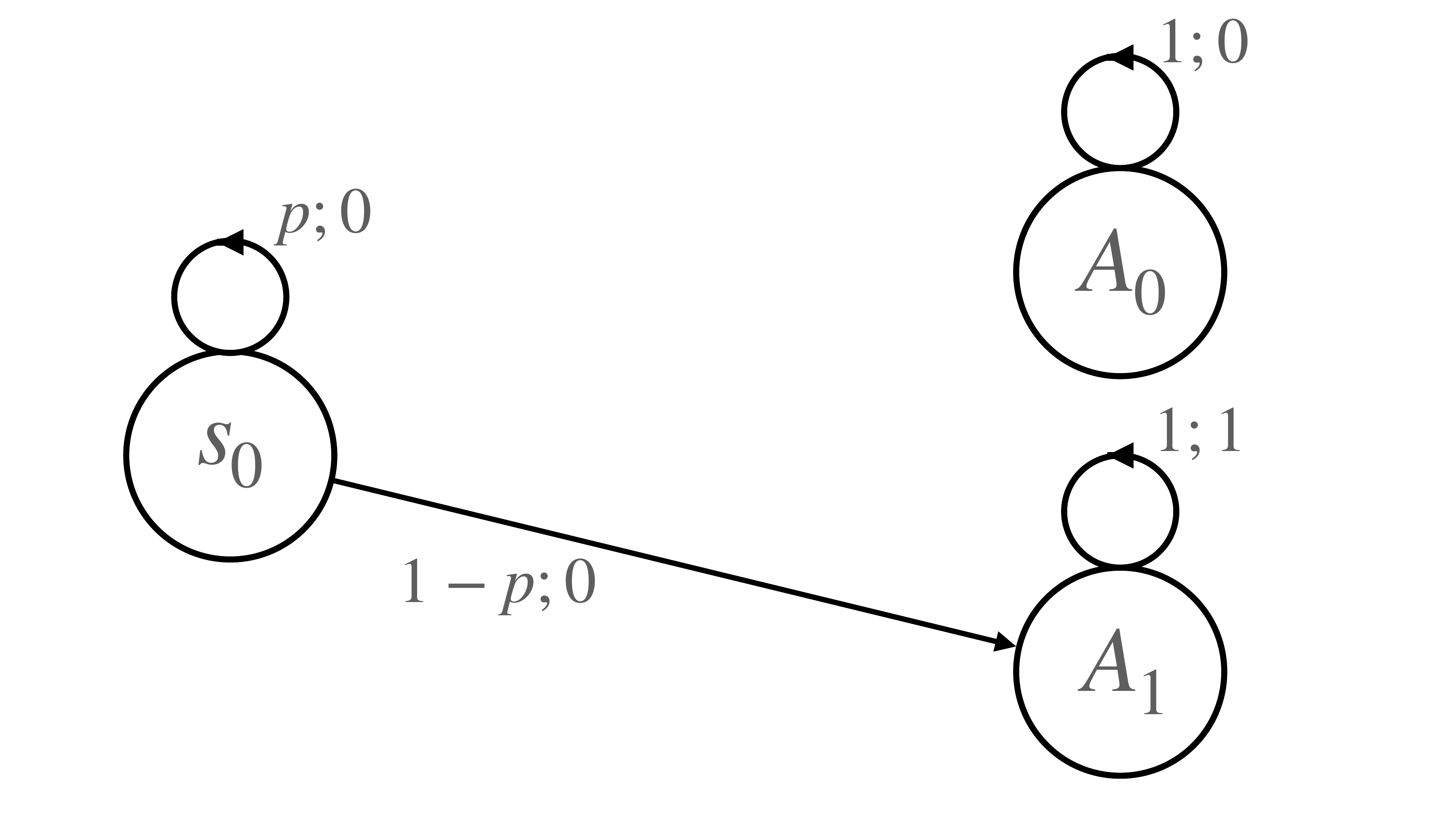}
    \caption{\tb{Transition for action $B$}}
    \label{fig:smex_B}
    \end{subfigure}
    \end{center}
\caption{\tb{Transitions and rewards for the MDP instance for Proposition . The adversary chooses $p \in [0,1]$ and the decision-maker chooses an action in $\{T,B\}$}.}\label{fig:smex}
\end{figure}
}

\tb{{\bf P1.} We now show that for any $\pi \in \PiH$, we have $\inf_{\bm{P} \in \U} \Ravg(\pi,\bm{P})<1$. Let $\pi$ be a history-dependent policy. Note that the decisions chosen become irrelevant once the system leaves state $s_{0}$ and reaches an absorbing state. Therefore, we can consider only histories up to period $t$ where the decision-maker remains in state $s_{0}$. In this case, a history consists of the sequence of actions chosen by the decision-maker until period $t-1$. However, if the system is still in state $s_{0}$ at period $t$ (before the decision-maker chooses their action), then the decision-maker must have chosen action $B$ in all the previous periods. Therefore, a history up to period $t$ actually consists in only knowing the period $t$, i.e., we can identify history-dependent policies with Markovian policies. Therefore, a history-dependent $\pi\in \Pi_{\sf H}$ can be identified with the sequence of action probabilities $(x_t)$, where $x_t \in [0,1]$ is the probability of choosing $T$ in period $t$, conditional on the current state being the initial state $s_{0}$. There are two cases. If $x_t=0$ for each $t \in \N$, then $R_{\mathrm{avg}}(\pi,\bm{P})= 0$, when $\bm{P}$ corresponds to $p=1$. If $x_t>0$ for some $t \in \N$, then $R_{\mathrm{avg}}(\pi,\bm{P})<1$ for any $\bm{P}$ such that $p \in (0,1)$. Either way, one has  $\inf_{\bm{P} \in \U} \Ravg(\pi,\bm{P})<1$.
}
\hfill \halmos \endproof
\tb{The non-existence (in general) of average optimal policies for s-rectangular is in contrast with the case of sa-rectangular uncertainty. This result is quite surprising, since for the discounted return sa-rectangular and s-rectangular RMDPs share similar properties, e.g. both admit optimal policies that can be chosen stationary. We also note that while average optimal policies may not exist in the RMDP from Figure \ref{fig:smex}, the value of the optimal average return is $1$ and {\em stationary} policies can achieve $1-\epsilon>0$ for any $\epsilon \in (0,1)$. The next example shows that this is not always the case.}

\begin{figure}
\begin{center}
    \begin{subfigure}{0.4\textwidth}
\includegraphics[width=\linewidth]{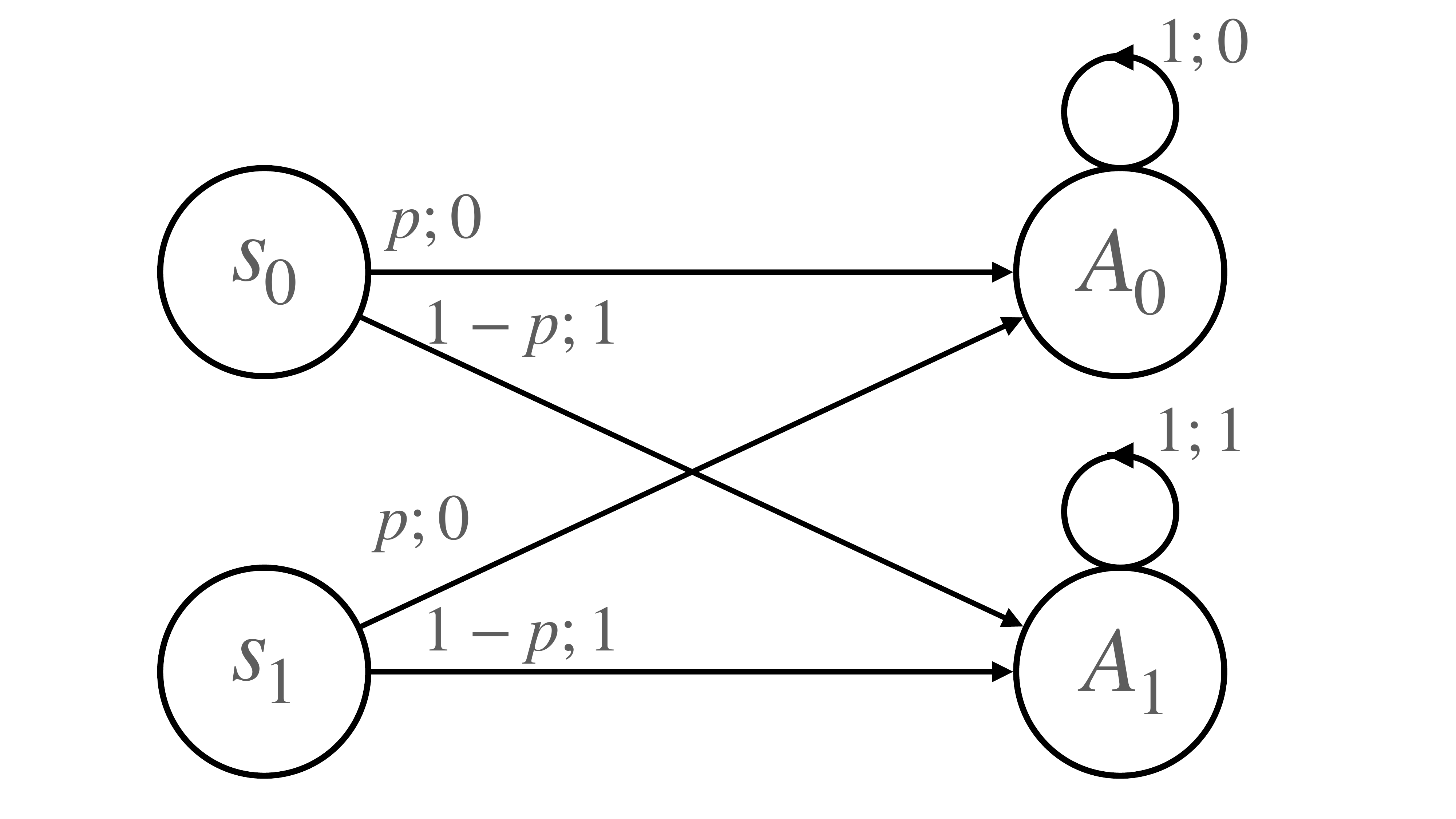}
    \caption{Transition for action $T$}
    \label{fig:tbm_T}
    \end{subfigure}
    \begin{subfigure}{0.4\textwidth}
\includegraphics[width=\linewidth]{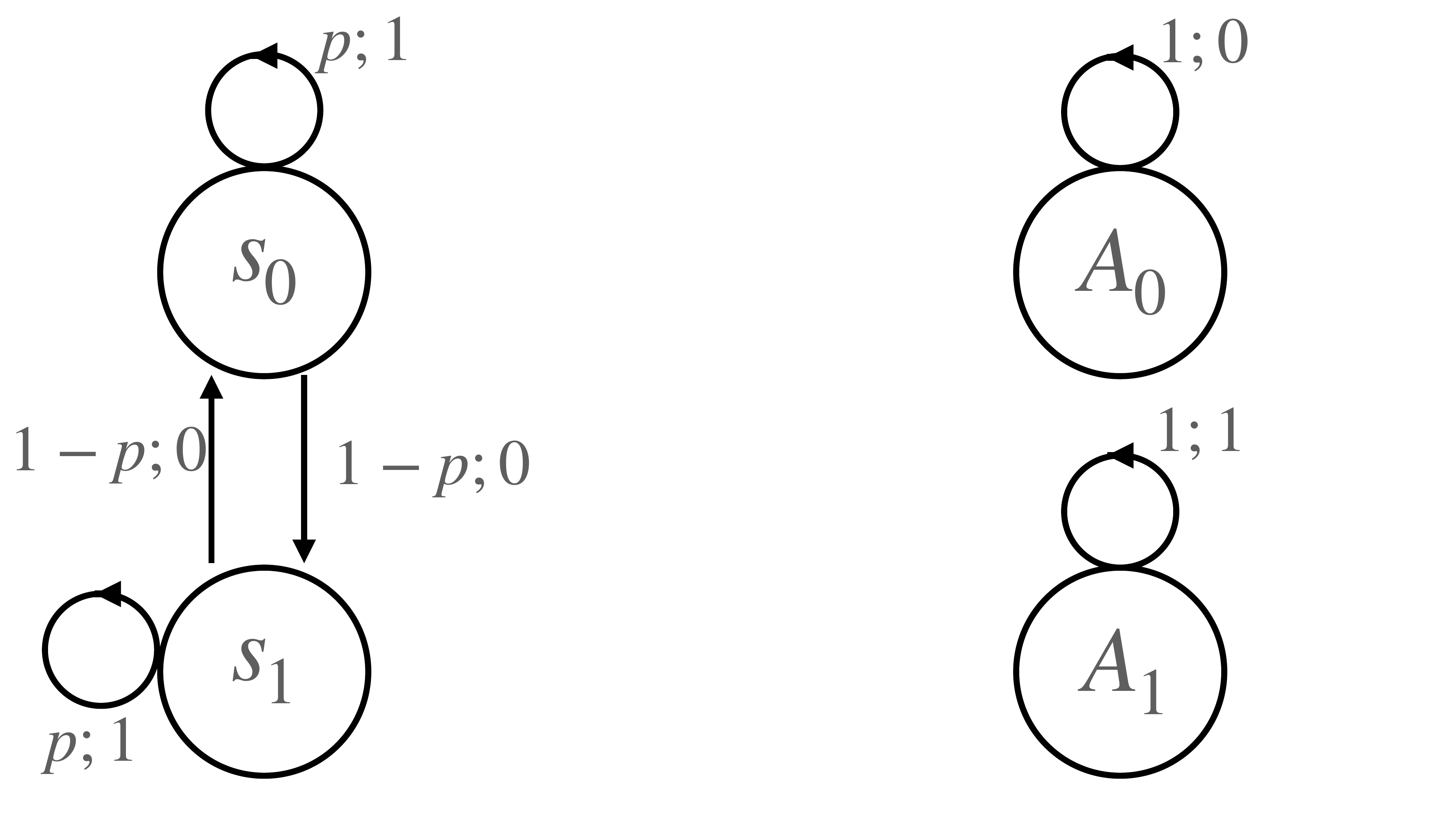}
    \caption{Transition for action $B$}
    \label{fig:tbm_B}
    \end{subfigure}
    \end{center}
\caption{Transitions and rewards for The Big Match~\cite{blackwell1968big,gillette1957stochastic} reformulated as an s-rectangular RMDP. The adversary chooses $p \in [0,1]$ and the decision-maker chooses an action in $\{T,B\}$.}\label{fig:TBM-s-rec}
\end{figure}

\tb{
{\bf An example where average optimal policies are Markovian.}
We now show that for s-rectangular RMDPs, when average optimal policies exist, they may need to be {\em history-dependent} (Markovian), while stationary policies perform very poorly. Again, this is in contrast with the case of sa-rectangular models where average optimal policies can be chosen stationary and deterministic (see our Theorem \ref{th:avg-reward-main}).
}
We consider the RMDP represented in Figure \ref{fig:TBM-s-rec}, which is adapted from a stochastic game known as \tbm~\citep{blackwell1968big}.  There are two absorbing states $A_{0}$, and two non-absorbing states,  $s_0$ and $s_1$. The action set of the decision-maker is $\{T,B\}$, and the initial state is $s_0$. Once the absorbing state $A_k$ ($k \in \{0,1\}$) is reached (if ever), the game remains there,  and all payoffs are equal to $k$. 
For $s=s_0,s_1$, the set $\U_s\subset \Delta(\X)^\A$ of adversary choices in state $s$ is obtained as follows.
 The adversary chooses a scalar $p \in [0,1]$. If action $T$ is chosen, the game moves either to $A_0$ or $A_1$, with probabilities $p$ and $1-p$, and thus ends immediately; if action $B$ is chosen, the game moves either to $s_0$ or $s_1$, with probabilities $1-p$ and $p$. Note that $\U_s$  and $r_{sas'}$ are independent of $s\in \{s_0,s_1\}$, so that $s_0$ and $s_1$ are duplicates. 
 One could alternatively replace them with a single state $s$, if payoffs $r_{sas'}$ were allowed to depend on the decisions of the adversaries.

The intuition behind this game is the following. The higher the value of $p \in [0,1]$, the higher the probability of ending up in the absorbing state $A_0$ if the decision maker chooses to end the game (action $T$), but the worse the current payoff if the decision maker chooses to continue the game (action $B$). If the adversary could perfectly anticipate the decision maker's choices, they would choose $p=1$ when anticipating $T$, and $p=0$ when anticipating $B$.

\begin{proposition}\label{prop:tbm-avg-rew}
For the polyhedral s-rectangular robust MDP instance represented in Figure \ref{fig:TBM-s-rec}, the following statements hold:
\begin{description}
\item[{\bf P0}.]  \tb{There exists $\bm{P} \in \U$ such that $\Ravg(\pi,\bm{P}) \leq 1/2$  for any history-dependent policy $\pi \in \Pi_{\sf H}$.}
\item[{\bf P1}.] $\max_{\pi \in \Pi_{\sf S}}\inf_{\bm{P} \in \U} \Ravg(\pi,\bm{P}) = 0$ while $\min_{\bm{P} \in \U}\max_{\pi \in \Pi_{\sf S}} \Ravg(\pi,\bm{P}) = \frac12$.
\item[{\bf P2}.] $\max_{\pi \in \Pi_{\sf M}}\inf_{\bm{P} \in \U} \Ravg(\pi,\bm{P}) = \min_{\bm{P} \in \U}\max_{\pi \in \Pi_{\sf M} }\Ravg(\pi,\bm{P}) = \frac12$.
\end{description}
\end{proposition}

Thus,  there is a duality gap when the decision-maker is restricted to stationary policies (\textbf{P1}) but not
when the decision-maker can pick Markovian strategies (\textbf{P2}).
\medskip

\proof{Proof.}
Claims \textbf{P0} and \textbf{P1} are known in the SG literature, but not \textbf{P2}. In particular, \textbf{P0} follows from Blackwell's original argument (See Theorem 1 in \cite{blackwell1968big}), and \textbf{P1} is Lemma 1 in Chapter 12 of \cite{neyman2003stochastic}. We provide proofs for completeness. 
    \begin{description}
        \item[{\bf P0}.] Assume that adversary chooses $p=\frac12$, and let $\pi\in \Pi_{\sf H}$ be arbitrary. In the event where the DM always chooses $B$, rewards in the different stages are iid, equal to 0 or 1 with equal probabilities. In the event where the DM eventually chooses $T$, the game ends in $A_0$ and $A_1$ with equal probabilities. This implies that $\Ravg(\pi,\bm{P})=1/2$.
        \item[{\bf P1}.] We claim that $\max_{\pi \in \Pi_{\sf S} }\Ravg(\pi,\bm{P}) \geq \frac12$  for each $\bm{P}\in \U$. Together with \textbf{P0}, this will imply the second equality of \textbf{P1}. Let be given $p\in [0,1]$. If $p\leq \frac12$, the stationary strategy $\pi= T$ ensures that the game ends immediately, and that the probability of ending in $A_1$ is $1-p$, so that $\Ravg(\pi,\bm{P})=1-p\geq \frac12$. If $p>\frac12$, the strategy $\pi=B$ ensures that the game never ends, with $\Ravg(\pi,\bm{P})=p>\frac12$. This proves the claim. For the first equality, let $\pi\in \Pi_{\sf S}$ be given. There are two cases. If under $\pi$ the probability of choosing $T$ is 0, the game never ends. With $p= 1$, one has $\Ravg(\pi,\bm{P})=0$. If instead the probability of choosing $T$ is positive (and independent of the period of decisions since $\pi\in \Pi_{\sf S}$), then with $p=0$,  the game ends up in state $A_0$ with probability $1$, so that $\Ravg(\pi,\bm{P})=0$. This proves the first equality.
        \item[{\bf P2}.] Consider the Markovian policy $\pi_*$ that chooses both actions $T$ and $B$ with equal probabilities in the first period $t=0$, then always chooses $B$, and let $p\in [0,1]$ be arbitrary. If the DM chooses $T$ in the first period, the game ends in $A_0$ or $A_1$, with probabilities $p$ and $1-p$ respectively; if instead the DM chooses $B$ in the first period, average payoffs converge to $1-p$. Hence $\Ravg(\pi_*,\bm{P})= \frac{1}{2}\times p + \frac{1}{2} \times (1-p) = \frac{1}{2}$. \tb{This implies that 
         $\sup_{\pi \in \PiM} \inf_{\bm{P} \in \U} \Ravg(\pi,\bm{P}) \geq 1/2$. On the other hand, one has  
  \[\sup_{\pi \in \PiM} \inf_{\bm{P} \in \U} \Ravg(\pi,\bm{P}) \leq 
  \inf_{\bm{P} \in \U} \sup_{\pi \in \PiH} \Ravg(\pi,\bm{P})\leq \frac12,\]
  where the first inequality follows from weak duality and $\PiM\subset \PiH$, and the second one follows from \textbf{P0}. Therefore, one has
  \[\sup_{\pi \in \Pi_{\sf M}}\inf_{\bm{P} \in \U} \Ravg(\pi,\bm{P}) = \inf_{\bm{P} \in \U}\sup_{\pi \in \Pi_{\sf M} }\Ravg(\pi,\bm{P}) = \frac12.\]}
    \end{description}
    \tb{The fact that the $\sup$ is reached on the left-hand side follows from the construction of $\pi_*$. On the right-hand side, the fact that the $\inf$ is reached follows from \textbf{P0}, and the fact that the $\sup$ is reached for each $\bm{P}$ follows from the existence of optimal stationary strategies in nominal MDPs. This proves \textbf{P2}.}
\hfill \halmos \endproof

Overall, Proposition \ref{prop:avg-s-rec no existence opt policies} and
Proposition \ref{prop:tbm-avg-rew} show the contrast between sa-rectangular RMDPs, where average optimal policies always exist and may be chosen stationary and deterministic (Theorem \ref{th:avg-reward-main}), and s-rectangular RMDPs, where average optimal policies may not exist, may have to be Markovian and randomized, and where stationary policies may not even be $\epsilon$-average optimal. \tb{These examples also show that discount optimal policies may perform very poorly for the average return. Indeed, it is well-known that for the discounted return an optimal policy in \tbm{} can be chosen {\em stationary} and chooses action $T$ with probability $1-\gamma$ and action $B$ with probability $\gamma$.  From Proposition \ref{prop:tbm-avg-rew}, any stationary policy will achieve a worst-case average return of $0$, so that the discount optimal policy also achieves a worst-case average return of $0$, even for very large discount factors.}

We note that the case of a history-dependent adversary is the case studied in the stochastic game literature. In particular, \cite{blackwell1968big} show for \tbm{} that 
\[ \sup_{\pi \in \Pi_{\sf H}} \inf_{\bm{P} \in \U_{\sf H}} \Ravg(\pi,\bm{P}) =   \inf_{\bm{P} \in \U_{\sf H}} \sup_{\pi \in \Pi_{\sf H}} \Ravg(\pi,\bm{P}) = 1/2.\]
In addition, Markovian policies fail to achieve a worst-case average return higher than $0$ against a history-dependent adversary: $\inf_{\bm{P} \in \U_{\sf H}} \Ravg(\pi,\bm{P}) = 0$ for $\pi \in \Pi_{\sf M}$. We refer to Chapter 12 in \cite{neyman2003stochastic} for a modern exposition.
\tb{
\begin{remark}
    The reader may be surprised by the contrast between the properties of average optimal policies for sa-rectangular and for s-rectangular uncertainty sets. In particular, a recent paper~\citep{ramani2024family} shows an intriguing connection between these two models of uncertainty for the {\em discounted} return, and proves that discounted value functions for some s-rectangular uncertainty sets may coincide with the discounted value functions of some sa-rectangular uncertainty sets. We note that the results in \cite{ramani2024family} solely focus on the numerical values of the discounted value functions for different uncertainty sets, and in particular \cite{ramani2024family} do not touch upon the properties of optimal policies, which may still differ for s-rectangular and sa-rectangular models, even in the discounted case, where discount optimal policies can always be chosen deterministic for sa-rectangular models while they may have to be randomized for s-rectangular models.
\end{remark}}
\section{Robust MDPs with Blackwell optimality}\label{sec:blackwell-optimality}
We now study Blackwell optimality, which provides an adequate optimality criterion when there is no natural notion of discounting. We first introduce Blackwell optimality in Section \ref{sec:rmdp-blackwell}, then focus on sa-rectangular RMDPs in Section \ref{sec:sa-rec-rmdp-blackwell} and on s-rectangular RMDPs in Section \ref{sec:s-rec-rmdp-blackwell}.
\subsection{Blackwell optimality and $\epsilon$-Blackwell optimality}\label{sec:rmdp-blackwell}
An important limitation of the average return is that it ignores any reward obtained in finite time, which may be problematic. For instance, when optimizing patient trajectories over time, practitioners are concerned with the long-term goals (typically, survival at discharge) but also with the given patient condition at any point in time.
A policy is \emph{Blackwell optimal} if it is discount optimal for all discount factors sufficiently close to $1$~\citep{puterman2014markov}, thus balancing long-term and short-term goals.
\begin{definition}\label{def:blackwell-optimality}
\tb{Let $\U$ be compact.}
    A policy $\pi \in \Pi_{\sf H}$ is Blackwell optimal if there exists $\gamma_{0} \in (0,1)$ such that
\begin{equation}
    \inf_{\bm{P} \in \U} R_{\gamma}(\pi,\bm{P}) \geq \sup_{\tb{\pi'} \in \Pi_{\sf H}} \inf_{\bm{P} \in \U} R_{\gamma}(\pi',\bm{P}), \quad
    \forall \; \gamma \in (\gamma_{0},1).
\end{equation}
\end{definition}
For nominal MDPs, i.e. when $\U$ is a singleton and the sets $\X,\A$ are finite, there always exists a stationary Blackwell optimal policy~\citep{puterman2014markov}. We are interested in the existence (or not) of Blackwell optimal policies for rectangular robust MDPs. 
\tb{We also introduce the following definition of $\epsilon$-Blackwell optimality, where a policy remains approximately discount optimal for all discount factors sufficiently large.}
\begin{definition}[$\epsilon$-Blackwell optimality.]\label{def:blackwell-opt-eps}
\tb{Let $\U$ be compact.}
  A policy $\pi \in \Pi_{\sf H}$ is $\epsilon$-Blackwell optimal for $\epsilon > 0$ if there exists a discount factor $\gamma_{\epsilon} \in (0,1)$ such that
    \begin{equation}
      \inf_{\bm{P} \in \U} \, (1-\gamma)R_{\gamma}(\pi,\bm{P})
      \; \geq\; 
      \sup_{\tb{\pi'} \in \Pi_{\sf H}} \inf_{\bm{P} \in \U} \, (1-\gamma)R_{\gamma}(\pi',\bm{P}) - \epsilon,
      \quad
      \forall \; \gamma \in (\gamma_{\epsilon},1).
    \end{equation}
\end{definition}
\tb{
Compared to Definition \ref{def:blackwell-optimality}, note that Definition \ref{def:blackwell-opt-eps} involves normalized discounted returns $\gamma \mapsto (1-\gamma)R_{\gamma}(\pi,\bm{P})$, since the discounted return $\gamma \mapsto R_{\gamma}(\pi,\bm{P})$ may diverge to $+\infty$ or $-\infty$. We also note that since the number of states and the number of actions are finite, the normalized discounted return always belongs to the interval $[\min_{(s,a,s')} r_{sas'}, \max_{(s,a,s')} r_{sas'}]$. Therefore, for $\epsilon = \max_{(s,a,s')} |r_{sas'}|$, every policy $\pi \in \PiH$ is $\epsilon$-Blackwell optimal, and the notion of $\epsilon$-Blackwell optimality is only useful when $\epsilon$ is close to $0$, in a similar way that the usual definition of $\epsilon$-optimal solution to an optimization problem is only useful when $\epsilon$ is close to $0$. In particular, in this section we will study the existence and the properties of policies that remain $\epsilon$-Blackwell optimal {\em for all values of $\epsilon>0$.}
}
To the best of our knowledge, we are the first to introduce and study  $\epsilon$-Blackwell optimality in the context of robust MDPs. In the rest of this section, we repeatedly use the following result pertaining to the existence of $\epsilon$-Blackwell optimal policy in the adversarial MDP for any choice of $\epsilon>0$, which we adapt from the literature on MDPs with compact action sets.
\begin{theorem}[Corollary 5.26, \cite{sorin2002first}]\label{th:eps-blackwell-compact-mdp}
Let $\U$ be a compact s-rectangular uncertainty set. Let $\pi \in \Pi_{\sf S}$ and $\epsilon>0$. Then there exist $\gamma_{\epsilon} \in (0,1)$ and $\bm{P}_{\epsilon} \in \U$ such that
\[\min_{\bm{P} \in \U} (1-\gamma)R_{\gamma}(\pi,\bm{P}) \leq (1-\gamma) R_{\gamma}(\pi,\bm{P}_{\epsilon}) \leq \min_{\bm{P} \in \U} (1-\gamma)R_{\gamma}(\pi,\bm{P}) + \epsilon, \forall \; \gamma \in (\gamma_{\epsilon},1).\]
\end{theorem}
\subsection{The case of sa-rectangular robust MDPs}\label{sec:sa-rec-rmdp-blackwell}
In this section, we provide a complete analysis of Blackwell optimality for sa-rectangular RMDPs. In Section \ref{sec:counter-example}, we first show that surprisingly, Blackwell optimal policies may fail to exist, \tb{although there exists a stationary and deterministic policy that is $\epsilon$-Blackwell optimal {\em for all values of $\epsilon>0$}}. Additionally, a policy with this property is average optimal, as we show in Section \ref{sec:sa-rec-relation-blackwell-opt-avg-opt}. Finally, in Section \ref{sec:definable-sa-rec} we introduce the notion of {\em definable uncertainty sets}, a very general class of uncertainty for which Blackwell optimal policies exist.
\subsubsection{Existence and non-existence results}\label{sec:counter-example}
We contrast here the existence properties of Blackwell optimal policies and $\epsilon$-Blackwell optimal policies (as introduced in Definition \ref{def:blackwell-opt-eps}). We first show that some sa-rectangular robust MDPs do not have  Blackwell optimal policies.
\begin{theorem}\label{th:sa-rec-no-blackwell}
There exists an sa-rectangular robust MDP instance, with a compact convex uncertainty set $\U$, and with no Blackwell optimal policy:
    \[\forall \; \pi \in \Pi_{\sf H}, \forall \; \gamma \in (0,1), \exists \; \gamma' \in (\gamma,1), \min_{\bm{P} \in \U} R_{\gamma'}(\pi,\bm{P}) < \sup_{\pi' \in \Pi_{\sf S}} \min_{\bm{P} \in \U} R_{\gamma'}(\pi',\bm{P}). \]
\end{theorem}
We view this new Theorem \ref{th:sa-rec-no-blackwell} as surprising because the existence of Blackwell optimal policies has been established in various related frameworks, including nominal MDPs~\citep{puterman2014markov}, or sa-rectangular RMDPs under some additional assumptions~\citep{goyal2022robust,wang2023robust}. For the sake of conciseness, we defer the proof of Theorem \ref{th:sa-rec-no-blackwell} to Appendix \ref{app:proof-th-sa-rec-no-Blackwell}, and we only provide here some intuition for this negative result.
 
 \noindent
 The proof of Theorem \ref{th:sa-rec-no-blackwell} is based on a variant of the simple example used for Proposition \ref{prop:inf-is-not-min-avg-rew}. We will consider the same instance, except that the DM has two actions $a_{1}$ and $a_{2}$ at the non-absorbing state $s_{0}$, instead of one. We identify actions $a_{1}, a_{2}$  with deterministic stationary policies. The crucial property is that the robust value functions $\gamma \mapsto v^{a_{1},\U}_{s_{0},\gamma}$ and $\gamma \mapsto v^{a_{2},\U}_{s_{0},\gamma}$ have an oscillatory behavior when $\gamma$ approaches $1$. As these robust value functions oscillate more and more often as $\gamma \rightarrow 1$, they intersect infinitely often on any interval close to $1$, so that there are no discount factors close enough to $1$ after which $a_{1}$ (or $a_{2}$) always remains a discount optimal policy. 
 To obtain the oscillatory behaviors of the robust value functions, we construct two convex compact uncertainty sets $\U_{s_{0}a_{1}}$ and $\U_{s_{0}a_{2}}$. At a high level,  these uncertainty sets, 
 while  distinct  subsets of $\Delta(\X)$, are such that  their boundaries overlap and intersect infinitely often. As $\gamma$ increases to 1, the worst-case transition probabilities $ \bm{P}^{a_{i},\star}_{\gamma}$ for the stationary policy $a_{i}$  ($i=1,2$)  varies  along the boundary of $\U_{s_{0}a_{i}}$. For each value of $\gamma$ such that the two worst-case transitions $\bm{P}^{a_{1},\star}_{\gamma}$ and $\bm{P}^{a_{2},\star}_{\gamma}$ coincide, the optimal policy changes, from $a_1$ to $a_2$ or vice-versa.
Our analysis also explains why such a pathological behavior cannot arise in nominal MDPs, nor when the uncertainty sets are polyhedral. In the former case,  the transition probabilities are fixed, and $\gamma \mapsto v^{\pi,\bm{P}}_{s,\gamma}$ is a well-behaved (rational) function. In the latter case, intersection points of the uncertainty sets are isolated, and worst-case transition probabilities are constant for $\gamma$ large enough, as shown in \citep{goyal2022robust}.
 

\vspace{2mm}
 \noindent \tb{Given Theorem \ref{th:sa-rec-no-blackwell}, it is natural to ask whether $\epsilon$-Blackwell optimal policies always exist, even when $\epsilon>0$ is very small. We answer this question in the affirmative in the next theorem. In fact, we show a significantly stronger result: we can find a stationary deterministic policies that remains $\epsilon$-Blackwell optimal policy for all values of $\epsilon>0$. The proof of the next theorem is in Appendix \ref{app:sa-rec-epsilon-blackwell-opt-same}. The main idea is to leverage Theorem \ref{th:eps-blackwell-compact-mdp}, i.e., the existence of a policy that is $\epsilon$-Blackwell optimal policies in the adversarial MDP.}
\begin{theorem}\label{th:sa-rec-epsilon-blackwell-opt-same}
Let $\U$ be an sa-rectangular compact uncertainty set. Then there exists a stationary deterministic policy that is $\epsilon$-Blackwell optimal for all $\epsilon>0$, i.e., $\exists \; \pi \in \Pi_{\sf SD}, \forall \; \epsilon>0,\exists \; \gamma_{\epsilon} \in (0,1)$ such that
\[  \min_{\bm{P} \in \U} (1-\gamma)R_{\gamma}(\pi,\bm{P}) \geq \sup_{\pi' \in \Pi_{\sf S}} \min_{\bm{P} \in \U} (1-\gamma)R_{\gamma}(\pi',\bm{P}) - \epsilon, \forall \; \gamma \in (\gamma_{\epsilon},1).\]
\end{theorem}
  Theorem \ref{th:sa-rec-no-blackwell} shows that the conclusion of Theorem \ref{th:sa-rec-epsilon-blackwell-opt-same} does not hold for $\epsilon=0$.  \tb{The reader may be surprised by the strength of this conclusion, since the {\em same} policy is $\epsilon$-Blackwell optimal, for all values of $\epsilon>0$. We will provide some more intuition on this type of policies and their relation with average optimality in the next section.}
\subsubsection{Limit behavior and connection with average optimality}\label{sec:sa-rec-relation-blackwell-opt-avg-opt}
 We now highlight the connection between Blackwell optimality and average optimality. Intuitively, $\epsilon$-Blackwell optimal policies remain approximately discount optimal for all $\gamma$ close to $1$. Since the discount factor captures the willingness of the decision-maker to wait for future rewards, we expect that the discounted return resembles more and more the average return as $\gamma \rightarrow 1$. We show that this intuition is correct for sa-rectangular uncertainty sets in the following theorem. 
\begin{theorem}\label{th:sa-rec-relation-avg-rew-blackwell}
\tb{Consider an sa-rectangular robust MDP with a compact uncertainty set $\U$. Let $\pi \in \PiSD$. Then 
\begin{center}
    $\pi$ is $\epsilon$-Blackwell optimal for all $\epsilon>0 \iff \pi$ is average optimal.
\end{center}} 
\end{theorem}
The proof of Theorem \ref{th:sa-rec-relation-avg-rew-blackwell} is presented in Appendix \ref{app:proof-th-relation-blackwell-avg}. The proof relies on a careful inspection of the limit behavior of the discounted return as $\gamma$ approaches $1$. We provide an outline of the main steps here.

Recall that for every $ (\pi,\bm{P}) \in \Pi_{\sf S} \times \U$ one has 
\[ \lim_{\gamma \rightarrow 1} (1-\gamma)R_{\gamma}(\pi,\bm{P}) = \Ravg(\pi,\bm{P}).\]
The first step is to observe that this equality still holds when taking the infimum over the transition probabilities. 
\begin{lemma}\label{lem:aux-1}
    In the setting of Theorem \ref{th:sa-rec-relation-avg-rew-blackwell}, let $\pi \in \Pi_{\sf S}.$ Then $\min_{\bm{P} \in \U} (1-\gamma) R_{\gamma}(\pi,\bm{P})$ admits a limit as $\gamma \rightarrow 1$ and
    \[ \lim_{\gamma \rightarrow 1} \min_{\bm{P} \in \U} (1-\gamma) R_{\gamma}(\pi,\bm{P}) = \inf_{\bm{P} \in \U} \Ravg(\pi,\bm{P}).\]
\end{lemma}
Lemma \ref{lem:aux-1} is a consequence of Corollary 5.26 in \cite{sorin2002first}. We still provide the proof for completeness in Appendix \ref{app:proof-th-relation-blackwell-avg}.
We next show that the same conclusion holds for the limit of the {\em optimal} discounted return.
\begin{lemma}\label{lem:aux-2}
    In the setting of Theorem \ref{th:sa-rec-relation-avg-rew-blackwell}, $\max_{\pi \in \Pi_{\sf SD}} \min_{\bm{P} \in \U} (1-\gamma)R_{\gamma}(\pi,\bm{P})$ admits a limit as $\gamma \rightarrow 1$ and
    \[ \lim_{\gamma \rightarrow 1} \max_{\pi \in \Pi_{\sf SD}} \min_{\bm{P} \in \U} (1-\gamma)R_{\gamma}(\pi,\bm{P}) = \max_{\pi \in \Pi_{\sf SD}} \inf_{\bm{P} \in \U} \Ravg(\pi,\bm{P}).\]
\end{lemma}
The last part of the proof relates average optimality and policies remaining $\epsilon$-Blackwell optimality for all $\epsilon>0$.
\begin{lemma}\label{lem:aux-3}
    In the setting of Theorem \ref{th:sa-rec-relation-avg-rew-blackwell}, let $\pi \in \Pi_{\sf SD}$ be $\epsilon$-Blackwell optimal for any $\epsilon>0$. Then  \[ \lim_{\gamma \rightarrow 1} \min_{\bm{P} \in \U} (1-\gamma) R_{\gamma}(\pi,\bm{P}) = \lim_{\gamma \rightarrow 1} \max_{\pi' \in \Pi_{\sf SD}} \inf_{\bm{P} \in \U} (1-\gamma)R_{\gamma}(\pi',\bm{P}).\]
\end{lemma}
Theorem \ref{th:sa-rec-relation-avg-rew-blackwell} follows when combining Lemma \ref{lem:aux-3} and Lemma \ref{lem:aux-2}.
Note that all the conclusions from Lemma \ref{lem:aux-1}, Lemma \ref{lem:aux-2} and Lemma \ref{lem:aux-3} are true when we replace $\Ravg$ as defined in \eqref{eq:definition-average-return-exp-limsup} by the other natural definitions $\eqref{eq:liminf-exp}, \eqref{eq:limsup-exp}, \mbox{ or } \eqref{eq:exp-liminf}$, since they only involve optimizing over $\Pi_{\sf S}$ and $\U$, where all definitions of the average return coincide.

\tb{
We conclude this section with the following theorem, which relates discount optimal policies (for large discount factors) and average optimal policies, in the same way that Theorem \ref{th:sa-rec-relation-avg-rew-blackwell} relates $\epsilon$-Blackwell optimal policies and average optimal policies.
\begin{theorem}\label{th:sa-rec-relation-avg-rew-discount}
    Consider an sa-rectangular robust MDP with a compact uncertainty set $\U$. There exists a discount factor $\gammaa \in [0,1)$ such that any stationary deterministic discount optimal policy for $\gamma \in (\gammaa,1)$ is also average optimal.
\end{theorem}
The proof of Theorem \ref{th:sa-rec-relation-avg-rew-discount} combines the fact that worst-case discount returns converge to worst-case average returns (Lemma \ref{lem:aux-1}) with the finiteness of the set of stationary deterministic policies, which are optimal for both the discount return and the average return. We provide the detailed proof in Appendix \ref{app:proof relation avg rew discount}. In principle, one could use the results from Theorem \ref{th:sa-rec-relation-avg-rew-discount} to compute an average optimal policy by solving a discount robust MDP with a discount factor larger than $\gammaa$. However, our proof of Theorem \ref{th:sa-rec-relation-avg-rew-discount} only shows the existence of such $\gammaa \in [0,1)$, it is not constructive. We leave designing an efficient algorithm for finding $\gammaa$ (or an upper bound on $\gammaa$) as an interesting open question.
}
\begin{remark}
    Lemma \ref{lem:aux-2} and Lemma \ref{lem:aux-3} are present as Corollary 4 in \cite{tewari2007bounded}, under the assumption that $\U$ is based on $\ell_{\infty}$-distance. This considerably simplifies the proof, as in this setting, the number of extreme points of $\U$ is finite, so the infimum over $\U$ can be reduced to a minimization over a finite number of elements, and we can exchange the minimization with the limit.
    \tb{Lemma \ref{lem:aux-1} and Theorem \ref{th:sa-rec-relation-avg-rew-discount} are also present in \cite{wang2023robust} (see Theorem 2 and Theorem 5 in \cite{wang2023robust}), under some additional assumptions ({\em unichain} compact sa-rectangular RMDPs). Our results hold in more generality and show that unichain assumption is unnecessary here}.
\end{remark}

\subsubsection{Definable robust Markov decision processes}\label{sec:definable-sa-rec}
In this section, we introduce a  general class of uncertainty sets that encompasses virtually all existing examples in the RMDP literature, and for which stationary deterministic Blackwell optimal policies exist. Indeed, the classical approach is to construct uncertainty sets based on simple functions, like affine maps, $\ell_{p}$-balls, or Kullback-Leibler divergence, see Section \ref{sec:RMDPs}. For such simple functions,  we intuitively do not expect that the robust value functions oscillate and intersect infinitely often. Our main contribution in this section is to formalize this intuition with the notion of {\em definability}~\citep{van1998minimal,coste2000introduction,bolte2015definable} and to prove Theorem \ref{th:definable-sa-rec-blackwell-opt}, which states that for definable sa-rectangular RMDPs there always exists a stationary deterministic Blackwell optimal policy. 

\vspace{2mm}

\noindent
{\bf A concise introduction to definability.}
We start with the following definition. Intuitively, a set is {\em definable} if it can be ``constructed" based on polynomials, the exponential function, and canonical projections (elimination of variables). 
\begin{definition}[Definable set and definable function]\label{def:definability}
    A subset of $\R^{n}$ is {\em definable} if it is the image, under a canonical projection $\R^{n+k} \rightarrow \R^{n}$ that eliminates any set of $k$ variables, of a set of the form
    \begin{equation}\label{eq:P-exp-def}
        \{ \bm{x} \in \R^{n+k} \; | \; {\sf Poly}(x_{1},...,x_{n+k},\exp(x_{1}),...,\exp(x_{n+k})) =0\}
    \end{equation}
where ${\sf Poly}(\cdot)$ is a real polynomial in $2(n+k)$ variables. 

For $\Omega \subset \R^{n}$ a definable subset of $\R^{n}$, a function $f: \Omega \rightarrow \R^{m}$ is definable if its graph $\{ (\bm{x},\bm{y}) \in \Omega \times \R^{m} \; | \; \bm{y}=f(\bm{x})\}$ is a definable subset of $\R^{n+m}$.
\end{definition}
We refer to \cite{bolte2015definable} and \cite{akian2019operator} for concise introductions to definability and to \cite{coste2000introduction} for a more in-depth treatment. It is instructive to consider a few simple examples. 
\begin{example}[Affine functions.]\label{ex:affine-maps-definable}
    Consider an affine function $f:\bm{x} \mapsto \bm{a}\tr\bm{x} + b$ for some $\bm{a} \in \R^{\X},b \in \R$. The graph of $f$ is $\{ (\bm{x},y) \in \R^{\X} \times \R \; | \; \sum_{s \in \X} a_{s}x_{s} + b - y = 0\}$, which can be written as \eqref{eq:P-exp-def}. Therefore, affine functions are definable.
\end{example}
\tb{
\begin{example}[Bilinear functions.]\label{ex:bilinear-maps-definable}
    Consider a bilinear function $f:(\bm{z},\bm{x}) \mapsto \bm{z}\tr\bm{x}$. The graph of $f$ is $\{ (\bm{z},\bm{x},y) \in \R^{\X} \times \R^{\X} \times \R \; | \; \sum_{s \in \X} z_{s}x_{s} - y = 0\}$, which can be written as \eqref{eq:P-exp-def} for ${\sf Poly}(\bm{z},\bm{x},y) = \sum_{s \in \X} z_{s}x_{s} - y$. Therefore, bilinear functions are definable.
\end{example}
}
Definable functions and definable sets are well-behaved under many useful operations, as shown in the following lemma. It follows from the definition and from some properties shown in \cite{bolte2015definable}. For completeness, we provide the precise references (and some proofs when necessary) in Appendix \ref{app:proof lemma definability stability}.
\begin{lemma}\label{lem:definability-stability}
        \begin{enumerate}
        \item The only definable subsets of $\R$ are the {\em finite} union of open intervals and singletons.
         \item If $A,B \subset \R^{n}$ are definable sets, then $A \cup B,A \cap B$ and $\R^{n} \setminus A$ are definable sets.
        \item Let $f,g$ be definable functions. Then $f \circ g$, $-f$, $f+g$, $f \times g$ are definable. 
        \item For $A,B$ two definable sets and $g: A \times B \rightarrow \R$ a definable function, then the functions $\bm{x} \mapsto \inf_{\bm{y} \in B} g(\bm{x},\bm{y})$ and $\bm{x} \mapsto \sup_{\bm{y} \in B} g(\bm{x},\bm{y})$ (defined over $A$) are definable functions. 
        \item 
        If $A,B\subseteq \R$ and $g:A\rightarrow \R$ are definable then $g^{-1}(B)$ is a definable set.
        \end{enumerate}
\end{lemma}
Based on Lemma \ref{lem:definability-stability}, we provide some examples of definable and non-definable functions below.
\begin{example}[$\ell_{p}$-norms are definable.]\label{ex:norms-definable}
    Consider an $\ell_{p}$-norm for $p \in \N$. Then its graph is $\{ (\bm{x},y) \in \R^{\X} \times \R \; | \; \sum_{s \in \X} |x_{s}|^{p} - y^{p} = 0, y \geq 0\}$, which is a definable set. Therefore, $\ell_{p}$-norms are definable.
\end{example}
\begin{example}[Logarithm and exponential are definable.]
    It is straightforward that $x \mapsto \exp(x)$ is a definable function. The graph of the logarithm is $\{ (x,y) \in \R_{+}^{*} \times \R \; | \; y = \log(x)\}$. This can be rewritten $\{ (x,y) \in \R_{+}^{*} \times \R \; | \; \exp(y) = x\}$, which is a definable set. Therefore, the logarithm is a definable function.
\end{example}

\begin{example}[Some functions based on entropy are definable.]\label{ex:KL-burg}
   Fix $\hat{\bm{p}} \in \Delta(\X)$ with $\hat{p}_{s} >0, \forall \; s \in \X$, and consider the Kullback-Leibler divergence: $\bm{p} \mapsto \sum_{s \in \X} p_{s} \log(p_{s}/\hat{p}_{s})$, defined over $\bm{p} \in \Delta(\X)$.
    Recall that $x \mapsto \log(x)$ is definable, so that $x \mapsto x\log(x)$ is also definable. By summation, the Kullback-Leibler divergence is a definable function. The case of the Burg entropy $\bm{p} \mapsto \sum_{s \in \X} \hat{p}_{s} \log(\hat{p}_{s}/p_{s})$ for $\bm{p} \in \Delta(\X)$ is similar.
\end{example}
\tb{
\begin{example}[Wasserstein distance is definable.]\label{ex:wasserstein}
   Fix $\hat{\bm{p}} \in \Delta(\X)$ and consider the Wasserstein distance from $\hat{\bm{p}}$, defined as $W:(\cdot,\hat{\bm{p}}):\bm{p} \mapsto \min \left\{ \sum_{(s,s') \in \X \times \X} A_{ss'}c_{ss'} \; | \; \bm{A} \in \R^{\X \times \X}_{+},\bm{A}\bm{e}=\bm{p},\bm{A}\tr\bm{e} = \hat{\bm{p}}\right\}$, for a given cost function $
\bm{c} \in \R^{\X \times \X}_{+}$. From Kantorovich duality (e.g. Chapter 1 in \cite{villani2021topics}), we obtain that $W:(\bm{p},\hat{\bm{p}}) = \max \left\{ \bm{p}\tr\bm{x} + \hat{\bm{p}}\tr\bm{y} \; | \; x_{s} + y_{s'} \leq c_{ss'}, \forall \; (s,s') \in \X \times \X \right\}.$ From Example \ref{ex:bilinear-maps-definable}, we know that $(\bm{x},\bm{y},\bm{p},\hat{\bm{p}}) \mapsto \bm{p}\tr\bm{x} + \hat{\bm{p}}\tr\bm{y}$ is definable, and from Example \ref{ex:affine-maps-definable}, we know that $\{(\bm{x},\bm{y}) \in \R^{\X} \times \R^{\X} \; | \;x_{s} + y_{s'} \leq c_{ss'}, \forall \; (s,s') \in \X \times \X\}$ is a definable set. Combining these facts with the fourth statement of Lemma \ref{lem:definability-stability}, we obtain that the Wasserstein distance is definable.
\end{example}
}

\tb{
\begin{example}[$\sin$ is not definable.]\label{ex:non-def-1}
Consider the sinus function: $\sin:\R \mapsto \R$ and note that $\{0\}$ is a definable set. If $\sin$ was definable, then $\sin^{-1}(\{0\})$ would be a definable set of $\R$, i.e. following the first point in Lemma \ref{lem:definability-stability} it would be a {\em finite} union of open intervals and singletons. However, $\sin^{-1}(\{0\}) = \{k \pi \; | \; k \in \N\}$ is infinite, and therefore $\sin$ is not definable. To extend this counterexample to the case of a function defined over a bounded interval, we can show in a similar fashion that $x\mapsto \sin(1/x)$, defined over $(0,1]$, is not definable.
\end{example}}

\tb{
\begin{example}[A convex function that is not definable.]
We can construct a convex, non-definable function as follows. Consider $F:\R_{+}\rightarrow \R, x \mapsto x^2$, which is convex and definable. We construct the function $f$ as follows: $f:
R_{+} \rightarrow \R$ is the piece-wise affine function that changes slope at every integer $k \in \N$ and which coincides with $F$ at every integer $k \in \N$. It is straightforward to see that $f$ is convex (since $F$ is convex). However, if $f$ was definable, then $f-F$ would be definable (from the third statement in Lemma \ref{lem:definability-stability}). This would imply that $(f-F)^{-1}(\{0\})$ is a definable set of $\R$. However, $(f-F)^{-1}(\{0\})=\N$, and $\N$ is not a definable set of $\R$ (it is not {\em finite} union of open intervals and singletons). This shows that $f$ is not definable. Similarly as for the $\sin$ function, we can extend this example to $(0,1]$ by considering a piece-wise affine function $(0,1] \mapsto \R$ that coincides with $F$ over $\{1/k \; | \; k \in \N, k \geq 1\}$.
\end{example}}

\tb{
The previous two counterexamples highlight an important feature of definable functions of one real variable: they cannot take the same value {\em a countably infinite number of times}. Since we have highlighted the oscillations of value functions as the main issue behind the potential non-existence of Blackwell optimal policies in Section \ref{sec:counter-example}, the notion of {\em definability} provides an adequate fix to these potential pathological oscillations. In particular, the most important result for us is the following monotonicity result for definable functions of a real variable. 
}

\begin{theorem}[Theorem 2.1, \cite{coste2000introduction}]\label{th:monotonicity-theorem}
Let $f:(a,b) \rightarrow \R$ be a definable function. Then there exists a finite subdivision of the interval $(a,b)$ as $a = a_{1} < a_{2} < \dots < a_{k} = b$ such that on each $(a_{i},a_{i+1})$ for $i=1,...,k-1$, $f$ is continuous and either constant or strictly monotone.
\end{theorem}
This result shows that definable functions over $\R$ cannot oscillate infinitely often on an interval. As we have identified oscillations of the robust value functions as the main issue potentially precluding the existence of Blackwell optimal policies (see discussion after Theorem \ref{th:sa-rec-no-blackwell}), it is not surprising that Theorem \ref{th:monotonicity-theorem} will play an important role in the proof of existence of Blackwell optimality for definable RMDPs, see  Theorem \ref{th:definable-sa-rec-blackwell-opt} below.

\begin{remark} The notion of definability is actually more general, see for instance~\cite{van1998minimal,coste2000introduction,bolte2015definable}. We chose to introduce only the notions needed here. We simply note that the notion of definability introduced in Definition \ref{def:definability} is usually referred to as {\em definable in the real exponential field}  in the literature.
\end{remark}
\vspace{2mm}

\noindent
{\bf Definable robust MDPs.} We now study sa-rectangular robust MDPs with a definable compact uncertainty set. We first note that this encompasses the vast majority of the uncertainty sets studied in the robust MDP literature. Indeed, from Lemma \ref{lem:definability-stability}, we know that the set $
        \U_{sa} = \{ \bm{p} \in \Delta(\X) \; | \; d_{sa}\left(\bm{p}\right)\leq \alpha_{sa}\}$
is definable as soon as $d_{sa}:\R^{\X} \rightarrow \R$ is definable. From the various examples introduced above, we obtain that sa-rectangular uncertainty sets based on $\ell_{\infty}$-norm~\citep{givan1997bounded}, $\ell_{2}$-norm~\citep{iyengar2005robust}, $\ell_{1}$-norm~\citep{ho2021partial}, Kullback-Leibler divergence and Burg entropy~\citep{iyengar2005robust,ho2022robust} are definable uncertainty sets.

We start by showing that robust value functions of definable RMDPs are themselves definable functions. 
\begin{proposition}\label{prop:definable-value-function}
    Assume that $\U$ is sa-rectangular and definable. Then for any policy $\pi \in \Pi_{\sf S}$, the function $\gamma \mapsto \bm{v}^{\pi,\U}_{\gamma}$ is definable.
\end{proposition}
%
\proof{Proof of Proposition \ref{prop:definable-value-function}.}
Recall that the robust value function $\bm{v}^{\pi,\U}_{\gamma}$ is the unique fixed-point of the  operator $\bm{v} \mapsto T^{\pi,\U}_{\gamma}(\bm{v})$, where $T^{\pi,\U}_{\gamma}:\R^{\X} \rightarrow \R^{\X}$ is defined as
\[T^{\pi,\U}_{\gamma,s}(\bm{v}) = \min_{\bm{P}_{s} \in \U_{s}} \sum_{a \in \A} \pi_{sa} \bm{p}_{sa}\tr\left( \bm{r}_{sa} + \gamma \bm{v}\right), \forall \; s \in \X, \forall \; \bm{v} \in \R^{\X}.\]
The proof proceeds in two steps. We first show that $T^{\pi,\U}_{\gamma}$ is definable when $\U$ is definable. We next prove that this implies  that the robust value functions are definable.
\paragraph{First step.} We first prove that $(\bm{v},\gamma) \mapsto T^{\pi,\U}_{\gamma}(\bm{v})$ is definable.

Since a function is definable if and only if each of its components is definable (Exercise 1.10, \cite{coste2000introduction}), we need to show that 
$(\bm{v},\gamma) \mapsto \left(\sum_{a \in \A} \pi_{sa}\min_{\bm{p}_{sa} \in \U_{sa}}\bm{p}_{sa}\tr\left(\bm{r}_{sa}+\gamma  \bm{v}\right)\right)$ is definable for fixed $s\in \X$.

The map $(\gamma,\bm{v}, \bm{P}_s=(\bm{p}_{sa})_{a\in \A})\in (0,1)\times \R^\X \times \U_s\mapsto 
\sum_{a \in \A} \pi_{sa}\bm{p}_{sa}\tr\left(\bm{r}_{sa}+\gamma  \bm{v}\right)
$ is polynomial, and defined on a definable set, hence it is a definable function. 

From Lemma \ref{lem:definability-stability}, this implies that 
\[(\gamma,v)\in (0,1)\times \R^{\X} \mapsto \min_{\bm{P}_s\in \U_s} \sum_{a \in \A} \pi_{sa}\bm{p}_{sa}\tr\left(\bm{r}_{sa}+\gamma  \bm{v}\right)\]
is also definable. 


\paragraph{Second step.} We now show that $\gamma \mapsto \bm{v}^{\pi,\U}_{\gamma}$ is definable.

We need to show that the graph of   $\gamma \mapsto \bm{v}^{\pi,\U}_{\gamma}$ is a definable set. Since $\bm{v}^{\pi,\U}_{\gamma}$ is the unique fixed-point of the operator $T^{\pi,\U}_{\gamma}$, 
this graph  $\{ \left(\bm{v},\gamma\right) \in \R^{\X} \times \R \; | \; \bm{v} = \bm{v}^{\pi,\U}_{\gamma} \}$ is also equal to $\{ \left(\bm{v},\gamma\right) \in \R^{\X} \times (0,1) \; | \; T^{\pi,\U}_{\gamma}(\bm{v}) = \bm{v}\}$, which is the projection over $\R^\X \times (0,1)$ of the intersection 
\[ \{\left(\bm{v},\gamma,\bm{z}\right) \in \R^{\X}  \times \R \times \R^{\X} \; | \; \bm{z} = T^{\pi,\U}_{\gamma}(\bm{v})\}\cap \{\left(\bm{v},\gamma,\bm{z}\right) \R^{\X}  \times \R \times \R^{\X} \; | \; \bm{z} = \bm{v}\}.\]
The first of these two sets is the graph of  $\left(\bm{v},\gamma\right) \mapsto T^{\pi,\U}_{\gamma}(\bm{v})$, which is definable using the first step. The second set is the set of zeroes of the affine function $(\bm{v},\bm{z}) \mapsto \bm{v}-\bm{z}$, and is it definable. Therefore their intersection is definable.
\hfill \halmos \endproof

\medskip

Proposition \ref{prop:definable-value-function} shows an appealing property of the robust value functions of definable robust MDPs. From the monotonicity theorem, and since the difference of definable functions is definable, it implies that two robust value functions can only intersect finitely many times (or be equal on an entire interval close to $1$), which precludes the pathological behavior underpinning 
Theorem \ref{th:sa-rec-no-blackwell}. 
We are now ready to state our main theorem for this section, which formalizes this intuition and shows that Blackwell optimal policies exist for the class of sa-rectangular definable RMDPs.
\begin{theorem}\label{th:definable-sa-rec-blackwell-opt}
    Consider an sa-rectangular robust MDP with a definable compact uncertainty set $\U$. Then there exists a stationary deterministic Blackwell optimal policy:
    \[\exists \; \pi \in \Pi_{\sf SD}, \exists \; \gamma_{0} \in (0,1), \forall \; \gamma \in (\gamma_{0},1), \min_{\bm{P} \in \U} R_{\gamma}(\pi,\bm{P}) \geq \sup_{\pi' \in \Pi_{\sf S}} \min_{\bm{P} \in \U} R_{\gamma}(\pi',\bm{P}).\]
\end{theorem}
\proof{Proof of Theorem \ref{th:definable-sa-rec-blackwell-opt}}
From Proposition \ref{prop:definable-value-function} and Lemma \ref{lem:definability-stability}, the function $[0,1) \rightarrow \R, \gamma \mapsto v^{\pi,\U}_{\gamma,s} - v^{\pi',\U}_{\gamma,s}$ is definable for each pair $(\pi,\pi')$  of stationary deterministic policies. \tb{From the monotonicity theorem (Theorem \ref{th:monotonicity-theorem})}, it follows that $\gamma \mapsto v^{\pi,\U}_{\gamma,s} - v^{\pi',\U}_{\gamma,s}$ does not change signs in a neighborhood $(1-\eta(\pi,\pi',s),1)$ with $\eta(\pi,\pi',s) > 0$. We then define $\bar{\gamma} = \max_{\pi,\pi' \in \Pi_{\sf SD},s \in \X} 1-\eta(\pi,\pi',s)$ with $\bar{\gamma}<1$ since there are finitely many stationary deterministic policies and finitely many states. Any policy that is discount optimal for $\gamma \in (\bar{\gamma},1)$ is Blackwell optimal, which shows the existence of a stationary deterministic Blackwell optimal policy.
\hfill \halmos \endproof
 The proof of Theorem \ref{th:definable-sa-rec-blackwell-opt} is concise, and proves in addition that above some  {\em Blackwell discount factor} $\gamma_{\sf bw} \in (0,1)$,  any discount optimal policy is also Blackwell optimal. We state this as an independent result.
 
\begin{theorem}\label{th:blackwell-discount-factor}
Consider an sa-rectangular robust MDP with a definable uncertainty set $\U$. Then there exists a {\em Blackwell discount factor} $\gammab \in (0,1)$, such that any policy that is discount optimal for some $\gamma \in (\gammab,1)$ is also Blackwell optimal.
\end{theorem}
The existence of the Blackwell discount factor is shown in \cite{grand2024reducing} for nominal MDPs and for sa-rectangular RMDPs with polyhedral uncertainty sets. Theorem \ref{th:blackwell-discount-factor} extends this result to the larger class of definable uncertainty sets. As a consequence, a Blackwell optimal policy can be computed by solving a discounted robust MDP with discount factor $\gamma = \hat{\gamma}$, as soon as an upper bound $\hat{\gamma}<1$ on $\gammab$ is known. We leave the determination of such an upper bound on $\gammab$ as an interesting future direction. \tb{We also note the similarities between Theorem \ref{th:blackwell-discount-factor} and Theorem \ref{th:sa-rec-relation-avg-rew-discount}, which shows the existence of a discount factor $\gammaa \in [0,1)$ such that any discount optimal policy is average optimal, provided that the discount factor is larger than $\gammaa$. We note that Theorem \ref{th:sa-rec-relation-avg-rew-discount} does not require the definability assumption, since it only highlights the connection between discount optimal policies and average optimal policies, both of which always exist. In contrast, Theorem \ref{th:blackwell-discount-factor} highlights the connection between average optimal policies and Blackwell optimal policies, and the latter may not always exist, as shown in Theorem \ref{th:sa-rec-no-blackwell}.}
\subsection{The case of s-rectangular robust MDPs}\label{sec:s-rec-rmdp-blackwell}
We here provide an instance of an s-rectangular RMDPs with no Blackwell optimal policy, even though the uncertainty set is definable. This may happen when all optimal policies are randomized, and the randomized optimal policies vary with the discount factor. 
\begin{proposition}\label{prop:no-BO-srec}
    There exists an s-rectangular robust MDP instance, with a compact, polyhedral uncertainty set $\U$, such that there is no Blackwell optimal policy:
    \[\forall \; \pi \in \Pi_{\sf H}, \forall \; \gamma \in (0,1), \exists \; \gamma' \in (\gamma,1), \min_{\bm{P} \in \U} R_{\gamma'}(\pi,\bm{P}) < \sup_{\pi' \in \Pi_{\sf H}} \min_{\bm{P} \in \U} R_{\gamma'}(\pi',\bm{P}). \]
\end{proposition}
\proof{Proof.}
    
\begin{figure}[htb]
\begin{center}
    \begin{subfigure}{0.4\textwidth}
\includegraphics[width=\linewidth]{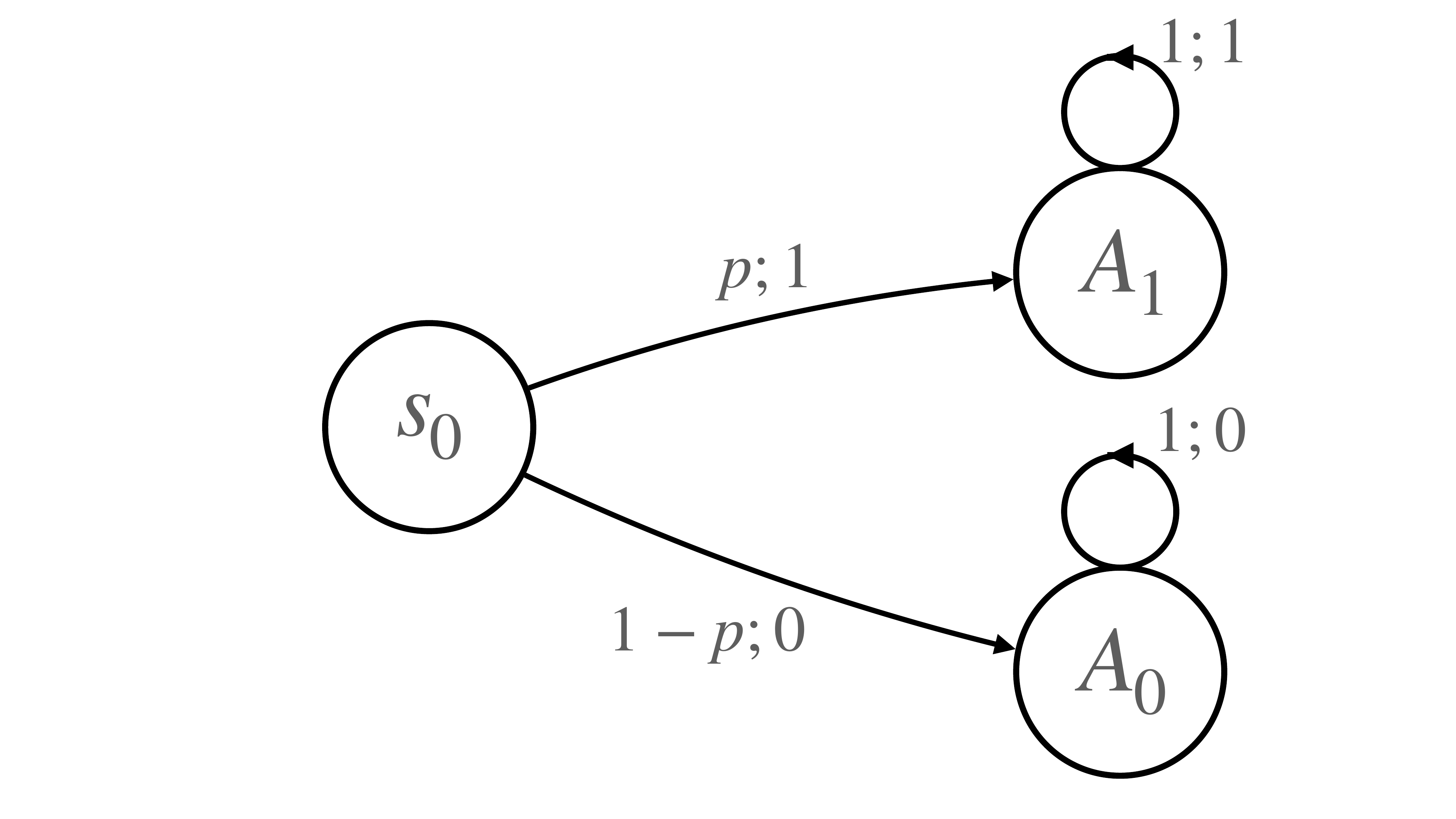}
    \caption{Transition for action $a_{1}$}
    \label{fig:rmdp_a1}
    \end{subfigure}
    \begin{subfigure}{0.4\textwidth}
\includegraphics[width=\linewidth]{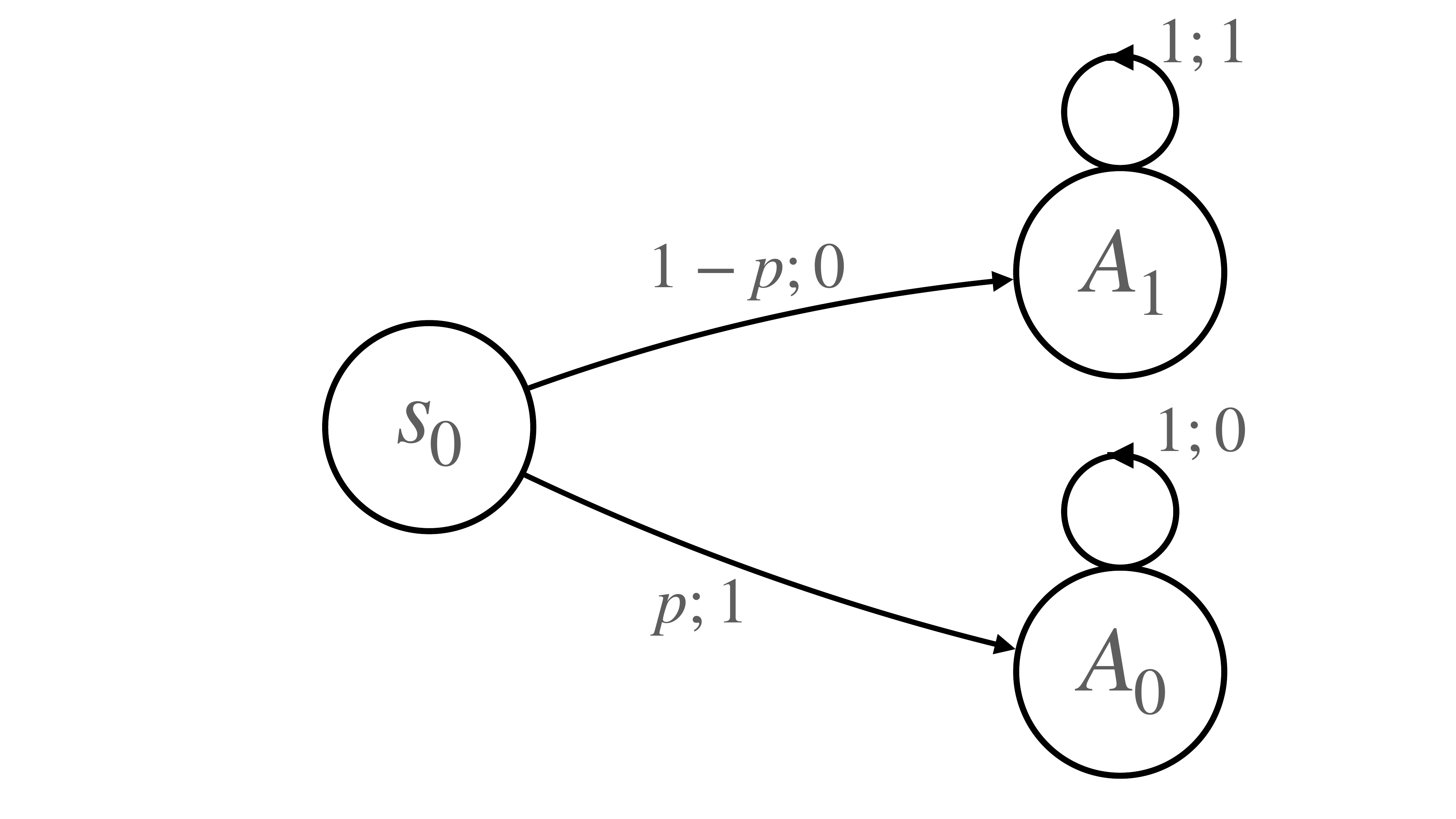}
    \caption{Transition for action $a_{2}$}
    \label{fig:rmdp_a2}
    \end{subfigure}
    \end{center}
\caption{Transitions and rewards for an s-rectangular RMDP instance with no Blackwell optimal policy.}    \label{fig:rmdp-srec-noblackwell}
\end{figure}

In Figure \ref{fig:rmdp-srec-noblackwell}, we adapt the example from Figure 3 in \cite{wiesemann2013robust}. There are three states $s_{0},A_{1},A_{0}$ and two actions $a_{1}$ and $a_{2}$.
We can parametrize the return by the probability $x \in [0,1]$ to play action $a_{1}$ and by the parameter $p \in [0,1]$ chosen by the adversary. We have
$R_{\gamma}(x,p) = x p  \left(1 + \frac{\gamma}{1-\gamma} \right) + (1-x) \left(p + (1-p) \frac{\gamma}{1-\gamma}\right).$
This can be reformulated
$R_{\gamma}(x,p) = p  \left(1 + x \frac{\gamma}{1-\gamma}\right) + (1-p) \left( (1-x) \frac{\gamma}{1-\gamma}\right).$
Let us compute the worst-case return for policy $x \in [0,1]$. We have $
    \min_{p \in [0,1]} R(x,p) = \min \{ 1 + x \frac{\gamma}{1-\gamma}, (1-x) \frac{\gamma}{1-\gamma}\}.
$
The optimal policy then maximizes $x \mapsto \min_{p \in [0,1]} R(x,p)$, i.e., it maximizes $x \mapsto \min \{ 1 + x \frac{\gamma}{1-\gamma}, (1-x) \frac{\gamma}{1-\gamma}\}$ over $[0,1]$. Assume that $\gamma \geq 1/2$ (we are interested in Blackwell optimal policies, i.e. in the case $\gamma \rightarrow 1$). Then the maximum of $x \mapsto \min \{ 1 + x \frac{\gamma}{1-\gamma}, (1-x) \frac{\gamma}{1-\gamma}\}$ is attained at $x^{*}(\gamma)$, solution of the equation $1 + x \frac{\gamma}{1-\gamma}= (1-x) \frac{\gamma}{1-\gamma}$. Therefore, there is a unique discount optimal policy and it depends on $\gamma$: $x\opt(\gamma) = 1 - \frac{1}{2 \gamma}$ for $\gamma \geq 1/2$.  Overall, no stationary policies remain discount optimal when $\gamma$ varies. Since the discount optimal policy is always unique and stationary in this RMDP instance, there are no stationary Blackwell optimal policies. 
\tb{Since the process moves to either state $A_{1}$ or state $A_{0}$ after one period, the actions of the DM are relevant in the first period only. Therefore, there is no history-dependent Blackwell optimal policy either.}
\hfill \halmos \endproof
We note that the counterexample in the proof of Proposition \ref{prop:no-BO-srec} is much simpler than the counterexample from Section \ref{sec:counter-example} for sa-rectangular RMDPs. In particular, Proposition \ref{prop:no-BO-srec} shows that Blackwell optimal policies may fail to exist for s-rectangular RMDPs, even in the simple setting where the uncertainty set is a polytope parametrized by $p \in [0,1]$. As a possible next step, it sounds promising to study the existence and tractability of $\epsilon$-Blackwell optimal policies for s-rectangular RMDPs.

\section{Algorithms}\label{sec:avg-alg}
In this section, we discuss various iterative algorithms to compute average optimal and Blackwell optimal policies. Our results for s-rectangular RMDPs from Section \ref{sec:s-rec-avg} and Section \ref{sec:s-rec-rmdp-blackwell} suggest that it may be difficult to compute average and Blackwell optimal policies in this case. Therefore we focus on sa-rectangular robust MDPs. Since for sa-rectangular RMDPs an average optimal policy may be chosen stationary and deterministic, we define the optimal gain $\bm{g}\opt \in \R^{\X}$ as
\begin{equation}\label{eq:optimal-gain-definition}
        g\opt_{s} = \max_{\pi \in \Pi_{\sf SD}} \inf_{\bm{P} \in \U} \lim_{T \rightarrow + \infty}  \E_{\pi,\bm{P}} \left[ \frac{1}{T+1} \sum_{t=0}^{T} r_{s_{t}a_{t}s_{t+1}} \; | \; s_{0} = s \right], \forall \; s \in \X.
    \end{equation}
We now introduce three algorithms to compute the optimal gain $\bm{g}\opt$ when $\U$ is sa-rectangular.
\subsection{Algorithms based on large discount factor}
\tb{ Theorem \ref{th:sa-rec-relation-avg-rew-discount} shows the existence of a discount factor $\gammaa \in (0,1)$, such that any policy that is $\gamma$-discount optimal for $\gamma \in (\gammaa,1)$ is also average optimal. Therefore, we can compute average optimal policies by solving a sequence of discounted RMDPs with discount factors increasing to $1$. An immediate consequence is the following theorem.
}
\begin{theorem}\label{th:convergence-alg-definable-uset-0}
    \tb{Let $\U$ be a definable compact sa-rectangular uncertainty set and let $\left(\bm{v}^{t}\right)_{t \geq 1}$ be the iterates of Algorithm \ref{alg:sarec-limit-discount-values}. Let $\left(\pi^{t}\right)_{t \in \N}$ be a sequence of policies such that $\pi^{t} \in \Pi_{\sf S}$ and $(1-\gamma_t)\bm{v}^{\pi^{t},\U}_{\gamma_{t-1}} = \bm{v}^{t}$. Then $\left(\bm{v}^{t}\right)_{t \geq 1}$ converges to $\bm{g}\opt \in \R^{\X}$, and $\pi^{t}$ is average optimal for $t \in \N$ large enough.}
\end{theorem}
\tb{Unfortunately, computing $\gammaa$ appears challenging, so that we are not able to provide a convergence rate for Algorithm \ref{alg:sarec-limit-discount-values}. }
\tb{In fact, since Blackwell optimal policies are also average optimal, then we can also compute average optimal policies by computing discount optimal policies for $\gamma > \gammab$, the Blackwell discount factor introduced in Theorem \ref{th:blackwell-discount-factor}. When $\U$ is polyhedral, \cite{grand2024reducing} obtains an upper bound on $\gammab$ (see Theorem 4.4 in \cite{grand2024reducing}), but this upper bound is too close to $1$ to be of practical use, even for very simple instance. We provide an example of this issue in Appendix \ref{app:issue-gamma-b}. We also note that Algorithm \ref{alg:sarec-limit-discount-values} requires solving a robust MDP at every iteration, which may be computationally intensive. We now describe two algorithms with a lower per-iteration complexity.}
\begin{algorithm}
\caption{Iterative algorithm for average optimality: limit of discounted values}\label{alg:sarec-limit-discount-values}
\begin{algorithmic}[1]
\State {\bf Input:} An increasing sequence of discount factors $\left(\gamma_t\right)_{t \in \N}$ with $\lim_{t \rightarrow + \infty} \gamma_t = 1$.
\State Initialize $\bm{v}^{0} = \bm{0} \in \R^{\X},\pi^{0} \in \Pi_{\sf S}$.

\For{$t \in \N$}
    \State $\bm{v}^{t+1} = (1-\gamma_{t})\bm{v}^{\star,\U}_{\gamma_t}$ 
\EndFor
\end{algorithmic}
\end{algorithm}
\subsection{Algorithms based on value iteration}
We now introduce two value iteration algorithms to compute the optimal worst-case average return when the uncertainty set $\U$ is sa-rectangular {\em and} definable. 

\vspace{2mm}

\noindent 
{\bf Algorithm based on increasing horizon.}
\begin{algorithm}
\caption{Iterative algorithm for average optimality: increasing horizon}\label{alg:sarec-increasing-horizon}
\begin{algorithmic}[1]
\State Initialize $\bm{v}^{0} = \bm{0} \in \R^{\X}$.

\For{$t \in \N$}
    \State $v^{t+1}_{s} = \max_{a \in \A} \min_{\bm{p}_{sa} \in \U_{sa}} \bm{p}_{sa}\tr\left(\bm{r}_{sa} + \bm{v}^{t}\right), \forall \; s \in \X$. 
\EndFor
\end{algorithmic}
\end{algorithm}
We start with Algorithm \ref{alg:sarec-increasing-horizon}. \tb{This algorithm computes the optimal value functions of a finite-horizon RMDP as the horizon increases to $+\infty$, which can be done by dynamic programming, as in Step 3 of this algorithm, where the optimal value function $\bm{v}^{t+1}$ in horizon $t+1$ is computed given the optimal value function $\bm{v}^{t}$ in horizon $t$.} \tb{Algorithm \ref{alg:sarec-increasing-horizon} is inspired from the vast literature on the {\em asymptotic approach} in stochastic games~\citep{everett1957recursive,kohlberg1974repeated,bewley1976asymptotic}, studying the convergence of the finite-horizon values as the players become more and more patients, i.e. as the horizon increases to $+\infty$.}
The following theorem shows that Algorithm \ref{alg:sarec-increasing-horizon} is correct.
\begin{theorem}\label{th:convergence-alg-definable-uset-1}
    Let $\U$ be a definable compact sa-rectangular uncertainty set and let $\left(\bm{v}^{t}\right)_{t \geq 1}$ be the iterates of Algorithm \ref{alg:sarec-increasing-horizon}. Then $\left(\frac{1}{t}\bm{v}^{t}\right)_{t \geq 1}$ converges to $\bm{g}\opt \in \R^{\X}$.
\end{theorem}
\proof{Proof.}
When $\U$ is a definable uncertainty set, the same steps as Proposition \ref{prop:definable-value-function} show that the robust Bellman operator \eqref{eq:operator-bellman} is definable. We also recall that we can cast an sa-rectangular RMDP as a perfect information SG. The authors in \cite{bolte2015definable} show that when the operator of an SG is definable, then $\left(\frac{1}{t}\bm{v}^{t}\right)_{t \geq 1}$ admits a limit and that this limit coincides with $\lim_{\gamma \rightarrow 1} (1-\gamma)\bm{v}\opt_{\gamma}$ (see corollary 4, point (i), in \cite{bolte2015definable}).
    We have shown in the third step of the proof of Theorem \ref{th:sa-rec-relation-avg-rew-blackwell} that $\lim_{\gamma \rightarrow 1} (1-\gamma)\bm{p}_{0}\tr\bm{v}\opt_{\gamma}$ is equal to $\max_{\pi \in \Pi_{\sf SD}} \min_{\bm{P} \in \U} \Ravg(\pi,\bm{P})$. By definition of $\bm{g}\opt$, 
this shows that $\lim_{t \rightarrow + \infty} \frac{1}{t} \bm{v}^{t} = \bm{g}\opt$, which concludes the proof of Theorem \ref{th:convergence-alg-definable-uset-1}.
\hfill \halmos \endproof
The proof of Theorem \ref{th:convergence-alg-definable-uset-1} is relatively short. This conciseness comes from invoking powerful results from the stochastic game literature~\citep{bolte2015definable}, \tb{themselves related to the literature on {\em nonexpansive} operators $f$ and the limit of their {\em escape rate}, defined as $\left(f^{k}(\bm{0})\right)_{k \in \N}$~\citep{neyman2003stochastic,akian2003spectral,akian2016uniqueness}.}
\vspace{2mm}

\noindent
{\bf Algorithm based on increasing discount factor.}
\begin{algorithm}
\caption{Iterative algorithm for average optimality: increasing discount factor}\label{alg:sarec-increasing-discount-factor}
\begin{algorithmic}[1]
\State {\bf Input:} an increasing sequence $\left(\omega_{t}\right)_{t \geq 0} \in \R^{\N}_{+}$ with $\lim_{t \rightarrow + \infty} \omega_{t} = + \infty, \lim_{t \rightarrow + \infty} \omega_{t}/\omega_{t+1} = 1$.
\State Initialize $\bm{v}^{0} = \bm{0} \in \R^{\X}$.

\For{$t \in \N$}
    \State $\gamma_{t} \;=\; \frac{\omega_{t}}{\omega_{t+1}} $
    \State $v^{t+1}_{s} \;=\; \max_{a \in \A} \min_{\bm{p}_{sa} \in \U_{sa}} \bm{p}_{sa}\tr\left((1-\gamma_{t})\bm{r}_{sa} + \gamma_{t}\bm{v}^{t}\right), \forall \; s \in \X$. 
\EndFor
\end{algorithmic}
\end{algorithm}
We now consider Algorithm \ref{alg:sarec-increasing-discount-factor}. 
\tb{This algorithm computes iterative Bellman updates, but where the discount factor $\gamma_t$ increases at every iteration $t$ and converges to $1$. The main rationales behind this algorithms are the following: (a) we know  that if $\gamma_t$ is fixed to a certain value $\gamma \in [0,1)$ at every iteration, then Algorithm \ref{alg:sarec-increasing-discount-factor} would exactly coincide with value iteration for discounted RMDPs and it would converge to an optimal discounted value function $(1-\gamma)\bm{v}^{\star,\U}_{\gamma}$ (see Theorem 4 in \cite{wiesemann2013robust}), and (b) we know that when $\gamma_t \rightarrow 1$, the optimal optimal discounted value function converges to the optimal average gain:  $\lim_{\gamma \rightarrow 1} (1-\gamma)\bm{v}^{\star,\U}_{\gamma} = \bm{g}\opt$ (see our Lemma \ref{lem:aux-2}). To the best of our knowledge, Algorithm \ref{alg:sarec-increasing-discount-factor} was introduced for the first time in \cite{tewari2007bounded}, with weights $\omega_{t} = t+1$ and for sa-rectangular polyhedral uncertainty sets. Algorithm \ref{alg:sarec-increasing-discount-factor} has also been studied in \cite{wang2023robust} for unichain sa-rectangular uncertainty sets. We prove below its convergence for definable sa-rectangular uncertainty sets. In particular, we have the following theorem.
}
\begin{theorem}\label{th:convergence-alg-definable-uset-2}
    Let $\U$ be a definable compact sa-rectangular uncertainty set and let $\left(\bm{v}^{t}\right)_{t \geq 1}$ be the iterates of Algorithm \ref{alg:sarec-increasing-discount-factor}. Then $\left(\bm{v}^{t}\right)_{t \geq 1}$ converges to $\bm{g}\opt \in \R^{\X}$.
\end{theorem}
The proof of Theorem \ref{th:convergence-alg-definable-uset-2} is deferred to Appendix \ref{app:proof convergence alg uset 2}. We only provide the main steps here. The first step is to show the following lemma, which decomposes the suboptimality gap of the current iterates of Algorithm \ref{alg:sarec-limit-discount-values}.
\begin{lemma}\label{lem:alg uset 2 decomposition error}
There exist a period $t_{0} \in \N$ and two sequences of non-negative scalars $\left(\epsilon_{t}\right)_{t \geq t_{0}}$ and $\left(\delta_{t}\right)_{t \geq t_{0}}$ such that
    \[ \| \bm{v}^{t} - \bm{g}\opt \|_{\infty} \leq \epsilon_{t} + \delta_{t}, \forall \; t \geq t_{0},\]
with $\lim_{t \rightarrow 0} \epsilon_{t} \rightarrow 0$ and with, for $t \geq t_{0}$ and $\pi^{B}$ a Blackwell optimal policy, 
\begin{align}
    \delta_{t+1} & = \gamma_{t} \delta_{t} + e_{t}, \label{eq:delta-t-induction} \\
    e_{t} & = \| (1-\gamma_{t+1})\bm{v}^{\pi^B,\U}_{\gamma_{t+1}} - (1-\gamma_{t})\bm{v}^{\pi^B,\U}_{\gamma_{t}} \|_{\infty}. \label{eq: e-t definition}
\end{align}
\end{lemma}
The difficult part of the proof is to prove that the sequence $\left(\delta_{t}\right)_{t \geq t_{0}}$ converges to $0$. 
To show $\lim_{t \rightarrow + \infty}\delta_{t} \rightarrow 0$, we use the inductive step \eqref{eq:delta-t-induction} and $\gamma_{t} = \omega_{t}/\omega_{t+1} $ to show
\begin{equation}\label{eq:delta-t-unfolded}
    \delta_t = \frac{1}{\omega_{t+1}} \sum_{j=t_0}^{t} \omega_{j+1} e_{j}, \forall \; t \geq t_{0}.
\end{equation}
Equation \eqref{eq:delta-t-unfolded} reveals that the expression of $\delta_{t}$ resembles a weighted Cesaro average of the sequence $\left(e_{t}\right)_{t \geq t_{0}}$. It is known that the Cesaro average of a converging sequence also converges to the same limit; we prove Theorem \ref{th:convergence-alg-definable-uset-2} by extending this result to weighted Cesaro averages and by invoking existing results from SGs to obtain that $\lim_{t \rightarrow + \infty} e_{t} = 0.$ 

{\bf Comparison with previous work.}
We now compare Theorem \ref{th:convergence-alg-definable-uset-2} with related results from the RMDP literature. The authors in \cite{tewari2007bounded} prove
Theorem \ref{th:convergence-alg-definable-uset-2} with $\omega_{t} = t+1,$ i.e., with $\gamma_{t} = (t+1)/(t+2)$ and for the case of sa-rectangular uncertainty sets with $\ell_{\infty}$-balls. In this setting, \cite{tewari2007bounded} show the Lipschitz continuity of the robust value functions to simplify the term $e_{t}$ from Lemma \ref{lem:alg uset 2 decomposition error}. In particular, for sa-rectangular polyhedral $\U$, for $t$ large enough, there exists a constant $L>0$ such that
\begin{equation}\label{eq:Lip-e-t}
        e_{t} = \| (1-\gamma_{t+1})\bm{v}^{\pi^B,\U}_{\gamma_{t+1}} - (1-\gamma_{t})\bm{v}^{\pi^B,\U}_{\gamma_{t}} \|_{\infty} \leq L  \cdot \left(\gamma_{t+1}-\gamma_{t}\right).
    \end{equation}
    Once~\eqref{eq:Lip-e-t} is established, one can conclude that 
    \begin{equation}\label{eq:delta-t-induction-tewari}
        \delta_{t+1} = \gamma_t \delta_t + L \cdot \left(\gamma_{t+1}-\gamma_t\right)
    \end{equation}
    and Lemma 8 in \cite{tewari2007bounded} shows that a sequence $\left(\delta_t\right)_{t \geq 0}$ satisfying~\eqref{eq:delta-t-induction-tewari} converges to $0$. \tb{Our analysis of Algorithm \ref{alg:sarec-increasing-discount-factor} shows that the convergence holds for the much more general class of {\em definable} sa-rectangular uncertainty sets, and not solely for sa-rectangular uncertainty sets based on $\ell_{\infty}$-balls as proved in~\cite{tewari2007bounded}.}

    \tb{
    The authors in \cite{wang2023robust} consider Algorithm \ref{alg:sarec-increasing-discount-factor} with $\gamma_{t} = (t+1)/(t+2)$ for the case of general sa-rectangular uncertainty set (potentially non-polyhedral), with the assumption that the Markov chain associated with any transition probabilities and any policy is unichain. The analysis of Algorithm \ref{alg:sarec-increasing-discount-factor} in \cite{wang2023robust} follows the same lines as the proof in \cite{tewari2007bounded}, and the authors conclude by proving that the Lipschitzness property~\eqref{eq:Lip-e-t} holds for {\em unichain} sa-rectangular RMDPs. However, we do not assume that the RMDP is unichain, and in all generality the normalized robust value function $\gamma  \mapsto (1-\gamma)\bm{v}^{\pi,\U}_{\gamma}$ may not be Lipschitz continuous for a fixed $\pi \in \Pi_{\sf S}$, as we prove in Appendix \ref{app:proof counterexample not lip}. We sidestep this issue by leveraging existing results from SGs~\citep{bolte2015definable} and elementary results from real analysis.
    }

\tb{
\subsection{Convergence rates}
\subsubsection{Theoretical convergence rates}
We note that Theorem \ref{th:convergence-alg-definable-uset-0}, Theorem \ref{th:convergence-alg-definable-uset-1} and Theorem \ref{th:convergence-alg-definable-uset-2} shows the asymptotic convergence of Algorithm \ref{alg:sarec-limit-discount-values}, Algorithm \ref{alg:sarec-increasing-horizon} and Algorithm \ref{alg:sarec-increasing-discount-factor} to the optimal gain $\bm{g}\opt$. However, we are not able to provide a {\em theoretical convergence rate} for our algorithms, i.e. to bound the suboptimality of the vectors returned by our algorithms at every iteration. This situation is in contrast with the case of algorithms for the discounted return, for which precise convergence rates can be obtained e.g. for value iteration~\citep{wiesemann2013robust} and policy iteration~\citep{hansen2013strategy}. All the algorithms for computing discount optimal policies rely on two main properties:
\begin{itemize}
    \item Property 1: the optimal discounted value is the unique solution to a fixed-point optimality equation $\bm{v}\opt = T(\bm{v}\opt)$, with $T$ the Bellman operator as defined in \eqref{eq:operator-bellman}, 
    \item Property 2: the Bellman operator $T$ is a contraction of the rate $\gamma$ where $\gamma$ is the discount factor, that is, 
    \[ \| T(\bm{v}) - T(\bm{w}) \|_{\infty} \leq \gamma \| \bm{v} - \bm{w} \|_{\infty}, \forall \; \bm{v},\bm{w} \in \R^{\X}.\]
    This implies in particular that $\left(T^{t}(\bm{0})\right)_{t \in \N}$ converges to the fixed-point $\bm{v}\opt$ of $T$ at a rate of $\gamma^t$:
    \begin{equation}\label{eq:gamma-rate-vi}
        \| T^{t}(\bm{0}) - \bm{v}\opt \|_{\infty} \leq \gamma^t \| \bm{v}\opt\|_{\infty}, \forall \; t \geq 1.
    \end{equation}
\end{itemize}
In contrast, to the best of our knowledge, there are no known fixed-point optimality equations for the optimal gain $\bm{g}\opt$ in all generality. Worse, even though under some additional assumptions some fixed-point optimality equations can be constructed for $\bm{g}\opt$ (e.g. see Theorem 7 in \cite{wang2023robust} for unichain RMDPs), the optimality operator appearing in the fixed-point equations is not contracting. This absence of contraction property should not be surprising since we have repeatedly shown that average optimality can be interpreted as the limit of discount optimality when the discount factor $\gamma$ approaches $1$, in which case the bound \eqref{eq:gamma-rate-vi} becomes vacuous. The absence of Property 1 (in all generality) and Property 2 (even in some special cases) makes it difficult to provide theoretical convergence rates for our algorithms. We note that previous works studying variants of Algorithm \ref{alg:sarec-increasing-discount-factor} also do not provide theoretical convergence rate, and for perfect information mean payoff SGs, which is akin to sa-rectangular RMDPs with average return, there is no known {\em polynomial-time} algorithm for solving this class of SGs~\cite{condon1992complexity,zwick1996complexity,andersson2009complexity} even after more than six decades since their introduction in the game theory literature~\cite{gillette1957stochastic}, and this is regarded as one of the major open questions in algorithmic game theory.
}

\tb{
\subsubsection{Empirical convergence rates}
Given the difficulty of obtaining {\em theoretical} convergence rates for our algorithms, we study this question numerically. To do so, we run our algorithms for several different MDP instances and evaluate their errors as a function of the number of iterations. We describe our empirical setup in detail below.
}

\tb{
\vspace{2mm}
\noindent {\bf Test instances.}  We consider four different nominal MDP instances. In the machine replacement problem~\cite{delage2010percentile,wiesemann2013robust}, the goal is to compute a replacement and repair schedule for a line of machines. In the forest management instance~\cite{possingham1997application,pymdp}, a forest grows at every period, and the goal is to balance the revenue from wood cutting and the risk of wildfire. In the instance inspired from healthcare~\cite{goyal2022robust}, the goal is to plan the treatment of a patient, avoiding the mortality state while reducing the invasiveness of the treatment. Our last MDP instance consists of Garnet MDPs~\citep{archibald1995generation}, a widely used family of random MDP instances. We refer to the nominal transition probabilities in these instances as $\bm{P}_{\sf nom}$, and we represent them in Appendix \ref{app:simu-instance}. We instantiate all instances with $S=20$ states and a uniform initial distribution. We renormalize the rewards in all instances to $[0,1]$ to make the errors comparable.
}

\tb{
\vspace{2mm}
\noindent {\bf Construction of the uncertainty sets.} For each MDP instances, we consider two types of uncertainty sets. The first uncertainty set is based on box inequalities:
\begin{equation}\label{eq:uncertainty-set-box}
    \U_{sa}^{\sf box} = \{ \bm{p} \in \Delta(\X) \; | \; \bm{p}_{sa}^{\sf low} \leq \bm{p} \leq \bm{p}_{sa}^{\sf up} \}, \forall \; (s,a) \in \X \times \A
\end{equation}
with $\bm{p}_{sa}^{\sf low},\bm{p}_{sa}^{\sf up} \in [0,1]^{\X}$ two vectors such that $\bm{p}_{sa}^{\sf low} \leq \bm{p}_{{\sf nom},sa} \leq \bm{p}_{sa}^{\sf up}$ for $(s,a) \in \X \times \A$.
The second uncertainty set is based on ellipsoid constraints:
\begin{equation}\label{eq:uncertainty-set-ellipsoid}
\U_{sa}^{\sf \ell_2} = \{ \bm{p} \in \Delta(\X) \; | \; \| \bm{p} -\bm{p}_{{\sf nom},sa}\|_{2} \leq \alpha \}, \forall \; (s,a) \in \X \times \A
\end{equation}
where $\alpha>0$ is a scalar. Box inequalities have been used in applications of RMDPs in healthcare~\citep{goh2018data}, while uncertainty based on $\ell_{2}$-distance can be seen as conservative approximations of sets based on relative entropy~\citep{iyengar2005robust}.
Note that both $\U^{\ell_{2}}$ and $\U^{\sf box}$ are definable sets, with $\U^{\sf box}$ being polyhedral. Additionally, there exist efficient algorithms to evaluate the robust Bellman operator as in~\eqref{eq:operator-bellman} for $\U^{\sf box}$ (e.g. proposition 3 in \cite{goh2018data}). For $\U^{\sf \ell_2}$, evaluating the robust Bellman operator requires solving a convex program. Still, we can obtain a closed-form expression assuming that $\alpha>0$ is sufficiently small, as we detail in Appendix \ref{app:simu-uncertainty sets}. 
For $\U^{\sf \ell_2}$ we choose a radius of $\alpha = 0.05$ and for $\U^{\sf box}$ we choose the upper and lower bound on each coefficient to allow for $5 \%$ deviations from the nominal distributions.
}

\tb{
\vspace{2mm}
\noindent {\bf Practical implementation.}
We implement all algorithms in Python 3.8.8.
For Algorithm \ref{alg:sarec-limit-discount-values} we choose a sequence of discount factors $\gamma_t = (t+1)/(t+2)$. To solve a robust MDP at every iteration in Algorithm \ref{alg:sarec-limit-discount-values}, we implement the two-player strategy iteration~\cite{hansen2013strategy} with warm-starts, see more details in Appendix \ref{app:simu-2pl-PI}. 
For Algorithm \ref{alg:sarec-increasing-discount-factor}, we choose $\omega_{t} = t+1$, so that the resulting sequence of discount factors in Algorithm \ref{alg:sarec-increasing-discount-factor} coincides with the sequence of discount factors in Algorithm \ref{alg:sarec-limit-discount-values}.  
}

\tb{
\vspace{2mm}
\noindent {\bf Empirical setup.} Our goal is to compute the the difference between the vectors returned by the algorithms at each iteration and the optimal gain. To do so, we first estimate the optimal gain by running Algorithm \ref{alg:sarec-limit-discount-values} for a large number of iterations, which we then compare to the iterates of the algorithms. In particular, we proceed as follows. 
\begin{enumerate}
    \item We first run Algorithm \ref{alg:sarec-limit-discount-values} for $\hat{T} = 5 \times 1,000$ iterations to obtain an estimate $\hat{\bm{g}}= (1-\gamma_{\hat{T}})\bm{v}_{\gamma_{\hat{T}}}^{\star,\U}$ of the optimal gain $\bm{g}\opt$, and therefore an estimate $\bm{p}_{0}\tr\hat{\bm{g}}$ of the optimal average return $\bm{p}_{0}\tr\bm{g}\opt$. We verify that for the last iterations close to $\hat{T}$, the iterates of Algorithm \ref{alg:sarec-limit-discount-values} vary only marginally by at most $10^{-5}$.
    \item We then run each algorithm for $T=1,000$ iterations. At every iteration $t$, the reported error $\err_{1,t}$ for Algorithm \ref{alg:sarec-limit-discount-values}, the reported error $\err_{2,t}$ for Algorithm \ref{alg:sarec-increasing-horizon} and the reported error $\err_{3,t}$ for Algorithm \ref{alg:sarec-increasing-discount-factor} are
    \begin{align*}
        \err_{1,t} & = |\bm{p}_{0}\tr\hat{\bm{g}} - (1-\gamma_t)\bm{p}_{0}\tr\bm{v}_{\gamma_{t}}^{\star,\U}|,\\
        \err_{2,t} & = |\bm{p}_{0}\tr\hat{\bm{g}}-(1/t)\bm{p}_{0}\tr\bm{v}^{t}|, \\
        \err_{3,t} & = |\bm{p}_{0}\tr\hat{\bm{g}}-\bm{p}_{0}\tr\bm{v}^{t}|.
    \end{align*}
\end{enumerate}
For Garnet MDPs, we display the mean errors over $25$ samples of the nominal transition probabilities as well as the associated 95 \% confidence intervals.
}

\begin{figure}[htb]
\begin{center}
\includegraphics[width=\linewidth]{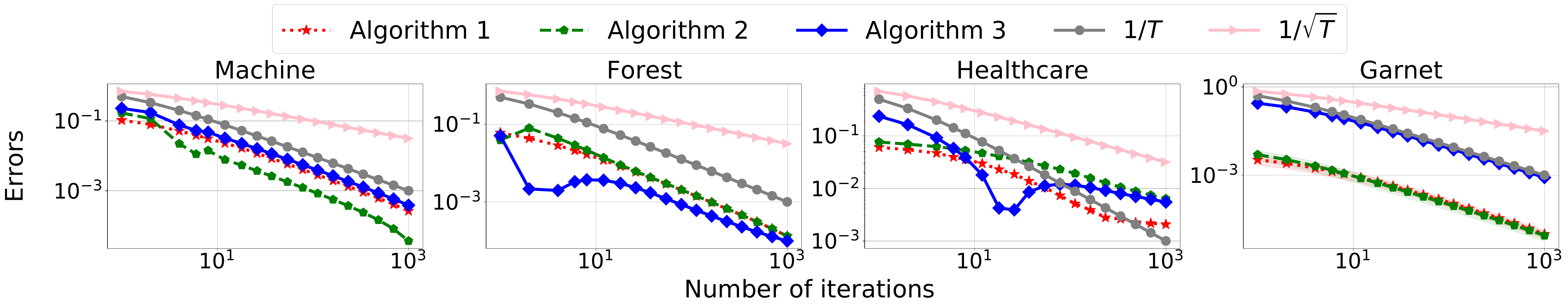}
\end{center}
\caption{Errors to the optimal gain for a box uncertainty set~\eqref{eq:uncertainty-set-box}.}
\label{fig:sim:hypercube_error_to_optimal_gain}
\end{figure}

\begin{figure}[htb]
\begin{center}
\includegraphics[width=\linewidth]{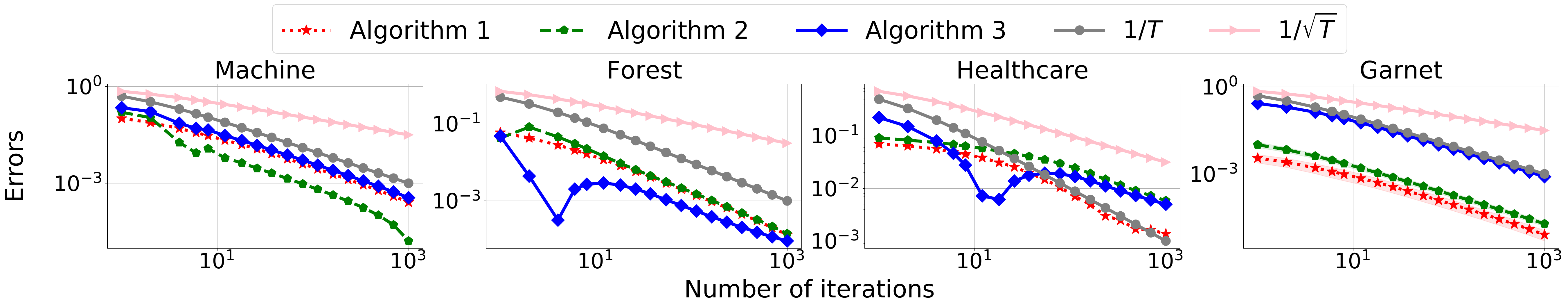}
\end{center}
\caption{Errors to the optimal gain for an ellipsoid uncertainty set~\eqref{eq:uncertainty-set-ellipsoid}.}
\label{fig:sim:ellipsoid_error_to_optimal_gain}
\end{figure}

\tb{
\vspace{2mm}
\noindent {\bf Numerical results.} Our numerical results for the errors of our algorithms at every iteration are presented in Figure \ref{fig:sim:hypercube_error_to_optimal_gain} (with box uncertainty set~\eqref{eq:uncertainty-set-box}) and Figure \ref{fig:sim:ellipsoid_error_to_optimal_gain} (with ellipsoid uncertainty set~\eqref{eq:uncertainty-set-ellipsoid}). We use log-log plots to make it easier to read the empirical convergence rate directly on the figures. For reference, we also display the lines corresponding to $1/T$ and $1/\sqrt{T}$, where $T$ is the number of iterations. All our algorithms exhibit a $O(1/T)$ empirical convergence rate in the Machine, Forest, and Garnet MDP instances, with no algorithms clearly outperforming the others across all instances and at every iteration. The results in the Healthcare instance are more surprising: the empirical convergence rate of Algorithm \ref{alg:sarec-increasing-horizon} and Algorithm \ref{alg:sarec-increasing-discount-factor} appears close to $O(1/\sqrt{T})$ while the empirical convergence rate of Algorithm \ref{alg:sarec-limit-discount-values} remains close to $O(1/T)$. Additionally, the errors may not decrease monotonically, as evident for Algorithm \ref{alg:sarec-increasing-discount-factor} on the Forest and Healthcare MDP instances.
}


\tb{
\section{Discussion and main insights}\label{sec:discussion}
We now discuss the main insights from our paper. In particular, the results in our paper highlight several pitfalls and subtle differences between discounted return and average return and emphasize the superior practical properties of sa-rectangular models based on natural distance-based uncertainty sets compared to s-rectangular uncertainty.
Below, we provide a more detailed discussion of our methodological takeaways.
}

\noindent
\tb{{\bf On the interplay between discounted and average return.}
The interplay between average and discounted returns for large discount factors has been widely studied in the nominal MDP literature, but much less in the robust MDP literature. Our paper provides several important results in this direction, both positive and negative.
}

\tb{
On the positive side, for sa-rectangular uncertainty sets, the worst-case average {\em values} are indeed limits of worst-case discounted {\em values} (Lemma \ref{lem:aux-1}):
    \[ \lim_{\gamma \rightarrow 1} \min_{\bm{P} \in \U} (1-\gamma) R_{\gamma}(\pi,\bm{P}) = \inf_{\bm{P} \in \U} \Ravg(\pi,\bm{P}), \forall \; \pi \in \PiS\]
and optimal worst-case average {\em values} are limits of optimal worst-case discounted {\em values} (Lemma \ref{lem:aux-2}):
    \[ \lim_{\gamma \rightarrow 1} \max_{\pi \in \Pi_{\sf SD}} \min_{\bm{P} \in \U} (1-\gamma)R_{\gamma}(\pi,\bm{P}) = \max_{\pi \in \Pi_{\sf SD}} \inf_{\bm{P} \in \U} \Ravg(\pi,\bm{P}).\]
    The same limit behavior is true for {\em policies} for sa-rectangular sets: discount optimal policies become average optimal for large enough discount factors (Theorem \ref{th:sa-rec-relation-avg-rew-discount}), and average optimal policies are $\epsilon$-Blackwell optimal for all $\epsilon>0$ (Theorem \ref{th:sa-rec-relation-avg-rew-blackwell}).
}

\tb{
On the negative side, we note that this limit behavior may fail for the worst-case {\em transition probabilities}, even for sa-rectangular uncertainty. Our results show that one should be very careful about the interpretation of the notion of {\em worst-case transition probabilities} when using average return, in contrast to the discounted return, where this notion is always well-defined.
Indeed, Proposition \ref{prop:inf-is-not-min-avg-rew} shows an example where the worst-case average value is not attained and where the limit of the sequence of worst-case transition probabilities for the discounted return is meaningless, and similarly for a sequence of transition probabilities achieving the worst-case average return in the limit. 
Even worse, for s-rectangular uncertainty, the discount optimal policies may be arbitrarily bad for the average return, as highlighted in the example of \tbm, where a discount optimal policy always achieves a worst-case average return of $0$ independently of the discount factor, whereas the optimal worst-case average return is $1/2$ (attained by a Markovian policy).
}

\tb{
Overall, our results provide a cautionary tale about the interplay between discounted and average returns. While for sa-rectangular uncertainty it is a correct intuition to view optimal average values and policies as limits of optimal discounted values and policies, this intuition is wrong for the worst-case transition probabilities, and for s-rectangular uncertainty, one cannot (a priori) rely on the discounted problem to approach the average problem. 
}

\noindent
\tb{{\bf The case for sa-rectangular uncertainty.} Recall that for the discounted return, s-rectangular and sa-rectangular models of uncertainty are ``mostly equivalent": both models guarantee the existence of stationary discount optimal policies that can be computed efficiently, the main difference being that discount optimal policies may be chosen deterministic for sa-rectangular models, whereas discount optimal policies may have to be randomized for s-rectangular models.
}

\tb{Contrary to the discounted case, our paper shows that there is a contrast between s-rectangularity and sa-rectangularity for the average return. In particular, only the sa-rectangular model guarantees the existence of average optimal policies, their stationarity, and their deterministic nature (Theorem \ref{th:avg-reward-main}). In contrast, for s-rectangular models, there is no guarantee that an average optimal policy exists (Proposition \ref{prop:avg-s-rec no existence opt policies}), and even in very simple problem instances where they exist, average optimal policies may have to be non-stationary (Proposition \ref{prop:tbm-avg-rew}), making them more difficult to deploy and interpret in real-world applications.}
\tb{ The sa-rectangular model also ensures that stationary adversaries are equivalent to history-dependent adversaries (Theorem \ref{th:avg-rew-history-dep-adversaries}) and that various notions of average return are equivalent (Corollary \ref{cor:avg-reward-main-other-ravg}). For all these reasons, we recommend choosing sa-rectangular models for practical applications of RMDPs with the average return.}

\noindent
\tb{{\bf The case for distance-based uncertainty sets.} Our analysis from Section \ref{sec:definable-sa-rec} shows that virtually all distance-based uncertainty sets, e.g., based on Kullback-Leibler divergence, $\ell_{2}$-distance or Wasserstein distance, fall into the class of {\em definable} sets. Our main results show that definable uncertainty sets imply several important properties, including the existence of Blackwell optimal policies (Theorem \ref{th:definable-sa-rec-blackwell-opt}), the existence of a Blackwell discount factor above which all discount optimal policies are also Blackwell optimal (Theorem \ref{th:blackwell-discount-factor}), and well-behaved (non-oscillatory) behaviors for the discounted value functions. For definable uncertainty sets, we are also able to provide algorithms to compute the optimal average value with small per-iteration complexity (Algorithm \ref{alg:sarec-increasing-horizon} and Algorithm \ref{alg:sarec-increasing-discount-factor}) compared to more general algorithm (Algorithm \ref{alg:sarec-limit-discount-values}). Note that distance-based uncertainty sets are the most used in practice since they naturally yield statistical guarantees~\citep{ramani2022robust} (and, in fact, we are not aware of any other types of uncertainty sets being used for RMDPs), and one of the main take-aways from our paper is that they can safely be used in the average return setting (combined with sa-rectangularity).}
\section{Conclusion}
 Our paper addresses several important issues in the existing literature and derives the fundamental properties of average optimal and Blackwell optimal policies. In particular, our work highlights important distinctions between the widely-studied framework of discounted RMDPs and the less-studied frameworks of RMDPs with average optimality and Blackwell optimality. We view the non-existence of average optimal policies for s-rectangular RMDPs and the non-existence of Blackwell optimal policies for sa-rectangular RMDPs (in all generality) as surprising results. The notion of {\em definable} uncertainty sets is also of independent interest, characterizing the cases where the value functions are well-behaved and iterative algorithms asymptotically converge to the optimal gain in the sa-rectangular case. 
 \tb{Our main results suggest that decision-makers should use distance-based sa-rectangular models of uncertainty when focusing on the average return, a model where one can (mostly) rely on the intuition built from the discounted case.}
 Finally, our work opens new research avenues for RMDPs. Among them, deriving an efficient algorithm for computing an average optimal policy for sa-rectangular RMDPs appears crucial, but this may be difficult since a similar question remains open after several decades in the literature on SGs. Important other research questions include studying in more detail the particular case of irreducible instances, which considerably simplifies the main technical challenges, and studying the case of distributionally robust MDPs. \tb{It would also be interesting to derive the optimality equations for the optimal gain in all generality, akin to the Bellman equation for the optimal discounted value functions.}

\bibliographystyle{alpha}
\bibliography{ref}

\newcommand{\etalchar}[1]{$^{#1}$}
\begin{thebibliography}{EPGMFBV03}

\bibitem[AG03]{akian2003spectral}
Marianne Akian and St{\'e}phane Gaubert.
\newblock Spectral theorem for convex monotone homogeneous maps, and ergodic
  control.
\newblock {\em Nonlinear Analysis: Theory, Methods \& Applications},
  52(2):637--679, 2003.

\bibitem[AGGCG19]{akian2019operator}
Marianne Akian, St{\'e}phane Gaubert, Julien Grand-Cl{\'e}ment, and
  J{\'e}r{\'e}mie Guillaud.
\newblock The operator approach to entropy games.
\newblock {\em Theory of Computing Systems}, 63(5):1089--1130, 2019.

\bibitem[AGN16]{akian2016uniqueness}
Marianne Akian, St{\'e}phane Gaubert, and Roger Nussbaum.
\newblock Uniqueness of the fixed point of nonexpansive semidifferentiable
  maps.
\newblock {\em Transactions of the American Mathematical Society},
  368(2):1271--1320, 2016.

\bibitem[AM09]{andersson2009complexity}
Daniel Andersson and Peter~Bro Miltersen.
\newblock The complexity of solving stochastic games on graphs.
\newblock In {\em International Symposium on Algorithms and Computation}, pages
  112--121. Springer, 2009.

\bibitem[AMT95]{archibald1995generation}
TW~Archibald, KIM McKinnon, and LC~Thomas.
\newblock On the generation of markov decision processes.
\newblock {\em Journal of the Operational Research Society}, 46(3):354--361,
  1995.

\bibitem[BB01]{baxter2001infinite}
Jonathan Baxter and Peter~L Bartlett.
\newblock Infinite-horizon policy-gradient estimation.
\newblock {\em journal of artificial intelligence research}, 15:319--350, 2001.

\bibitem[BCP{\etalchar{+}}16]{brockman2016openai}
Greg Brockman, Vicki Cheung, Ludwig Pettersson, Jonas Schneider, John Schulman,
  Jie Tang, and Wojciech Zaremba.
\newblock Openai gym.
\newblock {\em arXiv preprint arXiv:1606.01540}, 2016.

\bibitem[BF68]{blackwell1968big}
David Blackwell and Tom~S Ferguson.
\newblock The big match.
\newblock {\em The Annals of Mathematical Statistics}, 39(1):159--163, 1968.

\bibitem[BGV15]{bolte2015definable}
J{\'e}r{\^o}me Bolte, St{\'e}phane Gaubert, and Guillaume Vigeral.
\newblock Definable zero-sum stochastic games.
\newblock {\em Mathematics of Operations Research}, 40(1):171--191, 2015.

\bibitem[BH13]{bennett2013artificial}
Casey~C Bennett and Kris Hauser.
\newblock Artificial intelligence framework for simulating clinical
  decision-making: A {M}arkov decision process approach.
\newblock {\em Artificial intelligence in medicine}, 57(1):9--19, 2013.

\bibitem[Bie87]{bierth1987expected}
K-J Bierth.
\newblock An expected average reward criterion.
\newblock {\em Stochastic processes and their applications}, 26:123--140, 1987.

\bibitem[BK76]{bewley1976asymptotic}
Truman Bewley and Elon Kohlberg.
\newblock The asymptotic theory of stochastic games.
\newblock {\em Mathematics of Operations Research}, 1(3):197--208, 1976.

\bibitem[BPH21]{behzadian2021fast}
Bahram Behzadian, Marek Petrik, and Chin~Pang Ho.
\newblock Fast algorithms for $l_{\infty}$ constrained s-rectangular robust
  {MDP}s.
\newblock {\em Advances in Neural Information Processing Systems},
  34:25982--25992, 2021.

\bibitem[BR11]{bauerle2011markov}
Nicole B{\"a}uerle and Ulrich Rieder.
\newblock {\em Markov decision processes with applications to finance}.
\newblock Springer Science \& Business Media, 2011.

\bibitem[CGT15]{pymdp}
Steven Cordwell, Yasser Gonzalez, and Theja Tulabandhula.
\newblock {M}arkov {D}ecision {P}rocess ({MDP}) toolbox for python.
\newblock {\em https://github.com/sawcordwell/pymdptoolbox}, 2015.

\bibitem[CHJ{\etalchar{+}}22]{chae2022robust}
Jongseong Chae, Seungyul Han, Whiyoung Jung, Myungsik Cho, Sungho Choi, and
  Youngchul Sung.
\newblock Robust imitation learning against variations in environment dynamics.
\newblock In {\em International Conference on Machine Learning}, pages
  2828--2852. PMLR, 2022.

\bibitem[CN14]{choudary2014real}
Alla Dita~Raza Choudary and Constantin~P Niculescu.
\newblock {\em Real analysis on intervals}.
\newblock Springer, 2014.

\bibitem[Con92]{condon1992complexity}
Anne Condon.
\newblock The complexity of stochastic games.
\newblock {\em Information and Computation}, 96(2):203--224, 1992.

\bibitem[Cos00]{coste2000introduction}
Michel Coste.
\newblock {\em An introduction to o-minimal geometry}.
\newblock Istituti editoriali e poligrafici internazionali Pisa, 2000.

\bibitem[DBK{\etalchar{+}}16]{deng2016deep}
Yue Deng, Feng Bao, Youyong Kong, Zhiquan Ren, and Qionghai Dai.
\newblock Deep direct reinforcement learning for financial signal
  representation and trading.
\newblock {\em IEEE transactions on neural networks and learning systems},
  28(3):653--664, 2016.

\bibitem[DDE{\etalchar{+}}20]{dewanto2020average}
Vektor Dewanto, George Dunn, Ali Eshragh, Marcus Gallagher, and Fred Roosta.
\newblock Average-reward model-free reinforcement learning: a systematic review
  and literature mapping.
\newblock {\em arXiv preprint arXiv:2010.08920}, 2020.

\bibitem[DM10]{delage2010percentile}
Erick Delage and Shie Mannor.
\newblock Percentile optimization for {Markov} decision processes with
  parameter uncertainty.
\newblock {\em Operations research}, 58(1):203--213, 2010.

\bibitem[EPGMFBV03]{egozcue2003isometric}
Juan~Jos{\'e} Egozcue, Vera Pawlowsky-Glahn, Gl{\`o}ria Mateu-Figueras, and
  Carles Barcelo-Vidal.
\newblock Isometric logratio transformations for compositional data analysis.
\newblock {\em Mathematical geology}, 35(3):279--300, 2003.

\bibitem[Eve57]{everett1957recursive}
Hugh Everett.
\newblock Recursive games.
\newblock {\em Contributions to the Theory of Games}, 3(39):47--78, 1957.

\bibitem[FS12]{feinberg2012handbook}
Eugene~A Feinberg and Adam Shwartz.
\newblock {\em Handbook of {Markov} decision processes: methods and
  applications}, volume~40.
\newblock Springer Science \& Business Media, 2012.

\bibitem[FV12]{filar2012competitive}
Jerzy Filar and Koos Vrieze.
\newblock {\em Competitive {Markov} decision processes}.
\newblock Springer Science \& Business Media, 2012.

\bibitem[GBZ{\etalchar{+}}18]{goh2018data}
Joel Goh, Mohsen Bayati, Stefanos~A Zenios, Sundeep Singh, and David Moore.
\newblock Data uncertainty in {M}arkov chains: Application to
  cost-effectiveness analyses of medical innovations.
\newblock {\em Operations Research}, 66(3):697--715, 2018.

\bibitem[GCCGE23]{grand2023robustness}
Julien Grand-Cl{\'e}ment, Carri~W Chan, Vineet Goyal, and Gabriel Escobar.
\newblock Robustness of proactive intensive care unit transfer policies.
\newblock {\em Operations Research}, 71(5):1653--1688, 2023.

\bibitem[GCK21a]{grand2021conic}
Julien Grand-Cl{\'e}ment and Christian Kroer.
\newblock Conic blackwell algorithm: Parameter-free convex-concave saddle-point
  solving.
\newblock {\em Advances in Neural Information Processing Systems},
  34:9587--9599, 2021.

\bibitem[GCK21b]{grand2021first}
Julien Grand-Clement and Christian Kroer.
\newblock First-order methods for wasserstein distributionally robust {MDP}.
\newblock In {\em International Conference on Machine Learning}, pages
  2010--2019. PMLR, 2021.

\bibitem[GCK23]{grand2023solving}
Julien Grand-Cl{\'e}ment and Christian Kroer.
\newblock Solving optimization problems with blackwell approachability.
\newblock {\em Mathematics of Operations Research}, 2023.

\bibitem[GCP24a]{grand2024convex}
Julien Grand-Cl{\'e}ment and Marek Petrik.
\newblock On the convex formulations of robust {Markov} decision processes.
\newblock {\em Mathematics of Operations Research}, 2024.

\bibitem[GCP24b]{grand2024reducing}
Julien Grand-Cl{\'e}ment and Marek Petrik.
\newblock Reducing blackwell and average optimality to discounted mdps via the
  blackwell discount factor.
\newblock {\em Advances in Neural Information Processing Systems}, 36, 2024.

\bibitem[GGC22]{goyal2022robust}
Vineet Goyal and Julien Grand-Cl{\'e}ment.
\newblock Robust {M}arkov decision processes: Beyond rectangularity.
\newblock {\em Mathematics of Operations Research}, 2022.

\bibitem[Gil57]{gillette1957stochastic}
Dean Gillette.
\newblock Stochastic games with zero stop probabilities.
\newblock {\em Contributions to the Theory of Games}, 3:179--187, 1957.

\bibitem[GLD97]{givan1997bounded}
Robert Givan, Sonia Leach, and Thomas Dean.
\newblock Bounded parameter {M}arkov decision processes.
\newblock In {\em European Conference on Planning}, pages 234--246. Springer,
  1997.

\bibitem[HMZ13]{hansen2013strategy}
Thomas~Dueholm Hansen, Peter~Bro Miltersen, and Uri Zwick.
\newblock Strategy iteration is strongly polynomial for 2-player turn-based
  stochastic games with a constant discount factor.
\newblock {\em Journal of the ACM (JACM)}, 60(1):1--16, 2013.

\bibitem[HPW21]{ho2021partial}
Chin~Pang Ho, Marek Petrik, and Wolfram Wiesemann.
\newblock Partial policy iteration for l1-robust {Markov} decision processes.
\newblock {\em The Journal of Machine Learning Research}, 22(1):12612--12657,
  2021.

\bibitem[HPW22]{ho2022robust}
Chin~Pang Ho, Marek Petrik, and Wolfram Wiesemann.
\newblock Robust phi-divergence mdps.
\newblock {\em Advances in Neural Information Processing Systems},
  35:32680--32693, 2022.

\bibitem[Iye05]{iyengar2005robust}
G.~Iyengar.
\newblock Robust dynamic programming.
\newblock {\em Mathematics of Operations Research}, 30(2):257--280, 2005.

\bibitem[KDG{\etalchar{+}}23]{kumar2023policy}
Navdeep Kumar, Esther Derman, Matthieu Geist, Kfir~Y Levy, and Shie Mannor.
\newblock Policy gradient for rectangular robust markov decision processes.
\newblock {\em Advances in Neural Information Processing Systems},
  36:59477--59501, 2023.

\bibitem[Koh74]{kohlberg1974repeated}
Elon Kohlberg.
\newblock Repeated games with absorbing states.
\newblock {\em The Annals of Statistics}, pages 724--738, 1974.

\bibitem[Lei03]{leizarowitz2003algorithm}
Arie Leizarowitz.
\newblock An algorithm to identify and compute average optimal policies in
  multichain {Markov} decision processes.
\newblock {\em Mathematics of Operations Research}, 28(3):553--586, 2003.

\bibitem[LL69]{liggett1969stochastic}
Thomas~M Liggett and Steven~A Lippman.
\newblock Stochastic games with perfect information and time average payoff.
\newblock {\em Siam Review}, 11(4):604--607, 1969.

\bibitem[LS15]{laraki2015advances}
Rida Laraki and Sylvain Sorin.
\newblock Advances in zero-sum dynamic games.
\newblock In {\em Handbook of game theory with economic applications},
  volume~4, pages 27--93. Elsevier, 2015.

\bibitem[LSK23]{li2023policy}
Mengmeng Li, Tobias Sutter, and Daniel Kuhn.
\newblock Policy gradient algorithms for robust {MDP}s with non-rectangular
  uncertainty sets.
\newblock {\em arXiv preprint arXiv:2305.19004}, 2023.

\bibitem[LT07]{le2007robust}
Yann Le~Tallec.
\newblock {\em Robust, risk-sensitive, and data-driven control of {Markov}
  decision processes}.
\newblock PhD thesis, Massachusetts Institute of Technology, 2007.

\bibitem[LZL22]{li2022first}
Yan Li, Tuo Zhao, and Guanghui Lan.
\newblock First-order policy optimization for robust {Markov} decision process.
\newblock {\em arXiv preprint arXiv:2209.10579}, 2022.

\bibitem[MKS{\etalchar{+}}13]{mnih2013playing}
Volodymyr Mnih, Koray Kavukcuoglu, David Silver, Alex Graves, Ioannis
  Antonoglou, Daan Wierstra, and Martin Riedmiller.
\newblock Playing atari with deep reinforcement learning.
\newblock {\em arXiv preprint arXiv:1312.5602}, 2013.

\bibitem[MMX16]{mannor2016robust}
S.~Mannor, O.~Mebel, and H.~Xu.
\newblock Robust {MDP}s with k-rectangular uncertainty.
\newblock {\em Mathematics of Operations Research}, 41(4):1484--1509, 2016.

\bibitem[MSZ15]{mertens2015repeated}
Jean-Fran{\c{c}}ois Mertens, Sylvain Sorin, and Shmuel Zamir.
\newblock {\em Repeated games}, volume~55.
\newblock Cambridge University Press, 2015.

\bibitem[NG05]{nilim2005robust}
A.~Nilim and L.~El Ghaoui.
\newblock Robust control of {M}arkov decision processes with uncertain
  transition probabilities.
\newblock {\em Operations Research}, 53(5):780--798, 2005.

\bibitem[NSS03]{neyman2003stochastic}
Abraham Neyman, Sylvain Sorin, and S~Sorin.
\newblock {\em Stochastic games and applications}, volume 570.
\newblock Springer Science \& Business Media, 2003.

\bibitem[PK22]{panaganti2022sample}
Kishan Panaganti and Dileep Kalathil.
\newblock Sample complexity of robust reinforcement learning with a generative
  model.
\newblock In {\em International Conference on Artificial Intelligence and
  Statistics}, pages 9582--9602. PMLR, 2022.

\bibitem[PT97]{possingham1997application}
Hugh Possingham and G~Tuck.
\newblock Application of stochastic dynamic programming to optimal fire
  management of a spatially structured threatened species.
\newblock In {\em Proceedings International Congress on Modelling and
  Simulation, MODSIM}, pages 813--817, 1997.

\bibitem[Put14]{puterman2014markov}
Martin~L Puterman.
\newblock {\em Markov decision processes: discrete stochastic dynamic
  programming}.
\newblock John Wiley \& Sons, 2014.

\bibitem[Ren19]{renault2019tutorial}
J{\'e}r{\^o}me Renault.
\newblock A tutorial on zero-sum stochastic games.
\newblock {\em arXiv preprint arXiv:1905.06577}, 2019.

\bibitem[RG22]{ramani2022robust}
Sivaramakrishnan Ramani and Archis Ghate.
\newblock Robust markov decision processes with data-driven, distance-based
  ambiguity sets.
\newblock {\em SIAM Journal on Optimization}, 32(2):989--1017, 2022.

\bibitem[RG24]{ramani2024family}
Sivaramakrishnan Ramani and Archis Ghate.
\newblock A family of-rectangular robust mdps: Relative conservativeness,
  asymptotic analyses, and finite-sample properties.
\newblock {\em SIAM Journal on Optimization}, 34(2):1540--1568, 2024.

\bibitem[SB18]{sutton2018reinforcement}
Richard~S Sutton and Andrew~G Barto.
\newblock {\em Reinforcement learning: An introduction}.
\newblock MIT press, 2018.

\bibitem[SD17]{steimle2017markov}
Lauren~N Steimle and Brian~T Denton.
\newblock Markov decision processes for screening and treatment of chronic
  diseases.
\newblock {\em Markov Decision Processes in Practice}, pages 189--222, 2017.

\bibitem[Sha53]{shapley1953stochastic}
Lloyd~S Shapley.
\newblock Stochastic games.
\newblock {\em Proceedings of the national academy of sciences},
  39(10):1095--1100, 1953.

\bibitem[SL73]{satia1973markov}
J.K. Satia and R.L. Lave.
\newblock {M}arkov decision processes with uncertain transition probabilities.
\newblock {\em Operations Research}, 21(3):728--740, 1973.

\bibitem[Sor02]{sorin2002first}
Sylvain Sorin.
\newblock {\em A first course on zero-sum repeated games}, volume~37.
\newblock Springer Science \& Business Media, 2002.

\bibitem[TB07]{tewari2007bounded}
Ambuj Tewari and Peter~L Bartlett.
\newblock Bounded parameter {M}arkov decision processes with average reward
  criterion.
\newblock In {\em International Conference on Computational Learning Theory},
  pages 263--277. Springer, 2007.

\bibitem[VDD98]{van1998minimal}
Lou Van Den~Dries.
\newblock O-minimal structures and real analytic geometry.
\newblock {\em Current developments in mathematics}, 1998(1):105--152, 1998.

\bibitem[VdDM96]{van1996geometric}
Lou Van~den Dries and Chris Miller.
\newblock Geometric categories and o-minimal structures.
\newblock 1996.

\bibitem[VHK{\etalchar{+}}21]{viano2021robust}
Luca Viano, Yu-Ting Huang, Parameswaran Kamalaruban, Adrian Weller, and Volkan
  Cevher.
\newblock Robust inverse reinforcement learning under transition dynamics
  mismatch.
\newblock {\em Advances in Neural Information Processing Systems},
  34:25917--25931, 2021.

\bibitem[Vil21]{villani2021topics}
C{\'e}dric Villani.
\newblock {\em Topics in optimal transportation}, volume~58.
\newblock American Mathematical Soc., 2021.

\bibitem[WHP23]{wang2022convergence}
Qiuhao Wang, Chin~Pang Ho, and Marek Petrik.
\newblock Policy gradient in robust mdps with global convergence guarantee.
\newblock In {\em ICML}, volume 202 of {\em Proceedings of Machine Learning
  Research}, pages 35763--35797. PMLR, 2023.

\bibitem[WKR13]{wiesemann2013robust}
Wolfram Wiesemann, Daniel Kuhn, and Ber{\c{c}} Rustem.
\newblock Robust {M}arkov decision processes.
\newblock {\em Mathematics of Operations Research}, 38(1):153--183, 2013.

\bibitem[WVA{\etalchar{+}}23]{wang2023robust}
Yue Wang, Alvaro Velasquez, George Atia, Ashley Prater-Bennette, and Shaofeng
  Zou.
\newblock Robust average-reward markov decision processes.
\newblock In {\em Proceedings of the AAAI Conference on Artificial
  Intelligence}, volume~37, pages 15215--15223, 2023.

\bibitem[XM10]{xu2010distributionally}
Huan Xu and Shie Mannor.
\newblock Distributionally robust {Markov} decision processes.
\newblock {\em Advances in Neural Information Processing Systems}, 23, 2010.

\bibitem[ZP96]{zwick1996complexity}
Uri Zwick and Mike Paterson.
\newblock The complexity of mean payoff games on graphs.
\newblock {\em Theoretical Computer Science}, 158(1-2):343--359, 1996.

\bibitem[ZSD17]{zhang2017robust}
Yuanhui Zhang, L~Steimle, and BT~Denton.
\newblock Robust {M}arkov decision processes for medical treatment decisions.
\newblock {\em Optimization online}, 2017.

\end{thebibliography}



\begin{APPENDICES}

\section{Comparisons of SGs and RMDPs}\label{app:more-detailed-sgs-rmdps}
\tb{In this section, we show the connections between stochastic games and rectangular robust MDPs. 
}
\tb{
\subsection{Perfect information SGs and sa-rectangularity}
\label{app:more-detailed-sgs-rmdps-sarec}
}
\tb{
The class of perfect information stochastic games corresponds to SG instances where the state set $\Omega$ can be partitioned into two subsets $\Omega_{1},\Omega_{2}$ such that the first player entirely controls the states in $\Omega_{1}$ (in terms of instantaneous rewards and probabilities), and similarly for the second player for states in $\Omega_{2}$. We refer to chapter 4 in \cite{neyman2003stochastic} for a concise introduction to perfect information SGs.}

\tb{
We now briefly describe the reformulation of an sa-rectangular RMDP as a perfect information SG. We note that this construction is also described in section 5.2 in \cite{goyal2022robust} (without making the connection to stochastic games), in appendix A in \cite{nilim2005robust} and in Section 5 in \cite{iyengar2005robust}. We provide the reformulation here for completeness.
}
\tb{
Consider an sa-rectangular RMDP instance $\left(\X,\A,\bm{r},\U,\bm{p}_{0}\right)$ with a compact set $\U$. We define a corresponding perfect information stochastic game instance as follows. The set of states is $\Omega = \Omega_{1} \cup \Omega_{2}$ with $\Omega_{1} = \X,\Omega_{2} = \X \times \A$. The set of actions for the first player is $\A$, and the set of actions for the second player is the set of extreme points of $\U$. Note that with this construction, randomized strategies for the second player (in the stochastic game) correspond to convex combinations of extreme points of the set $\U$. 
When visiting a state in $\Omega_1$, there is an instantaneous payoff of $0$. When visiting a state $(s,a)$ in $\Omega_{2}$, there is an instantaneous payoff of $\bm{p}_{sa}\tr\bm{r}_{sa}$ when the second player chooses $\bm{P} \in \U$. Given a state $\omega \in \Omega_{1}$, the transition function attributes a probability of $1$ to the next state $(s,a) \in \Omega_{2}$ if the first-player chooses action $a \in \A$ and $\omega =s$ for some $s \in \X$. Given a state $\omega \in \Omega_{2}$, the transition function attributes a probability of $p_{sas'}$ to the next state $s' \in \Omega_{1}$ if $\omega = (s,a)$ and the second player chooses action $\bm{P}$. 
We represent this construction in Figure \ref{fig:sa-rec-as-sgs}.
    \begin{figure}[htb]
\begin{center}
\includegraphics[width=0.4\linewidth]{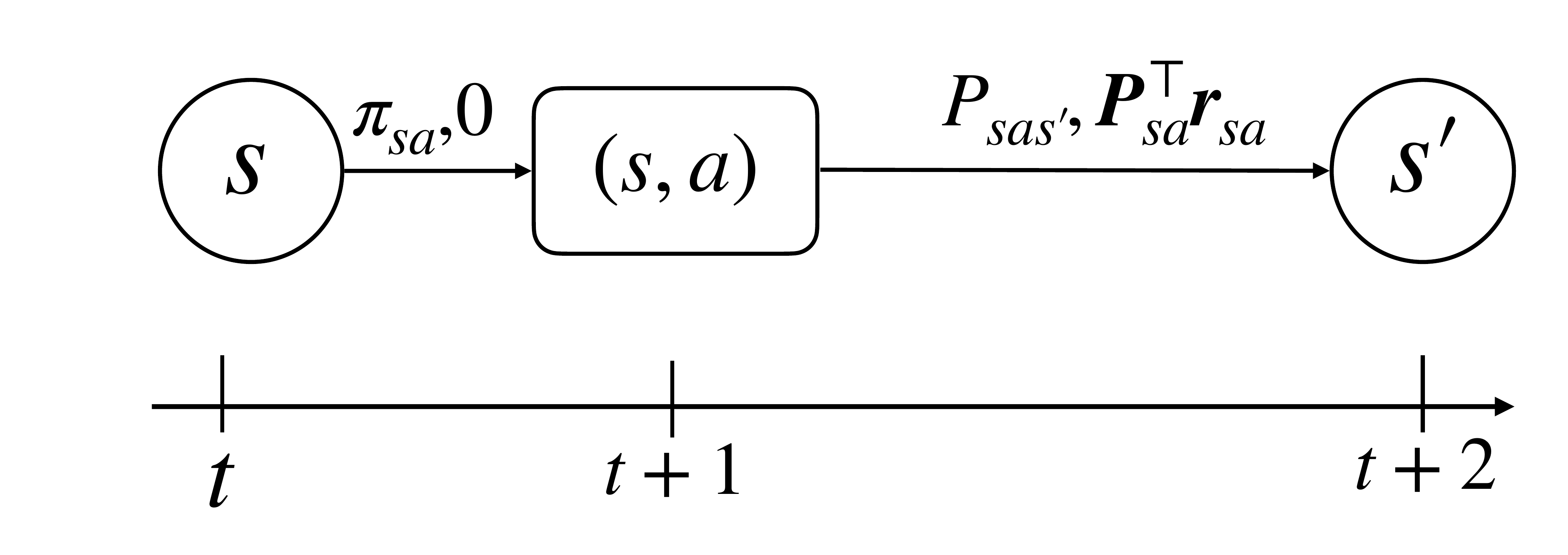}
    \end{center}
\caption{An sa-rectangular robust MDP as a perfect information game. The decision-maker controls the states $s \in \X$ and chooses an action $a \in \A$. The adversary controls the states $(s,a) \in \X \times \A$ and chooses transition probabilities $\bm{p}_{sa} \in \U_{sa}$. The pairs above the arcs represent the transition probabilities and instantaneous payoffs. }\label{fig:sa-rec-as-sgs}
\end{figure}
\begin{remark}[The case of polyhedral uncertainty]
    When the sa-rectangular set $\U$ is {\em polyhedral}, its set of extreme points is finite. From the construction highlighted above, actions of the second player (in the stochastic game) correspond to extreme points of the uncertainty set (in the robust MDP). Therefore, sa-rectangular robust MDPs with polyhedral uncertainty can be reformulated as perfect information stochastic games with {\em finitely} many actions for the second player.
\end{remark}
}
\tb{
\subsection{SGs and s-rectangularity}
\label{app:more-detailed-sgs-rmdps-srec}
We now describe in detail the reformulation of s-rectangular RMDP instances as stochastic games.
Consider an s-rectangular RMDP instance $\left(\X,\A,\bm{r},\U,\bm{p}_{0}\right)$ with a compact convex set $\U$. We consider a stochastic game instance constructed as follows. The set of states is $\Omega = \X$. The set of actions for the first player is $\A$, and the set of actions for the second player is the set of extreme points of $\U$. Similarly as for the construction in the previous section, randomized strategies for the second player correspond to convex combinations of extreme points of the set $\U$. Given a state $\omega \in \Omega$, with $\omega = s \in \X$, and given action $a \in \A$ for the first player and action $\bm{P} \in \U$ for the second player, the instantaneous payoff is $\bm{p}_{sa}\tr\bm{r}_{sa}$ and the transition probabilities to the next state $s' \in \X$ is $p_{sas'}$.
}

 \section{Proof of Proposition \ref{prop:discounted-return-min-inf}}\label{app:history dependent equivalent stationary - discount}
\begin{repeattheorem}[Proposition 2.2]
Let $\U$ be a convex compact s-rectangular uncertainty set. Then
    \begin{equation*}
     \sup_{\pi \in \Pi_{\sf H}} \inf_{\bm{P} \in \U_{\sf H}} R_{\gamma}(\pi,\bm{P}) = \sup_{\pi \in \Pi_{\sf H}} \min_{\bm{P} \in \U} R_{\gamma}(\pi,\bm{P}).
 \end{equation*}
\end{repeattheorem}
\proof{Proof of Proposition \ref{prop:discounted-return-min-inf}.}
Since $\U \subset \U_{\sf H}$, we always have 
$\sup_{\pi \in \Pi_{\sf H}} \inf_{\bm{P} \in \U_{\sf H}} R_{\gamma}(\pi,\bm{P}) \leq  \sup_{\pi \in \Pi_{\sf H}} \min_{\bm{P} \in \U} R_{\gamma}(\pi,\bm{P}).$
We now prove the converse inequality. We have
\begin{align}
    \sup_{\pi \in \Pi_{\sf H}} \min_{\bm{P} \in \U} R_{\gamma}(\pi,\bm{P}) & \leq  \min_{\bm{P} \in \U} \sup_{\pi \in \Pi_{\sf H}} R_{\gamma}(\pi,\bm{P}) \label{eq:proof-eq-loc-0} \\
    & = \min_{\bm{P} \in \U} \max_{\pi \in \Pi_{\sf S}} R_{\gamma}(\pi,\bm{P}) \label{eq:proof-eq-loc-1} \\
    & =  \max_{\pi \in \Pi_{\sf S}} \min_{\bm{P} \in \U} R_{\gamma}(\pi,\bm{P}) \label{eq:proof-eq-loc-2} \\
    & =  \max_{\pi \in \Pi_{\sf S}} \inf_{\bm{P} \in \U_{\sf H}} R_{\gamma}(\pi,\bm{P}) \label{eq:proof-eq-loc-3} \\
    & \leq \sup_{\pi \in \Pi_{\sf H}} \inf_{\bm{P} \in \U_{\sf H}} R_{\gamma}(\pi,\bm{P}) 
\end{align}
where~\eqref{eq:proof-eq-loc-0} follows from weak duality,~\eqref{eq:proof-eq-loc-1} follows from the inner maximization problem being an MDP, for which stationary optimal policies always exist~\citep{puterman2014markov},~\eqref{eq:proof-eq-loc-2} follows from strong duality for s-rectangular RMDPs,~\eqref{eq:proof-eq-loc-3} follows from the inner minimization being an MDP (the adversarial MDP), see~\citep{goyal2022robust,ho2021partial}, and the last inequality following from $\Pi_{\sf S} \subset \Pi_{\sf H}.$ 
\hfill \halmos \endproof

\section{Detail on Example \ref{ex:no-average-limit}}\label{app:detail-example-no-average-limit}
Consider an MDP instance with a single state and two actions $a_{0}$ and $a_{1}$.
Assume that the instantaneous reward for choosing $a_{0}$ is $0$ and the reward for choosing action $a_{1}$ is $1$.
Consider the following history-dependent (Markovian) policy $\pi$ such that at period $t \in \N$, $\pi$ chooses action $a_{1}$ if $2^{k}-1 \leq t \leq 2 \cdot (2^{k}-1)$ for some even integer $k \in \N$, or action $a_{0}$ otherwise. Intuitively, $\pi$ alternatively chooses $a_{0}$ and $a_{1}$, and the length of the interval where $\pi$ selects the same action doubles over time. For instance, the first $15=1+2+4+8$ instantaneous rewards obtained by $\pi$ are $1,0,0,1,1,1,1,0,0,0,0,0,0,0,0.$ We note that the average payoff up to period $T \in \N$ is the frequency of $1$'s during the first $T$ periods. We will call an {\em episode} for any interval of periods over which $\pi$ always chooses the same action. By definition the length of the $k$-th episode if $2^{k-1}$ (periods start at $t=0$).

After $2k$ episodes for $k \geq 1$, the number of periods is $\sum_{\ell=0}^{2k-1} 2^{\ell} = 2^{2k}-1 = 4^k -1$ while the number of $1$'s is $\sum_{\ell=0}^{k-1} 2^{2 \ell} = \sum_{\ell=0}^{k-1} 4^{\ell} = \frac{4^{k}-1}{3}$, so that the frequency of $1$'s is always $1/3$ after $2k$ episodes with $k \geq 1$.

However, after $2k+1$ episodes for $k \geq 1$, the number of periods is $\sum_{\ell=0}^{2k} 2^{\ell} = 2^{2k+1}-1= 2 \cdot 4^{k}-1$, while the number of $1$'s is $\sum_{\ell=0}^{k} 2^{2k} = \frac{4^{k+1}-1}{3} = \frac{4 \cdot 4^{k}-1}{3}$, so that the frequency of $1$'s after $2k+1$ episode with $k \geq 1$ is $\frac{1}{3}\frac{4 \cdot 4^{k}-1}{2 \cdot 4^{k}-1}$. This converges to $2/3$ as $k \rightarrow + \infty.$
    Therefore, we conclude that the average payoff does not have a limit as $T \rightarrow + \infty$.
\section{Some cases where the average worst-case transition probabilities exist}\label{app:ravg-inf-equal-min}
We have the following proposition.
\begin{proposition}\label{prop:infimum-attained}
    Let $\U$ be an s-rectangular compact robust MDP instance. Let $\pi \in \Pi_{\sf S}$. Then the infimum in $\inf_{\bm{P} \in \U} \Ravg(\pi,\bm{P})$ is attained under {\em any} of the following assumptions:
    \begin{enumerate}
        \item For any $(\pi,\bm{P}) \in \Pi_{\sf S} \times \U$, the Markov chain induced by $(\pi,\bm{P})$ over $\X$ is unichain.
        \item For any $(\pi,\bm{P}) \in \Pi_{\sf S} \times \U$, the Markov chain induced by $(\pi,\bm{P})$ over $\X$ is irreducible.
        \item The uncertainty set $\U$ is polyhedral.
    \end{enumerate}
\end{proposition}
\proof{Proof of Proposition \ref{prop:infimum-attained}.}
Under the first two assumptions, the average return $(\pi,\bm{P}) \mapsto \Ravg(\pi,\bm{P})$ is continuous on $\Pi_{\sf S} \times \U$ (exercise 3, section 5.6,~\cite{filar2012competitive} and lemma 5.1.4, \cite{filar2012competitive}). From Weierstrass theorem, $\bm{P} \mapsto \Ravg(\pi,\bm{P})$ attains its minimum over the compact set $\U$.

Let us now assume that $\U$ is polyhedral. Recall that the optimization problem $\inf_{\bm{P} \in \U} \Ravg(\pi,\bm{P})$ can be seen as the adversarial MDP~\citep{ho2021partial}. In this MDP, the set of deterministic policies corresponds to the set of extreme points of $\U$, which is finite. Therefore, $\inf_{\bm{P} \in \U} \Ravg(\pi,\bm{P})$ can be interpreted as an MDP with finitely many actions and finitely many states.  In this case, it is a classical result (see for instance chapter 8 and chapter 9 in \cite{puterman2014markov}) that a stationary deterministic policy is optimal for the average return, so that the infimum is attained in $\inf_{\bm{P} \in \U} \Ravg(\pi,\bm{P})$.
\hfill \halmos \endproof

\section{Discussion on Lemma \ref{lem:avg-rew-nom-eps-stationary}}\label{app:discussion-lemma-avg-rew-adversarial mdp}
In this appendix, we discuss the results in \cite{bierth1987expected}. In particular, the authors in \cite{bierth1987expected} focus on the case of nominal MDP with a finite set of states and compact action sets, where the rewards and the transitions are continuous in the chosen action. In this case, theorem 2.5 in \cite{bierth1987expected} shows that the worst-case average return over the set of history-dependent strategies is the same as the worst-case average return over the set of stationary deterministic strategies.

Recall that the {\em adversarial MDP} described in Section \ref{sec:RMDPs-discounted} has a finite set of states, a compact action set and continuous rewards and transitions. Therefore, we can use the results from theorem 2.5 in \cite{bierth1987expected} to obtain Lemma \ref{lem:avg-rew-nom-eps-stationary}. The only difference is the conclusion of theorem 2.5 in \cite{bierth1987expected} pertains to stationary {\em deterministic} policies, which correspond to the set of extreme points of $\U$. Instead, we wrote Lemma \ref{lem:avg-rew-nom-eps-stationary} with an infimum over $\bm{P} \in \U$, which corresponds to stationary {\em randomized} policy of the adversary. We note that since ${\sf ext}(\U) \subset \U \subset \U_{\sf H}$, Lemma \ref{lem:avg-rew-nom-eps-stationary} holds. 

\begin{remark}
    The authors in \cite{wang2023robust} consider a statement equivalent to Lemma \ref{lem:avg-rew-nom-eps-stationary} for the case of {\em Markovian} adversary in {\em unichain} sa-rectangular RMDPs, see theorem 9 in \cite{wang2023robust}. In particular, \cite{wang2023robust} consider a stationary policy $\pi \in \Pi_{\sf S}$ and study the worst-case of the average return $\lim_{T \rightarrow + \infty} \E_{\pi,\bm{P}} \left[ \frac{1}{T+1}  \sum_{t=0}^{T} r_{s_{t}a_{t}s_{t+1}} \; | \; s_{0} \sim \bm{p}_{0} \right]$ over the set of Markovian adversary $\U_{\sf M}$ defined as $\U_{\sf M} = \U^{\N}$. However, we have shown in Example \ref{ex:no-average-limit} that this limit may not exist when $\bm{P} \in \U_{\sf M}$, which makes the analysis in \cite{wang2023robust} incomplete. 
    We conclude by noting that Lemma \ref{lem:avg-rew-nom-eps-stationary} holds without the unichain assumption and that s-rectangularity is sufficient (instead of sa-rectangularity). 
\end{remark}

\section{Proof of Theorem \ref{th:sa-rec-strong-duality}}\label{app:proof sa-rec strong duality}

\begin{repeattheorem}[Theorem 3.5]
        Consider an sa-rectangular robust MDP with a compact uncertainty set $\U$. Then the following strong duality results hold:
    \begin{align*}
        \sup_{\pi \in \Pi_{\sf H}} \inf_{\bm{P} \in \U} \Ravg(\pi,\bm{P}) & =  \inf_{\bm{P} \in \U} \sup_{\pi \in \Pi_{\sf H}} \Ravg(\pi,\bm{P}) \\
        \max_{\pi \in \Pi_{\sf SD}} \inf_{\bm{P} \in \U} \Ravg(\pi,\bm{P}) & = \inf_{\bm{P} \in \U} \max_{\pi \in \Pi_{\sf SD}} \Ravg(\pi,\bm{P}) 
    \end{align*}
\end{repeattheorem}

\proof{Proof of Theorem \ref{th:sa-rec-strong-duality}.}
We note that \eqref{eq:strong-duality-max-min} is a consequence of \eqref{eq:strong-duality-sup-inf}. Indeed, if \eqref{eq:strong-duality-sup-inf} holds, then
\[ \inf_{\bm{P} \in \U} \max_{\pi \in \Pi_{\sf SD}} \Ravg(\pi,\bm{P}) \leq \inf_{\bm{P} \in \U} \sup_{\pi \in \Pi_{\sf H}} \Ravg(\pi,\bm{P}) =  \sup_{\pi \in \Pi_{\sf H}} \inf_{\bm{P} \in \U} \Ravg(\pi,\bm{P}) = \max_{\pi \in \Pi_{\sf SD}} \inf_{\bm{P} \in \U} \Ravg(\pi,\bm{P})\]
where the first inequality follows from $\Pi_{\sf SD} \subset \Pi_{\sf H}$ and the last equality follows from Theorem \ref{th:avg-reward-main}. We now dedicate our efforts to proving that \eqref{eq:strong-duality-sup-inf} holds.

First, when $\U$ is polyhedral, we have
\[ \inf_{\bm{P} \in \U} \sup_{\pi \in \Pi_{\sf H}} \Ravg(\pi,\bm{P}) = \inf_{\bm{P} \in \U} \max_{\pi \in \Pi_{\sf SD}} \Ravg(\pi,\bm{P}) = \max_{\pi \in \Pi_{\sf SD}} \inf_{\bm{P} \in \U}  \Ravg(\pi,\bm{P}) \leq \sup_{\pi \in \Pi_{\sf H}} \inf_{\bm{P} \in \U} \Ravg(\pi,\bm{P}) \]
where the first equality comes from the inner supremum problem being an MDP, the second equality comes from Theorem 1 in \cite{liggett1969stochastic}, and the inequality comes from $\Pi_{\sf SD} \subset \Pi_{\sf H}.$ This shows that \eqref{eq:strong-duality-sup-inf} holds when $\U$ is polyhedral.

We now turn to showing \eqref{eq:strong-duality-sup-inf} when $\U$ is a compact sa-rectangular uncertainty set. Let $\epsilon>0$ and consider the set $\U_{\sf f}$ defined similarly as in the second step of the proof of Theorem \ref{th:avg-reward-main}. We have
\begin{align}
     \inf_{\bm{P} \in \U} \sup_{\pi \in \Pi_{\sf H}} \Ravg(\pi,\bm{P}) & \leq  \inf_{\bm{P} \in \U_{\sf f}} \sup_{\pi \in \Pi_{\sf H}} \Ravg(\pi,\bm{P}) \label{eq:proof-sd-sup-inf-0} \\
     & =\sup_{\pi \in \Pi_{\sf H}} \inf_{\bm{P} \in \U_{\sf f}} \Ravg(\pi,\bm{P}) \label{eq:proof-sd-sup-inf-1} \\
     & = \max_{\pi \in \Pi_{\sf SD}} \inf_{\bm{P} \in \U_{\sf f}} \Ravg(\pi,\bm{P}) \label{eq:proof-sd-sup-inf-2} \\
     & \leq \max_{\pi \in \Pi_{\sf SD}} \inf_{\bm{P} \in \U} \Ravg(\pi,\bm{P}) \label{eq:proof-sd-sup-inf-3} + \epsilon \\
     & = \sup_{\pi \in \Pi_{\sf H}} \inf_{\bm{P} \in \U} \Ravg(\pi,\bm{P}) + \epsilon \label{eq:proof-sd-sup-inf-4}
\end{align}
where \eqref{eq:proof-sd-sup-inf-0} follows from $\U_{\sf f} \subset \U$, \eqref{eq:proof-sd-sup-inf-1} follows from the first two steps of this proof and $\U_{\sf f}$ having only finitely many extreme points, \eqref{eq:proof-sd-sup-inf-2} follows from Theorem \ref{th:avg-reward-main}, \eqref{eq:proof-sd-sup-inf-3} follows from the construction of $\U_{\sf f}$, and \eqref{eq:proof-sd-sup-inf-4} follows from Theorem \ref{th:avg-reward-main} again. Since the inequalities above are true for any choice of $\epsilon>0$, we conclude that 
$\sup_{\pi \in \Pi_{\sf H}} \inf_{\bm{P} \in \U} \Ravg(\pi,\bm{P}) = \inf_{\bm{P} \in \U} \sup_{\pi \in \Pi_{\sf H}} \Ravg(\pi,\bm{P}).$
\hfill \halmos \endproof
\section{Proof of Corollary \ref{cor:avg-reward-main-other-ravg}}\label{app:proof-corollaries}

\begin{repeattheorem}[Corollary 3.7]
       Consider an sa-rectangular robust MDP with a compact uncertainty set $\U$. Let $\hRavg$ as defined in \eqref{eq:liminf-exp}, \eqref{eq:limsup-exp}, or \eqref{eq:exp-liminf}.
    
    \begin{enumerate}
        \item[{\bf P0}.] There exists an average optimal policy that is stationary and deterministic, and which coincides with an average optimal policy for $\Ravg$ as in \eqref{eq:definition-average-return-exp-limsup}:
    \[ 
    \sup_{\pi \in \Pi_{\sf H}} \inf_{\bm{P} \in \U} \hRavg(\pi,\bm{P}) \;=\; 
    \max_{\pi \in \Pi_{\sf SD}} \inf_{\bm{P} \in \U} \hRavg(\pi,\bm{P}) = \max_{\pi \in \Pi_{\sf SD}} \inf_{\bm{P} \in \U} \Ravg(\pi,\bm{P}).
    \] 
    \item[{\bf P1}.] The strong duality results from Theorem \ref{th:sa-rec-strong-duality} still hold when  replacing $\Ravg$ by $\hRavg$.
    \item[{\bf P2}.] The equivalence between stationary adversaries and history-dependent adversaries from Theorem \ref{th:avg-rew-history-dep-adversaries} still hold when replacing $\Ravg$ by $\hRavg$.
    \end{enumerate}
\end{repeattheorem}
\proof{Proof of Corollary \ref{cor:avg-reward-main-other-ravg}.} 
We prove the three statements at once. The proof proceeds in three parts. Let $\hRavg$ as defined in \eqref{eq:liminf-exp}, \eqref{eq:limsup-exp}, or \eqref{eq:exp-liminf}. 

{\bf Part 1.}
Crucially, for any $(\pi,\bm{P}) \in \Pi_{\sf H} \times \U$, by the Fatou-Lebesgue theorem, we have
\begin{equation}\label{eq:ineq-lim-supinf-exp}
   \hRavg(\pi,\bm{P}) \leq \Ravg(\pi,\bm{P})
\end{equation}
and $\hRavg(\pi,\bm{P}) = \Ravg(\pi,\bm{P})$ when $(\pi,\bm{P}) \in \Pi_{\sf S} \times \U$ because the limits exist.
From this, we have that
\[\max_{\pi \in \Pi_{\sf SD}} \inf_{\bm{P} \in \U} \hRavg(\pi,\bm{P}) = \max_{\pi \in \Pi_{\sf SD}} \inf_{\bm{P} \in \U} \Ravg(\pi,\bm{P})= \sup_{\pi \in \Pi_{\sf H}} \inf_{\bm{P} \in \U} \Ravg(\pi,\bm{P})\geq \sup_{\pi \in \Pi_{\sf H}} \inf_{\bm{P} \in \U} \hRavg(\pi,\bm{P})\]
where the first equality holds from the optimization variables $(\pi,\bm{P})$ being in $\Pi_{\sf SD} \times \U$, the second equality holds from the first step of the proof, and the inequality follows from~\eqref{eq:ineq-lim-supinf-exp}. This shows that when $\hRavg$ is defined as in~\eqref{eq:liminf-exp}, \eqref{eq:limsup-exp}, or \eqref{eq:exp-liminf}, we have
$\sup_{\pi \in \Pi_{\sf H}} \inf_{\bm{P} \in \U} \hRavg(\pi,\bm{P}) = \max_{\pi \in \Pi_{\sf SD}} \inf_{\bm{P} \in \U} \hRavg(\pi,\bm{P}).$ Now since $\hRavg$ and $\Ravg$ coincides over $\Pi_{\sf S} \times \U$, we can conclude that $\max_{\pi \in \Pi_{\sf SD}} \inf_{\bm{P} \in \U} \hRavg(\pi,\bm{P}) = \max_{\pi \in \Pi_{\sf SD}} \inf_{\bm{P} \in \U} \Ravg(\pi,\bm{P})$.

\vspace{2mm}
{\bf Part 2.}
We now turn to proving that the following strong duality results hold:
    \begin{align}
        \sup_{\pi \in \Pi_{\sf H}} \inf_{\bm{P} \in \U} \hRavg(\pi,\bm{P}) & =  \inf_{\bm{P} \in \U} \sup_{\pi \in \Pi_{\sf H}} \hRavg(\pi,\bm{P}) \label{eq:strong-duality-sup-inf-other-ravg}\\
        \max_{\pi \in \Pi_{\sf SD}} \inf_{\bm{P} \in \U} \hRavg(\pi,\bm{P}) & = \inf_{\bm{P} \in \U} \max_{\pi \in \Pi_{\sf SD}} \hRavg(\pi,\bm{P}) \label{eq:strong-duality-max-min-other-ravg}
    \end{align}
    Similarly as for the proof of Theorem \ref{th:sa-rec-strong-duality}, it is sufficient to show \eqref{eq:strong-duality-sup-inf-other-ravg}. We have
\begin{align}
     \inf_{\bm{P} \in \U} \sup_{\pi \in \Pi_{\sf H}} \hRavg(\pi,\bm{P}) & \leq  \inf_{\bm{P} \in \U} \sup_{\pi \in \Pi_{\sf H}} \Ravg(\pi,\bm{P}) \label{eq:proof-sd-sup-inf-0-0} \\
     & =\sup_{\pi \in \Pi_{\sf H}} \inf_{\bm{P} \in \U} \Ravg(\pi,\bm{P}) \label{eq:proof-sd-sup-inf-1-0} \\
     & = \max_{\pi \in \Pi_{\sf SD}} \inf_{\bm{P} \in \U} \Ravg(\pi,\bm{P}) \label{eq:proof-sd-sup-inf-2-0} \\
     & \leq \max_{\pi \in \Pi_{\sf SD}} \inf_{\bm{P} \in \U} \hRavg(\pi,\bm{P}) \label{eq:proof-sd-sup-inf-3-0}\\
     & \leq \sup_{\pi \in \Pi_{\sf H}} \inf_{\bm{P} \in \U} \hRavg(\pi,\bm{P}) \label{eq:proof-sd-sup-inf-4-0}
\end{align}
where \eqref{eq:proof-sd-sup-inf-0-0} follows from $\hRavg(\pi,\bm{P}) \leq \Ravg(\pi,\bm{P})$ for $(\pi,\bm{P}) \in \Pi_{\sf H} \times \U$, \eqref{eq:proof-sd-sup-inf-1-0} follows from Theorem \ref{th:sa-rec-strong-duality}, \eqref{eq:proof-sd-sup-inf-2-0} follows from Theorem \ref{th:avg-reward-main}, \eqref{eq:proof-sd-sup-inf-3-0} follows from $\hRavg(\pi,\bm{P}) = \Ravg(\pi,\bm{P})$ for  $(\pi,\bm{P}) \in \Pi_{\sf S} \times \U$, and \eqref{eq:proof-sd-sup-inf-4-0} follows from $\U_{\sf SD} \subset \U_{\sf H}$.
Therefore,
$\sup_{\pi \in \Pi_{\sf H}} \inf_{\bm{P} \in \U} \hRavg(\pi,\bm{P}) = \inf_{\bm{P} \in \U} \sup_{\pi \in \Pi_{\sf H}} \hRavg(\pi,\bm{P}).$

\vspace{2mm}
{\bf Part 3.}
Finally, we prove that Theorem \ref{th:avg-rew-history-dep-adversaries} also holds for $\hRavg$, i.e. we show that
\[\sup_{\pi \in \Pi_{\sf H}} \inf_{\bm{P} \in \U} \hRavg(\pi,\bm{P}) = \sup_{\pi \in \Pi_{\sf H}} \inf_{\bm{P} \in \U_{H}} \hRavg(\pi,\bm{P}).\]
Note that proof for Theorem \ref{th:avg-rew-history-dep-adversaries} only relies on Theorem \ref{th:avg-reward-main}, Theorem \ref{th:sa-rec-strong-duality} and Lemma \ref{lem:avg-rew-nom-eps-stationary}. We have proved above that Theorem \ref{th:avg-reward-main} and Theorem \ref{th:sa-rec-strong-duality} still hold when we replace $\hRavg$ by $\Ravg$. The authors in \cite{bierth1987expected} also prove Lemma \ref{lem:avg-rew-nom-eps-stationary} for $\hRavg$ as in \eqref{eq:liminf-exp}, \eqref{eq:limsup-exp}, or \eqref{eq:exp-liminf}. Therefore, the proof of Theorem \ref{th:avg-rew-history-dep-adversaries} follows verbatim when replacing $\Ravg$ by $\hRavg$.
\hfill \halmos \endproof

\section{Proof of Theorem \ref{th:sa-rec-no-blackwell}}\label{app:proof-th-sa-rec-no-Blackwell}
This section is dedicated to proving Theorem \ref{th:sa-rec-no-blackwell}:
\begin{repeattheorem}[Theorem 4.4]
     There exists an sa-rectangular robust MDP instance, with a compact convex uncertainty set $\U$, and with no Blackwell optimal policy:
    \[\forall \; \pi \in \Pi_{\sf H}, \forall \; \gamma \in (0,1), \exists \; \gamma' \in (\gamma,1), \min_{\bm{P} \in \U} R_{\gamma'}(\pi,\bm{P}) < \sup_{\pi' \in \Pi_{\sf S}} \min_{\bm{P} \in \U} R_{\gamma'}(\pi',\bm{P}). \]
\end{repeattheorem}

We consider the following sa-rectangular robust MDP instance, which builds upon the counterexample for Proposition \ref{prop:inf-is-not-min-avg-rew}. There are three states $\{s_{0},A_{-1},A_{0}\}$, where $A_{-1}$ is an absorbing state with reward $-1$ and $A_{0}$ is an absorbing state with reward $0$. The decision-maker starts in state $s_{0}$ with probability $1$. The instantaneous reward in $s_{0}$ is $0$. There are two actions $\{a_{1},a_{2}\}$, and for $i \in \{1,2\}$,
\begin{align*}
    \U_{s_{0}a_{i}} = \{ \left(1-\alpha,\beta,\alpha-\beta\right) \; | \; \alpha \in [0,1], \beta \in [0,\alpha],\beta \leq D_{i}(\alpha) \}
\end{align*}
for $D_{1}$ and $D_{2}$ two concave functions such that $D_{1}(\alpha) \leq \alpha, D_{2}(\alpha) \leq \alpha.$ We will specify $D_{1}$ and $D_{2}$ later.
For a vector $\left(1-\alpha,\beta,\alpha-\beta\right)$ in the simplex, we can interpret $\alpha$ as the probability to leave state $s_{0}$ and $\beta$ as the probability to go to state $A_{-1}$. This is represented in Figure \ref{fig:simplex-alpha-beta}.
Since state $s_{0}$ is the only non-absorbing state, we can identify policies and the action ($a_{1}$ or $a_{2}$) chosen at $s_{0}$.
Since we are considering an sa-rectangular robust MDP, for every discount factor $\gamma \in (0,1)$ there is a deterministic stationary optimal policy in $\{a_{1},a_{2}\}$. The worst-case return of $a_{1}$ and $a_{2}$ are defined as, for $i \in \{1,2\}$
\begin{align*}
    R_{\gamma}(a_{i}) = \min_{\bm{P} \in \U} R_{\gamma}(a_{i},\bm{P}) = \frac{1}{1-\gamma} \min_{\alpha \in [0,1], \beta \in [0,\alpha],\left(1-\alpha,\beta,\alpha-\beta\right)  \in \U_{sa_{i}}} - \gamma \frac{\beta}{1-\gamma + \gamma\alpha}.
\end{align*}
From the definition of $\U_{sa_{i}}$ for $i \in \{1,2\}$, we can write $R_{\gamma}(a_{i})$ as
\begin{align*}
     R_{\gamma}(a_{i}) 
     = - \frac{1}{1-\gamma} \max_{\alpha \in [0,1]} \;  \gamma \frac{D_{i}(\alpha)}{1-\gamma + \gamma\alpha}.
\end{align*}
We will construct two sequences $\left(\gamma_{k}\right)_{k \in \N}$ and $\left(\gamma'_{k}\right)_{k \in \N}$ such that for any $k \in \N$, $a_{1}$ is the unique optimal policy for $\gamma = \gamma_{k}$ and $a_{2}$ is the unique optimal policy for $\gamma = \gamma'_{k}$, and $\gamma_{k} \rightarrow 1, \gamma'_{k} \rightarrow 1$ as $k \rightarrow + \infty$.
\paragraph{First step.} We first consider the following function $D:[0,1] \rightarrow [0,1]$ defined as $D(\alpha) = \alpha(1-\alpha)$, and consider the following optimization program:
\[\max_{\alpha \in [0,1]} \;  \gamma \frac{D(\alpha)}{1-\gamma + \gamma\alpha} = \max_{\alpha \in [0,1]} \;  \gamma \frac{\alpha (1-\alpha)}{1-\gamma + \gamma\alpha}.\]
The derivative of $\alpha \mapsto \frac{\alpha (1-\alpha)}{1-\gamma + \gamma\alpha}$ is 
$\frac{(1-2\alpha)(1-\gamma + \gamma\alpha) - \alpha(1-\alpha)\gamma}{(1-\gamma + \gamma\alpha)^2} = \frac{1-\gamma - 2 \alpha (1-\gamma) - \gamma\alpha^2}{(1-\gamma + \gamma\alpha)^2},$
which shows that the maximum of $\alpha \mapsto \frac{\alpha (1-\alpha)}{1-\gamma + \gamma\alpha}$ is attained at $\alpha\opt(\gamma)= \frac{\sqrt{1-\gamma}-(1-\gamma)}{\gamma}$. Note that $\gamma \mapsto \alpha\opt(\gamma)$ is a decreasing function on $(0,1]$ that we can continuously extend to $[0,1]$, with $\alpha\opt(1) = 0$ and $\alpha\opt(0)=1/2$.
\paragraph{Second step.} We now construct two piece-wise affine functions $D_{1}, D_{2}$, that are (pointwise) smaller than $D:\alpha \mapsto \alpha (1-\alpha)$ on $[0,1]$ , and that intersect infinitely many times on any interval $(0,\alpha_{0})$ for any $\alpha_{0} \in (0,1/2).$ In particular, we have the following lemma. We call {\em breakpoints} of a piecewise affine function the points where it changes slopes.
\begin{lemma}\label{lem:pathological-functions}
    Let $D_{1}:[0,1] \rightarrow \R,D_{2}:[0,1] \rightarrow \R$ defined as follows.
    \begin{enumerate}
        \item $D_{1}$ is piecewisewise affine, its set of breakpoints is $\{\frac{1}{2k} \; | \; k \geq 1\}$, where $D_{1}(\frac{1}{2k}) = D(\frac{1}{2k}), k \geq 1$.
        \item $D_{2}$ is piecewisewise affine, its set of breakpoints is $\{\frac{1}{2k+1} \; | \; k \geq 1\}$, where $D_{2}(\frac{1}{2k+1}) = D(\frac{1}{2k+1}), k \geq 1$.
    \end{enumerate}
    Then $D_{1},D_{2}$ are continuous, concave, and increasing functions on $[0,1/2]$. Additionally, $D_{1}$ and $D_{2}$ intersect infinitely many times on any interval $[0,\alpha_{0}]$ with $\alpha_{0} \in (0,1/2)$.
\end{lemma}

\proof{Proof of Lemma~\ref{lem:pathological-functions}.}
    We first show that $D_{1}$ is a concave, strictly increasing function. For $k \geq 1$ and $x \in \left(\frac{1}{2k+2},\frac{1}{2k}\right)$, there exists $a_{1,k},b_{1,k} \in \R$ such that $D_{1}(x) = a_{1,k}x+b_{1,k}$. In particular, we have
    \[a_{1,k} = \frac{D_{1}\left( \frac{1}{2k} \right) - D_{1}\left( \frac{1}{2k+1} \right)}{\frac{1}{2k}-\frac{1}{2k+2}}= \frac{\frac{1}{2k}\left(1-\frac{1}{2k}\right)-\frac{1}{2k+2}\left(1-\frac{1}{2k+2}\right)}{\frac{1}{2k}-\frac{1}{2k+2}} = \frac{2k^2-1}{2k^2+2k}.\]
    We also have $ a_{1,k} \frac{1}{2k} + b_{1,k}  = D_{1}\left(\frac{1}{2k}\right)$, so that $b_{1,k} = \frac{1}{4k(k+1)}.$
    Note that $k \mapsto \frac{2k^2-1}{2k^2+2k}$ is an increasing function, with $\lim_{k \rightarrow + \infty} \frac{2k^2-1}{2k^2+2k} = 1.$ Therefore, we have shown that $k \mapsto a_{1,k}$ is an increasing function. Since $D_{1}$ is continuous by construction, we conclude that $D_{1}$ is a concave, increasing function. 
     The same approach shows that $D_{2}$ is concave and strictly increasing.

    We now show that the two functions intersect infinitely many times on any interval $[0,\lambda_{0})$ with $\lambda_{0} > 0$. Let $k \geq 1$. Then there exists a zero of $x \mapsto D_{1}(x)-D_{2}(x)$ on any interval $\left(\frac{1}{k+1},\frac{1}{k}\right)$.
    From the concavity of $x \mapsto \sqrt{x}$, we have that $D_{2}\left(\frac{1}{2k+1}\right)-D_{1}\left(\frac{1}{2k+1}\right) \geq 0$ and $D_{2}\left(\frac{1}{2k+2}\right)-D_{1}\left(\frac{1}{2k+2}\right) \leq 0$. Since $D_{1},D_{2}$ are continuous, there exists a zero of $x \mapsto D_{1}(x)-D_{2}(x)$ on any interval $\left(\frac{1}{2k+2},\frac{1}{2k+1}\right)$ for $k \geq 1$. The same approach shows that there exists a zero $x \mapsto D_{1}(x)-D_{2}(x)$ on any interval $\left(\frac{1}{2k+1},\frac{1}{2k}\right)$ for $k \geq 1$.
\hfill \halmos \endproof
\paragraph{Third step.} We now construct a sequence of discount factors $\left(\gamma_{k}\right)_{k \in \N}$ such that $\gamma_{k} \rightarrow 1$ and $a_{1}$ is the unique optimal policy for $\gamma = \gamma_{k}$ for any $k \in \N$.
Let $k \in \N$ and $\gamma_{k} \in (0,1)$ such that $\alpha\opt(\gamma_{k}) = \frac{1}{2k+1}$. From the strict monotonicity of $\gamma \mapsto \alpha\opt(\gamma)$ we know that $\gamma_{k} \rightarrow 1$ as $k \rightarrow + \infty$.
By construction, 
\begin{align}
    R_{\gamma_{k}}(a_{2}) & = - \frac{1}{1-\gamma_k}\max_{\alpha \in [0,1]} \frac{D_{2}(\alpha)}{1-\gamma_{k}+\gamma_{k}\alpha} \nonumber \\
    & \geq - \frac{1}{1-\gamma_k} \max_{\alpha \in [0,1]} \frac{D(\alpha)}{1-\gamma_{k}+\gamma_{k}\alpha} \label{eq:prf-ctr-ex-0} \\
    & = - \frac{1}{1-\gamma_k} \frac{D(\frac{1}{2k+1})}{1-\gamma_{k} + \gamma_{k}\frac{1}{2k+1}} \label{eq:prf-ctr-ex-1} \\
    & = - \frac{1}{1-\gamma_k} \frac{D_{2}(\frac{1}{2k+1})}{1-\gamma_{k}+\gamma_{k}\frac{1}{2k+1}} \label{eq:prf-ctr-ex-2}
\end{align}
where~\eqref{eq:prf-ctr-ex-0} follows from $D_{2}(\alpha) \leq D(\alpha)$ for any $\alpha \in (0,1/2)$, where~\eqref{eq:prf-ctr-ex-1} follows from the definition of $\gamma_{k}$ such that $\alpha\opt(\gamma_{k}) = \frac{1}{2k+1}$, and~\eqref{eq:prf-ctr-ex-2} follows from $D_{2}(\frac{1}{2k+1}) = D(\frac{1}{2k+1})$ by construction of $D_{2}$ as in Lemma \ref{lem:pathological-functions}. This shows that the maximum in $R_{\gamma_{k}}(a_{2})$ is attained at $\alpha\opt(\gamma_{k})= \frac{1}{2k+1}$. 

We now show that $R_{\gamma_{k}}(a_{1}) > R_{\gamma_{k}}(a_{2})$.  
The function $\alpha \mapsto \frac{D_{1}(\alpha)}{1-\gamma_{k}+\gamma_{k}\alpha}$ is continuous on the compact set $[0,1/2]$, hence it attains its maximum. We distinguish two cases. 

{\em Case 1.} Suppose that this maximum is attained at $\alpha_{1} \in \{1/2k' \; | \; k' \in \N\}$. Then $D_{1}(\alpha_{1}) = D(\alpha_{1})$, and
\[\frac{D_{1}(\alpha_{1})}{1-\gamma_{k}+\gamma_{k}\alpha_{1}} = \frac{D(\alpha_{1})}{1-\gamma_{k}+\gamma_{k}\alpha_{1}} < \max_{\alpha \in [0,1/2]} \frac{D(\alpha_{1})}{1-\gamma_{k}+\gamma_{k}\alpha_{1}}\]
where the strict inequality follows from $\alpha\opt(\gamma_k) = 1/(2k+1) \neq \alpha_1.$

{\em Case 2.}
Otherwise, the maximum of $\alpha \mapsto \frac{D_{1}(\alpha)}{1-\gamma_{k}+\gamma_{k}\alpha}$ is attained at  $\alpha_1 \notin [0,1/2]$ but $\alpha_1 \notin \{1/2k' \; | \; k' \in \N\}$. In this case we have by construction that $D_{1}(\alpha_{1}) < D(\alpha_{1})$ so that 
\[\frac{D_{1}(\alpha_{1})}{1-\gamma_{k}+\gamma_{k}\alpha_{1}} < \frac{D(\alpha_{1})}{1-\gamma_{k}+\gamma_{k}\alpha_{1}} \leq \max_{\alpha \in [0,1/2]} \frac{D(\alpha_{1})}{1-\gamma_{k}+\gamma_{k}\alpha_{1}}.\]

Overall, we conclude that $R_{\gamma_k}(a_{1}) > R_{\gamma_k}(a_{2})$ in both case $1$ and case $2$. Threfore, we have constructed a sequence $\left(\gamma_k\right)_{k \in \N}$ such that $R_{\gamma_k}(a_{1}) > R_{\gamma_k}(a_{2}), \forall \; k \geq 1$. The construction of a sequence $\left(\gamma_{k}'\right)_{k \in \N}$ such that $R_{\gamma_{k}'}(a_{1}) > R_{\gamma_{k}'}(a_{2}), \forall \; k \geq 1$ is analogous, by choosing $\gamma'_{k}$ such that $\alpha\opt(\gamma'_k) = \frac{1}{2k}$.  This shows that there does not exist a stationary deterministic Blackwell optimal policy. Since for any $k \in \N$, at $\gamma=\gamma_k$ the optimal policy is unique {\em and} stationary deterministic, this shows that there does not exist a Blackwell optimal policy that is history-dependent. This concludes the construction of our counterexample.

\begin{remark}
    \cite{wang2023robust} show the existence of a stationary deterministic Blackwell optimal policy for sa-rectangular RMDPs, under two conditions: (a) the Markov chains induced by any pair of policy and transition probabilities are unichain, (b) there is a unique average optimal policy. The second assumption precludes two value functions from intersecting infinitely often as $\gamma \rightarrow 1$, which is exactly the problematic behavior in our counterexample. In fact, our counterexample violates both assumptions: there are two absorbing states (which violates the unichain assumption), and actions $a_{1}$ and $a_{2}$ are average optimal. 
\end{remark}
\section{Proof of Theorem \ref{th:sa-rec-epsilon-blackwell-opt-same}}\label{app:sa-rec-epsilon-blackwell-opt-same}

\begin{repeattheorem}[Theorem 4.5]
    Let $\U$ be an sa-rectangular compact uncertainty set. Then there exists a stationary deterministic policy that is $\epsilon$-Blackwell optimal for all $\epsilon>0$, i.e., $\exists \; \pi \in \Pi_{\sf SD}, \forall \; \epsilon>0,\exists \; \gamma_{\epsilon} \in (0,1)$ such that
\[  \min_{\bm{P} \in \U} (1-\gamma)R_{\gamma}(\pi,\bm{P}) \geq \sup_{\pi' \in \Pi_{\sf S}} \min_{\bm{P} \in \U} (1-\gamma)R_{\gamma}(\pi',\bm{P}) - \epsilon, \forall \; \gamma \in (\gamma_{\epsilon},1).\]
\end{repeattheorem}

\proof{Proof of Theorem \ref{th:sa-rec-epsilon-blackwell-opt-same}.}

We first fix a stationary policy $\pi \in \Pi_{\sf S}$ and a discount factor $\gamma \in (0,1)$. 
Theorem \ref{th:eps-blackwell-compact-mdp} shows that for each $\epsilon>0$ and each policy $\pi \in \Pi_{\sf S}$, there exists an $\epsilon$-Blackwell optimal policy for the adversarial MDP. In particular, for each $\epsilon>0$ and each policy $\pi \in \Pi_{\sf S}$, there exists a transition probability matrix $\bm{P}^{\pi}_{\epsilon} \in \U$ and a discount factor $\gamma^{\pi}_{\epsilon} \in (0,1),$ such that
\[
  (1-\gamma)\left(v^{\pi,\bm{P}^{\pi}_{\epsilon}}_{\gamma,s} - v^{\pi,\U}_{\gamma,s}\right) \leq \epsilon,
  \quad \forall \; \gamma > \gamma^{\pi}_{\epsilon}, \forall \; s \in \X.
\]
Assume a fixed $\epsilon>0$. For any discount factor $\gamma \in (0,1)$ consider the set $\{R_{\gamma}(\pi,\bm{P}^{\pi}_{\epsilon/2}) \; | \; \pi \in \Pi_{\sf SD}\}$. This set is finite since $\pi \in \Pi_{\sf SD}$ is finite. We then define $\hat{\pi}_{\gamma} \in \arg \max_{\pi \in \Pi_{\sf SD}} R_{\gamma}(\pi,\bm{P}^{\pi}_{\epsilon/2})$
Since $\hat{\pi}_{\gamma} \in \Pi_{\sf SD}, \forall \; \gamma \in (0,1)$ and $\Pi_{\sf SD}$ is finite, we can choose a sequence of discount factors $\left(\gamma_{n}\right)_{n \geq 1}$ such that for all $n \in \N$, we have $\hat{\pi}_{\gamma_{n}} = \pi \in \Pi_{\sf SD}$ for some $\pi \in \Pi_{\sf SD}$ and $\gamma_{n} \rightarrow 1$ as $n \rightarrow + \infty$ (see for instance lemma F.1 in \cite{goyal2022robust}). Next, we show that the policy $\pi$ is $\epsilon$-Blackwell optimal.

Indeed, suppose that $\pi$ is not $\epsilon$-Blackwell optimal: we can construct a sequence of discount factors $\left(\gamma'_{n}\right)_{n \in \N} \in (0,1)^\N$ and a sequence of states $\left(s_n\right)_{n \in \N} \in \X^\N$ such that $(1-\gamma'_{n})\left(v^{\star,\U}_{\gamma'_{n},s_n} - v^{\pi,\U}_{\gamma'_{n},s_n}\right) > \epsilon, \forall \; n \geq 1.$ and $\gamma'_n \rightarrow 1$ as $n \rightarrow + \infty$.
    Since for each $\gamma'_{n}$ the optimal discount policy can be chosen stationary and deterministic, i.e., in $\Pi_{\sf SD}$, we can find a policy $\pi'$ and a state $s \in \X$ such that 
$(1-\gamma'_{n})\left(v^{\pi',\U}_{\gamma'_{n},s} - v^{\pi,\U}_{\gamma'_{n},s} \right)> \epsilon, \forall \; n \geq 1.$
    By definition of the transition probabilities $\bm{P}^{\pi}_{\epsilon/2}$ and  $\bm{P}^{\pi'}_{\epsilon/2}$ we have
    \begin{align*}
      (1-\gamma)\left(v^{\pi',\bm{P}^{\pi'}_{\epsilon/2}}_{\gamma,s}- v^{\pi',\U}_{\gamma,s} \right) \leq \epsilon/2,
      \quad\forall \; \gamma > \gamma^{\pi'}_{\epsilon/2},\\
      (1-\gamma)\left(v^{\pi,\bm{P}^{\pi}_{\epsilon/2}}_{\gamma,s}- v^{\pi,\U}_{\gamma,s}\right) \leq \epsilon/2,
      \quad\forall \; \gamma > \gamma^{\pi}_{\epsilon/2}.
    \end{align*}
    Since $\gamma'_{n} \rightarrow 1$ as $n \rightarrow + \infty$, we conclude that for some $n_{0} \in \N$ we have
\begin{equation}\label{eq:prf-1}
(1-\gamma'_{n})\left(v^{\pi',\bm{P}^{\pi'}_{\epsilon/2}}_{\gamma'_{n},s} - v^{\pi,\bm{P}^{\pi}_{\epsilon/2}}_{\gamma'_{n},s}\right) > \epsilon/2 >0, \forall \; n \geq n_{0}.
    \end{equation}
    But by construction of $\pi$, we have $\pi \in \arg \max_{\hat{\pi} \in \Pi_{\sf SD}} v^{\hat{\pi},\bm{P}^{\hat{\pi}}_{\epsilon/2}}_{\gamma_{n}}, \forall \; n \geq 1.$
    From this we conclude that \begin{equation}\label{eq:prf-2}
v^{\pi',\bm{P}^{\pi'}_{\epsilon/2}}_{\gamma_{n},s} - v^{\pi,\bm{P}^{\pi}_{\epsilon/2}}_{\gamma_{n},s} \leq 0, \forall \; n \geq n_{0}.
    \end{equation}
    Combining~\eqref{eq:prf-1} and~\eqref{eq:prf-2}, we obtain that $\gamma \mapsto v^{\pi',\bm{P}^{\pi'}_{\epsilon/2}}_{\gamma,s} - v^{\pi,\bm{P}^{\pi}_{\epsilon/2}}_{\gamma,s}$ must cancel infinitely many times in the interval $(0,1)$, which is a contradiction since this is a rational function (lemma 10.1.3, \cite{puterman2014markov}).
    Therefore,  for any $\epsilon>0$ we can choose $\pi_{\epsilon} \in \Pi_{\sf SD}$ that is $\epsilon$-Blackwell optimal. Since $\Pi_{\sf SD}$ is finite, the same stationary deterministic policy can be chosen $\epsilon$-Blackwell optimal for all $\epsilon>0$.
\hfill \halmos \endproof

\section{Proof of Theorem \ref{th:sa-rec-relation-avg-rew-blackwell}}\label{app:proof-th-relation-blackwell-avg}

\begin{repeattheorem}[Theorem 4.6]
\tb{Consider an sa-rectangular robust MDP with a compact uncertainty set $\U$. Let $\pi \in \PiSD$. Then 
\begin{center}
    $\pi$ is $\epsilon$-Blackwell optimal for all $\epsilon>0 \iff \pi$ is an average optimal policy.
\end{center}} 
\end{repeattheorem}
\tb{Lemma \ref{lem:aux-1}, Lemma \ref{lem:aux-2} and Lemma \ref{lem:aux-3} are proved in Part 1 of our proof of Theorem \ref{th:sa-rec-relation-avg-rew-blackwell}.}
\proof{Proof of Theorem \ref{th:sa-rec-relation-avg-rew-blackwell}.} 
\tb{
\paragraph{Part 1:} Assume that $\pi \in \PiSD$ is $\epsilon$-Blackwell optimal for all $\epsilon>0$. We will prove that $\pi$ is also average optimal.}

The proof of this first part proceeds in four steps.
\paragraph{First step.} The first step proves to Lemma \ref{lem:aux-1}. We start by studying the limit behavior of the (worst-case) discounted return of a policy $\pi \in \Pi_{\sf S}$ as $\gamma \rightarrow 1$. Let $\epsilon>0$. From Theorem \ref{th:eps-blackwell-compact-mdp}, we know that there exists $\gamma_{0} \in (0,1)$ and $\bm{P}_{\epsilon} \in \U$ such that
\begin{equation}\label{eq:proof-app-1}
    (1-\gamma)R_{\gamma}(\pi,\bm{P}_{\epsilon}) - \epsilon \leq \min_{\bm{P} \in \U} (1-\gamma)R_{\gamma}(\pi,\bm{P}) \leq (1-\gamma)R_{\gamma}(\pi,\bm{P}_{\epsilon}), \forall \; \gamma > \gamma_{0}.
\end{equation}
Since $\pi \in \Pi_{\sf S}$ and $\bm{P}_{\epsilon} \in \U$, we know that $\lim_{\gamma \rightarrow 1} (1-\gamma)R_{\gamma}(\pi,\bm{P}_{\epsilon})$ exists and that \[\lim_{\gamma \rightarrow 1} (1-\gamma)R_{\gamma}(\pi,\bm{P}_{\epsilon}) = \Ravg(\pi,\bm{P}_{\epsilon}).\] This shows that for any $\epsilon>0$, we have 
\[ \limsup_{\gamma \rightarrow 1} (1-\gamma)\min_{\bm{P} \in \U} R_{\gamma}(\pi,\bm{P}) - \liminf_{\gamma \rightarrow 1} \min_{\bm{P} \in \U} (1-\gamma)R_{\gamma}(\pi,\bm{P}) \leq \epsilon.\]
Therefore, \[\limsup_{\gamma \rightarrow 1} \min_{\bm{P} \in \U} (1-\gamma)R_{\gamma}(\pi,\bm{P}) = \liminf_{\gamma \rightarrow 1} \min_{\bm{P} \in \U} (1-\gamma)R_{\gamma}(\pi,\bm{P}),\] i.e., $\lim_{\gamma \rightarrow 1} \min_{\bm{P} \in \U} (1-\gamma)R_{\gamma}(\pi,\bm{P})$ exists. For conciseness, let us write $R^{\pi} \in \R$ for this limit. Then by taking the limit as $\gamma \rightarrow 1$ in \eqref{eq:proof-app-1}, we have
\[ \Ravg(\pi,\bm{P}_{\epsilon}) - \epsilon \leq R^{\pi} \leq \Ravg(\pi,\bm{P}_{\epsilon}).\]
Since the inequality above is true for any $\epsilon>0$ and since $\Ravg(\pi,\bm{P}_{\epsilon}) \geq \inf_{\bm{P} \in \U} \Ravg(\pi,\bm{P})$, we proved that $\inf_{\bm{P} \in \U} \Ravg(\pi,\bm{P}) \leq R^{\pi}$. We now want to prove the converse inequality. To do this, let us consider a fixed $\bm{P} \in \U.$ We always have $(1-\gamma)R_{\gamma}(\pi,\bm{P}) \geq \inf_{\bm{P} \in \U} (1-\gamma)R_{\gamma}(\pi,\bm{P})$. Taking the limit as $\gamma \to 1$ on both sides, we have $\Ravg(\pi,\bm{P}) \geq R^{\pi}$. Taking the infimum over $\bm{P} \in \U$ on the right-hand side, we obtain that $\inf_{\bm{P} \in \U} \Ravg(\pi,\bm{P}) \geq R^{\pi}$. This shows that $\inf_{\bm{P} \in \U} \Ravg(\pi,\bm{P}) = R^{\pi}$, which concludes the first step.
\paragraph{Second step.} The second and third steps proves to Lemma \ref{lem:aux-2}. We now study the limit behavior of the optimal discounted return as $\gamma \rightarrow 1$. In particular, we show that $v$ exists, with \[v=\lim_{\gamma \rightarrow 1} \max_{\pi' \in \Pi_{\sf SD}} \inf_{\bm{P}' \in \U} (1-\gamma) R_{\gamma}(\pi',\bm{P}').\]
To show that $v$ exists, we consider $\pi$ a stationary deterministic policy that is $\epsilon$-Blackwell optimal for all $\epsilon>0$. Let $\epsilon>0$. We have, by the definition of $\pi$, that
\[ \min_{\bm{P} \in \U} (1-\gamma) R_{\gamma}(\pi,\bm{P}) + \epsilon \geq \max_{\pi' \in \Pi_{\sf SD}} \min_{\bm{P} \in \U} (1-\gamma) R_{\gamma}(\pi',\bm{P}) \geq \min_{\bm{P} \in \U} (1-\gamma) R_{\gamma}(\pi,\bm{P}), \forall \; \gamma > \gamma_{0}.\]
From the first step of our proof, $\lim_{\gamma \rightarrow 1} \min_{\bm{P} \in \U} R_{\gamma}(\pi,\bm{P})$ exists. Hence, for all $\epsilon>0$, we have
\[ \limsup_{\gamma \rightarrow 1} \max_{\pi' \in \Pi_{\sf SD}} \min_{\bm{P} \in \U} (1-\gamma) R_{\gamma}(\pi',\bm{P}) - \liminf_{\gamma \rightarrow 1} \max_{\pi' \in \Pi_{\sf SD}} \min_{\bm{P} \in \U} (1-\gamma) R_{\gamma}(\pi',\bm{P}) \leq \epsilon.\]
This shows that $v=\lim_{\gamma \rightarrow 1} \max_{\pi' \in \Pi_{\sf SD}} \min_{\bm{P} \in \U} (1-\gamma) R_{\gamma}(\pi',\bm{P})$ exists, and also that $v = \lim_{\gamma \rightarrow 1} \min_{\bm{P} \in \U} (1-\gamma) R_{\gamma}(\pi,\bm{P})$ for $\pi$ is an $\epsilon$-Blackwell optimal policy for any $\epsilon>0$.
\paragraph{Third step.} We now show that $v = \max_{\pi \in \Pi_{\sf SD}} \inf_{\bm{P} \in \U} \Ravg(\pi,\bm{P}).$
We start by showing that 
$ v \leq \max_{\pi \in \Pi_{\sf SD}} \min_{\bm{P} \in \U} \Ravg(\pi,\bm{P}).$
Indeed, let us consider the policy $\pi$ from Theorem \ref{th:sa-rec-epsilon-blackwell-opt-same}, that is stationary deterministic and $\epsilon$-Blackwell optimal, for all $\epsilon>0$:
\[ (1-\gamma) R_{\gamma}(\pi,\bm{P}) \geq \max_{\pi' \in \Pi_{\sf SD}} \min_{\bm{P}' \in \U} (1-\gamma) R_{\gamma}(\pi',\bm{P}') - \epsilon, \forall \; \epsilon>0,\forall \; \bm{P} \in \U,\forall \; \gamma > \gamma_{0}.\]
Since $\pi \in \Pi_{\sf SD}$, for $\bm{P} \in \U$ we have 
$ \lim_{\gamma \rightarrow 1} (1-\gamma) R_{\gamma}(\pi,\bm{P}) = \Ravg(\pi,\bm{P}).$
Hence for any $\bm{P} \in \U$ and $\epsilon>0$, we have
$\Ravg(\pi,\bm{P}) \geq v - \epsilon.$
Overall we have shown that $\inf_{\bm{P} \in \U} \Ravg(\pi,\bm{P}) \geq v$,
which implies that 
$\max_{\pi' \in \Pi_{\sf SD}} \inf_{\bm{P} \in \U} \Ravg(\pi',\bm{P}) \geq v.$

We now show that $\max_{\pi' \in \Pi_{\sf SD}} \inf_{\bm{P} \in \U} \Ravg(\pi',\bm{P}) \leq v$. We proceed by contradiction. Assume that there exists $\epsilon>0$ such that
\begin{equation}\label{eq:avg-rew-contradiction}
    \max_{\pi' \in \Pi_{\sf SD}} \inf_{\bm{P} \in \U} \Ravg(\pi',\bm{P}) > v + \epsilon
\end{equation}
and let $\pi$ be a stationary deterministic average optimal policy (which exists from Theorem \ref{th:avg-reward-main}).  Then
$\Ravg(\pi,\bm{P}) > v + \epsilon , \forall \; \bm{P} \in \U.$
Now let $\epsilon' < \epsilon$. From Theorem \ref{th:eps-blackwell-compact-mdp}, we know that there exists $\bm{P}_{\epsilon'} \in \U$ such that for any $\gamma > \gamma_{0}$, we have
\[ (1-\gamma) R_{\gamma}(\pi,\bm{P}_{\epsilon'}) \leq \min_{\bm{P} \in \U} (1-\gamma) R_{\gamma}(\pi,\bm{P}) + \epsilon' \leq \max_{\pi' \in \Pi_{\sf SD}} \min_{\bm{P} \in \U} (1-\gamma) R_{\gamma}(\pi',\bm{P}) + \epsilon'.\]
Taking the limit as $\gamma \rightarrow 1$ we obtain that
\[ \Ravg(\pi,\bm{P}_{\epsilon'}) \leq v + \epsilon' < v + \epsilon \leq \max_{\pi' \in \Pi_{\sf SD}} \inf_{\bm{P} \in \U} \Ravg(\pi', \bm{P}) \]
This is in contradiction with~\eqref{eq:avg-rew-contradiction}. Hence, we have shown that
$v = \max_{\pi' \in \Pi_{\sf SD}} \inf_{\bm{P} \in \U} \Ravg(\pi',\bm{P}).$
\paragraph{Fourth step.} The fourth step proves Lemma~\ref{lem:aux-3}. Let $\pi$ be a $\epsilon$-Blackwell optimal policy for any $\epsilon>0$. From the third step, we know that 
$\inf_{\bm{P} \in \U} \Ravg(\pi,\bm{P}) \geq v.$
From the third step, we also know that $v = \max_{\pi' \in \Pi_{\sf SD}} \inf_{\bm{P} \in \U} \Ravg(\pi',\bm{P}).$
Hence, we can conclude that \[\inf_{\bm{P} \in \U} \Ravg(\pi,\bm{P}) \geq \max_{\pi' \in \Pi_{\sf SD}} \inf_{\bm{P} \in \U} \Ravg(\pi',\bm{P}),\]
i.e., we can conclude that $\pi$ is optimal for the average return criterion. 

\tb{
\paragraph{Part 2:} Let $\pi\opt \in \PiSD$ be an average optimal policy. We will show that $\pi\opt$ is $\epsilon$-Blackwell optimal for all $\epsilon>0$.}

\tb{
We use Lemmas \ref{lem:aux-1}, \ref{lem:aux-2}, \ref{lem:aux-3} from the proof of Part 1.
In particular, let $\epsilon>0$. From Lemma \ref{lem:aux-1}, we know that there exists $\gamma_1 \in [0,1)$ such that 
\[|\min_{\bm{P} \in \U} (1-\gamma) R_{\gamma}(\pi\opt,\bm{P}) - \inf_{\bm{P} \in \U} \Ravg(\pi\opt,\bm{P})| \leq \epsilon/2,\forall \; \gamma \in (\gamma_1,1).\]
From Lemma \ref{lem:aux-2}, we know that there exists $\gamma_2 \in [0,1)$ such that
\[|\max_{\pi \in \Pi_{\sf SD}} \min_{\bm{P} \in \U} (1-\gamma)R_{\gamma}(\pi,\bm{P}) - \max_{\pi \in \Pi_{\sf SD}} \inf_{\bm{P} \in \U} \Ravg(\pi,\bm{P})| \leq \epsilon/2, \forall \; \gamma \in (\gamma_2,1).\]
Since $\pi$ is average optimal, we know that $\max_{\pi \in \Pi_{\sf SD}} \inf_{\bm{P} \in \U} \Ravg(\pi,\bm{P}) = \inf_{\bm{P} \in \U} \Ravg(\pi\opt,\bm{P})$.
Therefore,  we have
\[|\min_{\bm{P} \in \U} (1-\gamma) R_{\gamma}(\pi\opt,\bm{P}) - \max_{\pi \in \Pi_{\sf SD}} \min_{\bm{P} \in \U} (1-\gamma)R_{\gamma}(\pi,\bm{P})| \leq \epsilon, \forall \gamma \in (\max\{\gamma_1,\gamma_2\},1).\]
This shows that $\pi\opt$ is $\epsilon$-Blackwell optimal. Since the same proof works for any value of $\epsilon>0$, then $\pi$ is $\epsilon$-Blackwell optimal for all $\epsilon>0$, which concludes the proof of Theorem \ref{th:sa-rec-relation-avg-rew-blackwell}.}
\hfill \halmos \endproof

\tb{
\section{Proof of Theorem \ref{th:sa-rec-relation-avg-rew-discount}}\label{app:proof relation avg rew discount}}
\tb{\begin{repeattheorem}[Theorem 4.10]
    Consider an sa-rectangular robust MDP with a compact uncertainty set $\U$. There exists a discount factor $\gammaa \in [0,1)$ such that any stationary deterministic discount optimal policy for $\gamma \in (\gammaa,1)$ is also average optimal.
\end{repeattheorem}}
\tb{
\proof{Proof of Theorem \ref{th:sa-rec-relation-avg-rew-discount}.}
Let $\epsilon \geq 0$ be the minimum difference between the worst-case average return of an average optimal policy and any other non-average optimal {\em stationary deterministic} policy. Since $\PiSD$ is a finite set, we know that $\epsilon>0$.
}

\tb{
From Lemma \ref{lem:aux-2}, we know that $\lim_{\gamma \rightarrow 1} \max_{\pi \in \Pi_{\sf SD}} \min_{\bm{P} \in \U} (1-\gamma)R_{\gamma}(\pi,\bm{P}) = \max_{\pi \in \Pi_{\sf SD}} \inf_{\bm{P} \in \U} \Ravg(\pi,\bm{P}).$
By the definition of the limit, we can find $\gamma_1 \in [0,1)$ such that for all $\gamma \in (\gamma_1,1)$, we have
\begin{equation}\label{eq:prf-aux-1}
    |\max_{\pi \in \Pi_{\sf SD}} \min_{\bm{P} \in \U} \; (1-\gamma)R_{\gamma}(\pi,\bm{P}) - \max_{\pi \in \Pi_{\sf SD}} \inf_{\bm{P} \in \U} \Ravg(\pi,\bm{P})| \leq \epsilon/4.
\end{equation}
Now let $\pi \in \PiSD$. From Lemma \ref{lem:aux-1}, we know that $\lim_{\gamma \rightarrow 1} \min_{\bm{P} \in \U} \; (1-\gamma) R_{\gamma}(\pi,\bm{P}) = \inf_{\bm{P} \in \U} \Ravg(\pi,\bm{P}).$
By definition of the limit, we can find $\gamma_\pi$ such that for all $\gamma \in (\gamma_\pi,1)$ we have
\begin{equation}\label{eq:prf-aux-2}
|\min_{\bm{P} \in \U} (1-\gamma) R_{\gamma}(\pi,\bm{P}) - \inf_{\bm{P} \in \U} \Ravg(\pi,\bm{P})|\leq \epsilon/4.
\end{equation}
    Since $\PiSD$ is a finite set, we have that $\gamma_2 = \max_{\pi \in \PiSD} \gamma_\pi$ is such that $\gamma_2 \in [0,1)$ and \eqref{eq:prf-aux-2} holds for any $\gamma \in (\gamma_2,1)$ and any $\pi \in \PiSD$.
    }

\tb{Now let $\gammaa =\max\{\gamma_1,\gamma_2\}$ and let $\pi\opt \in \PiSD$ be discount optimal for $\gamma \in (\gammaa,1)$. We obtain
\begin{align*}
    \inf_{\bm{P} \in \U} \Ravg(\pi\opt,\bm{P}) & \geq \min_{\bm{P} \in \U} (1-\gamma) R_{\gamma}(\pi\opt,\bm{P}) - \epsilon/4 \\
    & = \max_{\pi \in \Pi_{\sf SD}} \min_{\bm{P} \in \U} (1-\gamma)R_{\gamma}(\pi,\bm{P}) + \epsilon/4 \\
    & \geq \max_{\pi \in \Pi_{\sf SD}} \inf_{\bm{P} \in \U} \Ravg(\pi,\bm{P}) - \epsilon/2
\end{align*}
where the first inequality follows from \eqref{eq:prf-aux-2}, the equality follows from $\pi\opt$ being discount optimal, and the second inequality follows from \eqref{eq:prf-aux-1}.
By the definition of $\epsilon>0$, this implies that $\pi\opt$ is average optimal, which concludes our proof.
}
    
\hfill \halmos \endproof

\section{Proof of Lemma \ref{lem:definability-stability}}\label{app:proof lemma definability stability}

\begin{repeattheorem}[Lemma 4.13]
            \begin{enumerate}
        \item The only definable subsets of $\R$ are the {\em finite} union of open intervals and singletons.
         \item If $A,B \subset \R^{n}$ are definable sets, then $A \cup B,A \cap B$ and $\R^{n} \setminus A$ are definable sets.
        \item Let $f,g$ be definable functions. Then $f \circ g$, $-f$, $f+g$, $f \times g$ are definable. 
        \item For $A,B$ two definable sets and $g: A \times B \rightarrow \R$ a definable function, then the functions $\bm{x} \mapsto \inf_{\bm{y} \in B} g(\bm{x},\bm{y})$ and $\bm{x} \mapsto \sup_{\bm{y} \in B} g(\bm{x},\bm{y})$ (defined over $A$) are definable functions. 
        \item 
    If $A,B,g:A\rightarrow \R$ are definable then $g^{-1}(B)$ is a definable set.
        \end{enumerate}
\end{repeattheorem}

The statements of Lemma \ref{lem:definability-stability} are either direct applications of results in \cite{bolte2015definable}, or textbook exercises from \cite{coste2000introduction}. We provide the precise reference below.
\proof{Proof of Lemma \ref{lem:definability-stability}.}
\begin{enumerate}
\item This follows from the axiomatic definition of definability~\citep{bolte2015definable}.
\item Similarly as the previous statement, this also follows from the definition of definability~\citep{bolte2015definable}.
\item Let $f,g:\R \rightarrow \R$ be two definable functions. The composition of two definable functions is still definable (exercise 1.11, \cite{coste2000introduction}). It is straightforward from the definition of definable functions that $-f$ is definable. Moreover,  the function $\phi:\R \rightarrow \R^{2},\phi(x) = (f(x),g(x))$ is definable (exercise 1.10, \cite{coste2000introduction}). The polynomial $P:\R^{2} \rightarrow \R,P(a,b) = a+b$ is definable. , hence $x \mapsto P \circ \phi(x)$ is definable, i.e., $x \mapsto f(x) + g(x)$ is definable.
    \item This is Proposition 1 in \cite{bolte2015definable}.
    \item This is example 1 in \cite{bolte2015definable}.
\end{enumerate}
\hfill \halmos \endproof

\tb{
\section{A large upper bound on $\gammab$}\label{app:issue-gamma-b}
For the sake of illustration, consider an MDP instance with rewards in $\{0,1\}$ and transition probabilities with entries in $\{0,1\}$, and $S$ states (the number of actions does not play a role). Then the bound in Theorem 4.4 in \cite{grand2024reducing} gives an upper bound $\bar{\gamma}$ on $\gamma_{\sf bw}$ such that
 \[ \bar{\gamma} = 1 - \frac{1}{2(2S-1)^{S+3/2}(2S4^S)^{2S-1}}.\]
 Note that the denominator is more than exponential in $S$, the number of states.
 As evident from this expression, $\bar{\gamma}$ is extremely close to $1$: for instance, for $S=2$, we obtain that $1-10^{-7}  < \bar{\gamma}$ and for $S=3$, we obtain $1-10^{-16} < \bar{\gamma}$.
 }

\section{Non-Lipschitzness of the robust value functions}\label{app:proof counterexample not lip}
We recall that a function $f:(0,1) \rightarrow \R^{n}$ is {\em Lipschitz continuous} if 
\begin{equation}\label{eq:lip-definition}
     \exists \; L > 0, \forall x,y \in (0,1), \| f(x) - f(y) \|_{\infty} \leq L |x-y|.
\end{equation}
Surprisingly, we show that robust value functions may not be Lipschitz continuous.
\begin{proposition}\label{prop:counter-ex-lip}
There exists an sa-rectangular definable uncertainty set $\U$ and a stationary deterministic policy $\pi$ such that the robust value function $\gamma \mapsto \bm{v}^{\pi,\U}_{\gamma}$ is not Lipschitz continuous and the normalized robust value function $\gamma \mapsto (1-\gamma)\bm{v}^{\pi,\U}_{\gamma}$ is not Lipschitz continuous.
\end{proposition}
 The proof builds upon our previous counterexample for the non-existence of Blackwell optimal policies from Proposition \ref{prop:inf-is-not-min-avg-rew}, which we adapt to show that the non-Lipschitzness of the value functions as in Proposition \ref{prop:counter-ex-lip}. Intuitively, in our counterexample the oscillations of the robust value function become more and more rapid as $\gamma$ approaches $1$, which makes them non-Lipschitz. 
\proof{Proof of Proposition \ref{prop:counter-ex-lip}}
    Consider the RMDP instance described in Section \ref{sec:RMDPs-average} for the proof of Proposition \ref{prop:inf-is-not-min-avg-rew}, where
    \[ \U_{s_{0}a_{0}} = \{ \left(1-\alpha,\beta,\alpha-\beta\right) \; | \; \alpha \in [0,1], \beta \in [0,\alpha],\beta \leq \alpha (1-\alpha)\}\]
    We have, for $\gamma \in (0,1)$,
    \[ R_{\gamma}(a_{0}) =\min_{\bm{P} \in \U} R_{\gamma}(a_{0},\bm{P}) = - \frac{1}{1-\gamma} \max_{\alpha \in [0,1]} \gamma \frac{\alpha(1-\alpha)}{1-\gamma+\gamma \alpha} = - \frac{1}{1-\gamma} \frac{1-2\sqrt{1-\gamma}+(1-\gamma)}{\gamma},\]
    where we use that the $\arg \max$ is attained at $\alpha\opt(\gamma) = \frac{\sqrt{1-\gamma}-(1-\gamma)}{\gamma}$.
    We now show that $\gamma \mapsto R_{\gamma}(a_{0})$ and $\gamma \mapsto (1-\gamma) R_{\gamma}(a_{0})$ are {\em not} Lipschitz continuous. From the definition of Lipschitz continuity as in~\eqref{eq:lip-definition}, a differentiable Lipschitz continuous function on $(0,1)$ has bounded derivatives on $(0,1)$. To show that the robust value functions and the normalized robust value functions are not Lipschitz continuous, we show that these functions are differentiable but their derivatives are unbounded in $(0,1)$.
    \begin{itemize}
        \item {\em The case of $\gamma \mapsto (1-\gamma)R_{\gamma}(a_3)$.} Let $f:\gamma \mapsto (1-\gamma)R_{\gamma}(a_3)$. In this case, $f(\gamma) = - \frac{1-2\sqrt{1-\gamma}+(1-\gamma)}{\gamma}.$
       Note that $f$ is differentiable in $(0,1)$ and
       \[ f'(\gamma) = -\frac{1}{\gamma} \left(\frac{1}{\sqrt{1-\gamma}}-1\right) + \frac{1-2\sqrt{1-\gamma}+(1-\gamma)}{\gamma^2}.\]
       We conclude that $f'(\gamma)  \sim \frac{-1}{\sqrt{1-\gamma}}$ as $\gamma \rightarrow 1$, so that $f'_{\gamma} \rightarrow - \infty$ as $\gamma \rightarrow 1$ and therefore $f$ is not Lipschitz continuous.
       \item {\em The case of $\gamma \mapsto R_{\gamma}(a_3)$.} Since $R(\gamma) = \frac{1}{1-\gamma}f(\gamma)$, we have 
       $R'_{\gamma}(a_3) = \frac{ f'(\gamma)(1-\gamma)+f(\gamma)}{(1-\gamma)^2} \sim \frac{-1}{(1-\gamma)^{3/2}}$ as $\gamma \rightarrow 1$. Therefore, $R'_{\gamma}(a_3) \rightarrow - \infty$ as $\gamma \rightarrow 1$ and $\gamma \mapsto R_{\gamma}(a_{3})$ is not Lipschitz continuous.
    \end{itemize}
\hfill \halmos \endproof
Despite their potential lack of Lipschitz continuity, robust value functions are always differentiable under some mild assumptions, as we show in the next theorem.
 Note that prior to this work, virtually nothing is known about the regularity of the robust value functions $\gamma \mapsto \bm{v}^{\pi,\U}_{\gamma}$ and of the optimal robust value function $\gamma \mapsto \bm{v}^{\star,\U}_{\gamma}.$ We recall that a function $f$ is of class $C^p$ for $p \in \N$ if $f$ is differentiable $p$-times and the $p$-th differential of $f$ is continuous.
\begin{theorem}\label{th:analytical-properties-value-function}
Consider an sa-rectangular robust MDP with a definable uncertainty set $\U$.
\begin{enumerate}
    \item Let $\pi \in \Pi_{\sf S}$ be a policy and $p \in \N$. Then there exists a finite subdivision of the interval $(0,1)$ as $0 = a_{1} < a_{2} < \dots < a_{k} = 1$ such that on each $(a_{i},a_{i+1})$ for $i=1,...,k-1$, the robust value function $\gamma \mapsto \bm{v}^{\pi,\U}_{\gamma}$ is of class $C^{p}$.
    \item Let $p \in \N$. Then there exists a finite subdivision of the interval $(0,1)$ as $0 = a_{1} < a_{2} < \dots < a_{k} = 1$ such that on each $(a_{i},a_{i+1})$ for $i=1,...,k-1$, the optimal robust value function $\gamma \mapsto \bm{v}^{\star,\U}_{\gamma}$ is of class $C^{p}$.
\end{enumerate}
\end{theorem} 
Theorem \ref{th:analytical-properties-value-function} is a straightforward consequence of a more general version of the monotonicity theorem (theorem 4.1 in \cite{van1996geometric}) and we omit the proof for conciseness.
Interesting properties related to H\"older continuity and \L ojasiewicz inequality can also be extended from Theorem 4.14 in \cite{van1996geometric}. 

\section{Proof of Theorem \ref{th:convergence-alg-definable-uset-2}}\label{app:proof convergence alg uset 2}

\begin{repeattheorem}[Theorem 5.3]
        Let $\U$ be a definable compact sa-rectangular uncertainty set and let $\left(\bm{v}^{t}\right)_{t \geq 1}$ be the iterates of Algorithm \ref{alg:sarec-increasing-discount-factor}. Then $\left(\bm{v}^{t}\right)_{t \geq 1}$ converges to $\bm{g}\opt \in \R^{\X}$.
\end{repeattheorem}

\proof{Proof of Theorem \ref{th:convergence-alg-definable-uset-2}.}
    Let $\pi^B \in \Pi_{\sf SD}$ be a Blackwell optimal policy. We have
\[ \| \bm{v}^{t} - \bm{g}\opt \|_{\infty} \leq \| \bm{v}^{t} - (1-\gamma_{t})\bm{v}^{\pi^B,\U}_{\gamma_{t}} \|_{\infty} + \| \bm{g}\opt -(1-\gamma_{t})\bm{v}^{\pi^B,\U}_{\gamma_{t}}\|_{\infty}.\]
The second term converges to $0$ as $\gamma_{t} \rightarrow 1$, as we proved in Theorem \ref{th:sa-rec-relation-avg-rew-blackwell}. In the rest of the proof, we show that the first term converges to $0$ as $t \rightarrow + \infty$.
We start with the following lemma, which is similar to lemma 7 in \cite{tewari2007bounded}.
\begin{lemma}\label{lem:interm-increasing-discount-1}
    There exists $t_{0} \in \N$ such that for any $t \geq t_{0}$, we have
    $\| \bm{v}^{t} - (1-\gamma_t)\bm{v}^{\pi^B,\U}_{\gamma_t} \|_{\infty} \leq \delta_{t}$
    with $\left(\delta_t\right)_{t \geq t_{0}}$ a sequence of scalars such that $\delta_{t+1} = \gamma_{t} \delta_{t} + e_{t}$ with $e_{t} = \| (1-\gamma_{t+1})\bm{v}^{\pi^B,\U}_{\gamma_{t+1}} - (1-\gamma_{t})\bm{v}^{\pi^B,\U}_{\gamma_{t}} \|_{\infty}.$
\end{lemma}
\proof{Proof of Lemma \ref{lem:interm-increasing-discount-1}}
Note that $\gamma_t \rightarrow 1$ as $t \rightarrow + \infty$. Therefore, we can choose $t_{0}$ large enough such that $\pi^B$ is an optimal policy, see Theorem \ref{th:blackwell-discount-factor}. We prove Lemma \ref{lem:interm-increasing-discount-1} by induction. Assume that  $\| \bm{v}^{t} - (1-\gamma_t)\bm{v}^{\pi^B,\U}_{\gamma_t} \|_{\infty} \leq \delta_t.$
We have
\begin{align*}
    \| \bm{v}^{t+1} - (1-\gamma_{t+1})\bm{v}^{\pi^B,\U}_{\gamma_{t+1}} \|_{\infty} & \leq \| \bm{v}^{t+1} - (1-\gamma_{t})\bm{v}^{\pi^B,\U}_{\gamma_{t}} \|_{\infty} + \| (1-\gamma_{t})\bm{v}^{\pi^B,\U}_{\gamma_{t}} - (1-\gamma_{t+1})\bm{v}^{\pi^B,\U}_{\gamma_{t+1}} \|_{\infty}.
\end{align*}
By definition, $e_{t} = \| (1-\gamma_{t})\bm{v}^{\pi^B,\U}_{\gamma_{t}} - (1-\gamma_{t+1})\bm{v}^{\pi^B,\U}_{\gamma_{t+1}} \|_{\infty}$. We now turn to bounding $\| \bm{v}^{t+1} - (1-\gamma_{t})\bm{v}^{\pi^B,\U}_{\gamma_{t}} \|_{\infty}$. We have, for any $s \in \X$,  $(1-\gamma_{t})v^{\pi^B,\U}_{\gamma_{t},s}  = (1-\gamma_t) T_{s,\gamma_t}\left(\bm{v}^{\pi^B,\U}_{\gamma_{t}}\right)$ because $\bm{v}^{\pi^B}_{\gamma_t,\U}$ is the fixed-point of the operator $T_{\gamma_t}$ since $\pi^B$ is Blackwell optimal. This shows that 
\[
    (1-\gamma_{t})v^{\pi^B,\U}_{\gamma_{t},s} =  \max_{a \in \A} \min_{\bm{p}_{sa} \in \U_{sa}} \bm{p}_{sa}\tr\left((1-\gamma_t)\bm{r}_{sa} + \gamma_{t}(1-\gamma_t)\bm{v}^{\pi^B,\U}_{\gamma_{t}}\right) \leq \max_{a \in \A} \min_{\bm{p}_{sa} \in \U_{sa}} \bm{p}_{sa}\tr\left((1-\gamma_t)\bm{r}_{sa} + \gamma_{t}(1-\gamma_t)\bm{v}^{t}\right) + \gamma_{t} \delta_t \]
where the first equality is by definition of $T$, and the inequality is because of the induction hypothesis and $\sum_{s' \in \X} p_{sas'} = 1, \forall \; \bm{p}_{sa} \in \U_{sa}$. Overall, we have proved that $(1-\gamma_{t})v^{\pi^B,\U}_{\gamma_{t},s} \leq v^{t+1}_{s} + \gamma_t \delta_t.$
We can prove the converse inequality similarly. From this, we conclude that $\| \bm{v}^{t} - (1-\gamma_t)\bm{v}^{\pi^B,\U}_{\gamma_t} \|_{\infty} \leq \delta_{t}$
    with $\delta_{t+1} = \gamma_t \delta_t + e_t$.
\hfill \halmos \endproof
We now turn to proving that $\delta_{t} \rightarrow 0$ as $t \rightarrow + \infty$.
From the induction $\delta_{t+1} = \gamma_t \delta_t + e_t$ we obtain that for any $t \geq t_{0}$ we have
\[ \delta_{t} = \sum_{j=t_{0}}^{t} \left(\prod_{i=j+1}^{t} \gamma_{i} \right) e_{j}.\]
Additionally, $\gamma_{i} = \frac{\omega_{i}}{\omega_{i+1}}$, therefore
\[\prod_{i=j+1}^{t} \gamma_{i} = \frac{\omega_{t}}{\omega_{t+1}} \frac{\omega_{t-1}}{\omega_{t}} \dots \frac{\omega_{j+2}}{\omega_{j+3}} \frac{\omega_{j+1}}{\omega_{j+2}} = \frac{\omega_{j+1}}{\omega_{t+1}}\] so that 
\[ \delta_t = \frac{1}{\omega_{t+1}} \sum_{j=t_0}^{t} \omega_{j+1} e_{j}.\]
We now prove that $\lim_{t \rightarrow + \infty} \frac{1}{\omega_{t+1}} \sum_{j=t_0}^{t} \omega_{j+1} e_{j} =0$. 
The most fundamental point is to note that from \cite{bolte2015definable} (Theorem 3 and Corollary 4), we have $\sum_{t=t_0}^{+\infty}  e_{t} < +\infty$ for any increasing sequence $\left(\gamma_t\right)_{t \in \N}$, see also \cite{renault2019tutorial}, corollary 1.6, point (2) for an equivalent result for SGs with finitely many actions.

Let us now write $E_j = \sum_{i = j}^{+\infty} e_j$.
 Since $\sum_{j=t_0}^{+\infty} e_j < +\infty$ and $e_j \geq 0$, we know that $E_j \rightarrow 0$ as $t \rightarrow + \infty$. Now we have
  \[  \sum_{j=t_{0}}^{t} \omega_{j+1} e_j  = \sum_{j=t_{0}}^{t} \omega_{j+1} \left(E_j - E_{j+1}\right) 
     = \omega_{t_{0}+1}E_{t_{0}} - \omega_{t+1}E_{t+1} + \sum_{j=t_0+1}^{t} (\omega_{j+1}-\omega_{j})E_j.\]
First, since $\omega_{t_{0}+1}E_{t_{0}}$ is constant and $\lim_{t \rightarrow + \infty} \omega_{t} = + \infty$, it is clear that $\omega_{t_{0}+1}E_{t_{0}}/\omega_{t+1} \rightarrow 0.$ Second, $\omega_{t+1}E_{t+1} /\omega_{t+1} = E_{t+1} \rightarrow 0$. There remains to show that $\frac{\sum_{j=t_0+1}^{t} (\omega_{j+1}-\omega_{j})E_j}{\omega_{t+1}} \rightarrow 0$ as $t \rightarrow + \infty$.
This can be done using the following simple lemma from real analysis (see theorem 2.7.2 in \cite{choudary2014real}). 
\begin{lemma}[Stolz-Cesaro theorem]\label{lem:real-analysis}
    Let $(A_{t})_{t \in \N}$ and $\left(B_{t}\right)_{t \in \N}$ two sequences of real numbers, with $\left(B_{t}\right)_{t \in \N}$ increasing with $\lim_{t \rightarrow + \infty} B_{t} = + \infty$. Let $\ell =\lim_{t \rightarrow + \infty} (A_{t+1}-A_{t})/(B_{t+1}-B_{t}), \ell \in \R \cup \{ + \infty\}$.
    Then the limit $\lim_{t \rightarrow + \infty} A_{t}/B_{t}$ exists and $\lim_{t \rightarrow + \infty} A_{t}/B_{t} = \ell$.
\end{lemma}
Applying Lemma \ref{lem:real-analysis} to $A_{t} = \sum_{j=t_{0}+1}^{t} (\omega_{j+1}-\omega_{j})E_j $ and $B_{t} = \sum_{j=t_{0}+1}^{t} \omega_{j+1}-\omega_{j} = \omega_{t+1} - \omega_{t_{0}+1}$, we obtain that
\[ \frac{A_{t+1}-A_{t}}{B_{t+1}-B_{t}} = \frac{(\omega_{t+2} - \omega_{t+1})E_{t+1}}{\omega_{t+2}-\omega_{t+1}} = E_{t+1} \rightarrow 0 \]
which shows that $A_{t}/B_{t} = \frac{\sum_{j=t_0+1}^{t} (\omega_{j+1}-\omega_{j})E_j}{\omega_{t+1}}$ converges to $0$.
 Overall, we have proved that $\delta_{t} \rightarrow 0$ as $t \rightarrow + \infty$, which concludes the proof of Theorem \ref{th:convergence-alg-definable-uset-2}.
\hfill \halmos \endproof

\section{Details on our numerical experiments}\label{app:simu}
\subsection{Nominal MDP instances}\label{app:simu-instance}
We start by describing the nominal transition probabilities and instantaneous rewards for the three MDP instances considered in our numerical experiments.
\subsubsection{Forest management instance.}
The states represent the growth of the forest. An optimal policy finds the right balance between maintaining the forest, earning revenue by selling cut wood, and the risk of wildfires. A complete description may be found at \cite{pymdp}.  This instance is inspired from the application of dynamic programming to optimal fire management \cite{possingham1997application}.
\paragraph{States.} There are $S$ states:  State $1$ is the youngest state for the forest, State $S$ is the oldest state. 
\paragraph{Actions.} The two actions are  \textit{wait} and \textit{cut \& sell}.
\paragraph{Transitions.} If the forest is in State $s$ and the action is \textit{wait}, the next state is State $s+1$ with probability $1-p$ (the forest grows) and $1$ with probability $p$ (a wildfire burns the forest down). If the forest is in State $s$ and the action is \textit{cut \& sell}, the next state is State $1$ with probability $1$. The probability of wildfire $p$ is chosen at $p=0.1$.
\paragraph{Rewards.} There is a reward of $4$ when the forest reaches the oldest state ($S$) and the chosen action is \textit{wait}.  There is a reward of $0$ at every other state if the chosen action is \textit{wait}. When the action is \textit{cut \& sell}, the reward at the youngest state $s=1$ is $0$, there is a reward of $1$ in any other state $s \in \{1, ..., S-1\}$, and a reward of $2$ in $s=S$. 
\subsubsection{Machine replacement instance}
This instance is represented in Figure \ref{fig:machine-instance}, in the case of $8$ states representing the operative conditions of the machine.  There are additional two repair states. To build larger instances, we add new operative states for the machine, with the same transition. We now describe in detail the states, actions, transitions, and rewards in this instance.
\paragraph{States.} The machine replacement problem involves a machine whose set of possible conditions are described by $S$ states.  The first $S-2$ states are operative states. The states $1$ to $S-2$ model the condition of the machine,  with $1$ being the perfect condition and $S-2$ being the worst condition.  The last two states $S-1$ and $S$ are states representing when the machine is being repaired. 

\paragraph{Actions.} There are two actions: {\sf repair} and {\sf wait}.

\paragraph{Transitions.} The transitions are detailed in Figures \ref{fig:Machine-MDP-1}-\ref{fig:Machine-MDP-2}.  When the action is {\sf wait}, the machine is likely to deteriorate toward the state $S-2$, or may stay in the same condition. When the action is {\sf repair}, the decision-maker brings the machine to the states $S-1$ and $S-2$.

\paragraph{Rewards.}There is a cost of 0 for states $1, ..., S-3$; letting the machine reach the worst operative state $S-2$ is penalized with a cost of $20$. 
The state $S-1$ is a standard repair state and has a cost of 2, while the last state $S$ is a longer and more costly repair state and has a cost of 10. We turn costs into rewards by flipping the signs.
\begin{figure}[htb]
\begin{center}
    \begin{subfigure}{0.4\textwidth}
\includegraphics[width=\linewidth]{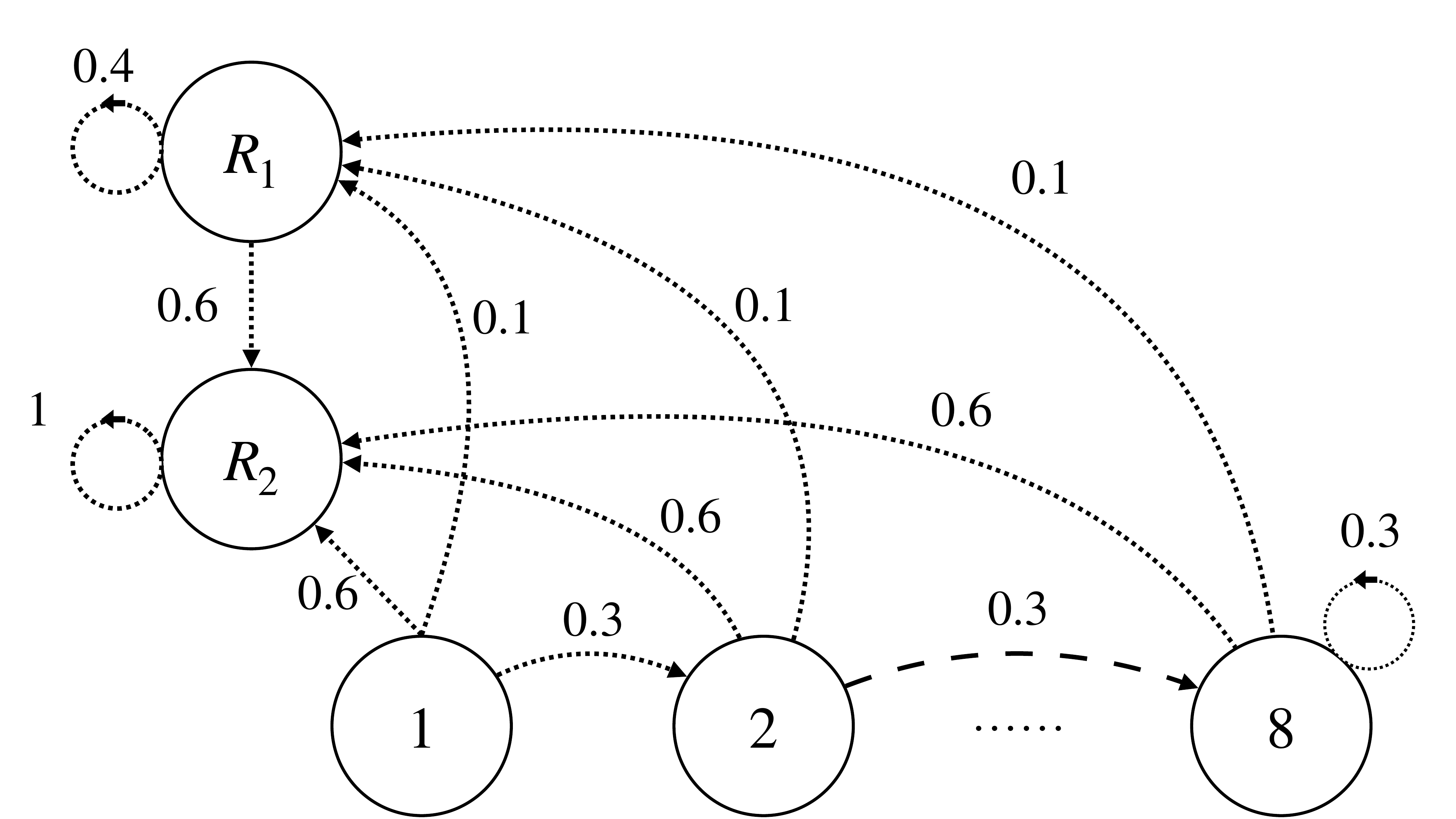}
         \caption{Action {\sf repair}.}
           \label{fig:Machine-MDP-1}
    \end{subfigure}
    \begin{subfigure}{0.4\textwidth}
\includegraphics[width=\linewidth]{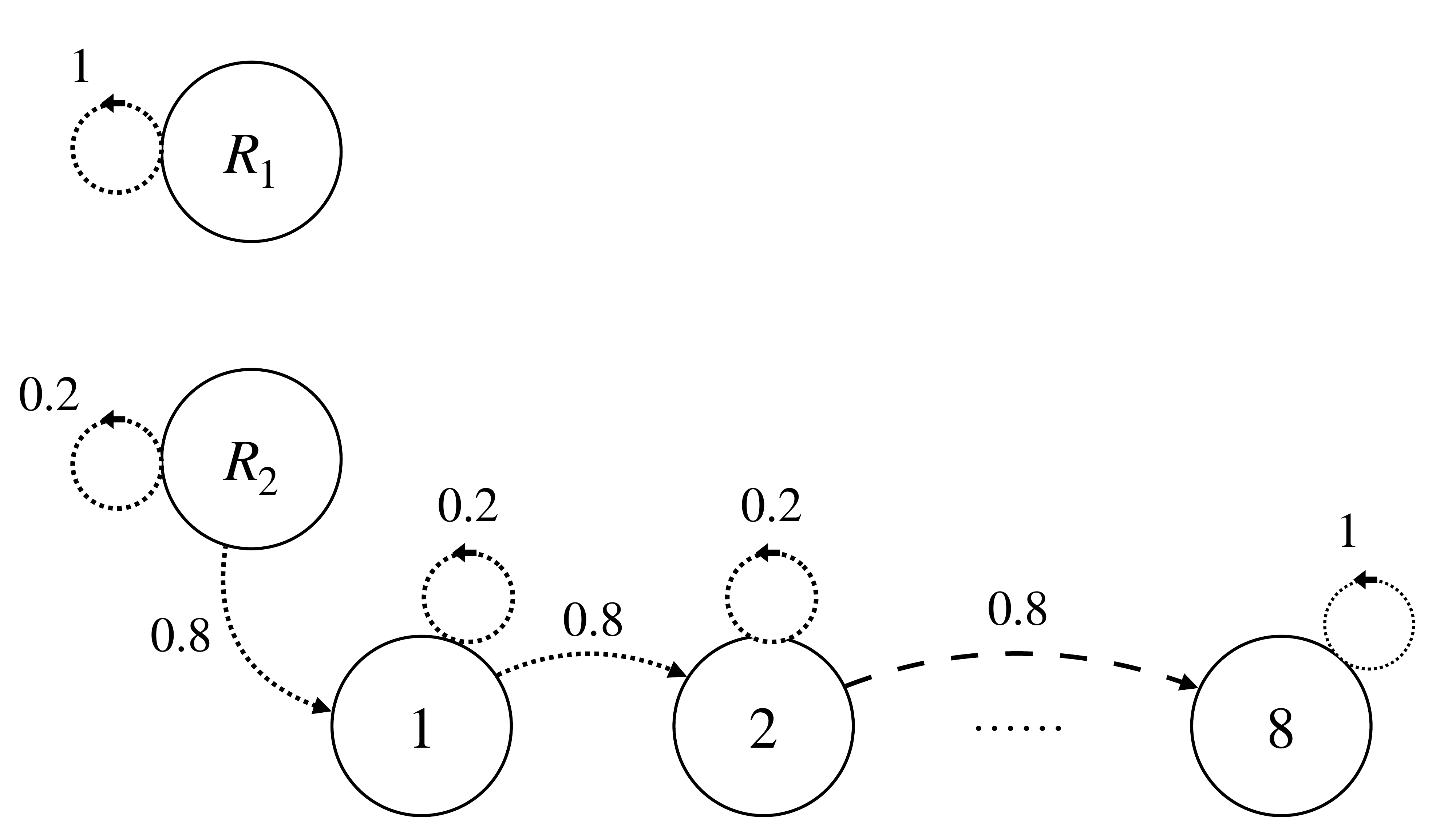}
         \caption{Action {\sf wait}.}
         \label{fig:Machine-MDP-2}
    \end{subfigure}
    \end{center}
\caption{Nominal transitions for the machine replacement instance.}    
\label{fig:machine-instance}
\end{figure}

\subsubsection{Instance inspired by healthcare}
We present in Figure \ref{fig:healthcare-instance} the nominal transitions for this instance, in the case of $8$ health conditions states. This goal is to minimize the mortality rate of the patient while reducing the invasiveness of the drug dosage (low, medium or high) prescribed at each state.
Similarly, as for the machine replacement instance, the instances for a larger number of states are constructed by adding some health condition states for the patient.

\paragraph{States.} There are $S$ states. The first $S-1$ states represent the health conditions of the patient, with $S-1$ being the worst condition before the mortality absorbing state $m$.
\paragraph{Actions.} The action set is $\A = \{\sf low, medium, high\}$, representing the drug dosage at each state. 
\paragraph{Transitions.} The transitions are represented in Figure \ref{fig:healthcare-instance}. They capture the fact that the patient is more likely to recover (i.e., transitioning early states) under the high-intensity treatment (Figure \ref{fig:healthcare-3}) than under the low-intensity treatment (Figure \ref{fig:healthcare-1}).
\paragraph{Rewards.} The rewards penalize using an intense treatment plan for the patient. In a health condition state $s \in \{1,...,S-1\}$, the reward is $10$ for choosing action {\sf low}, $8$ for choosing action {\sf medium}, and $6$ for choosing action {\sf high.}
\begin{figure}[htb]
\begin{center}
    \begin{subfigure}{0.3\textwidth}
\includegraphics[width=\linewidth]{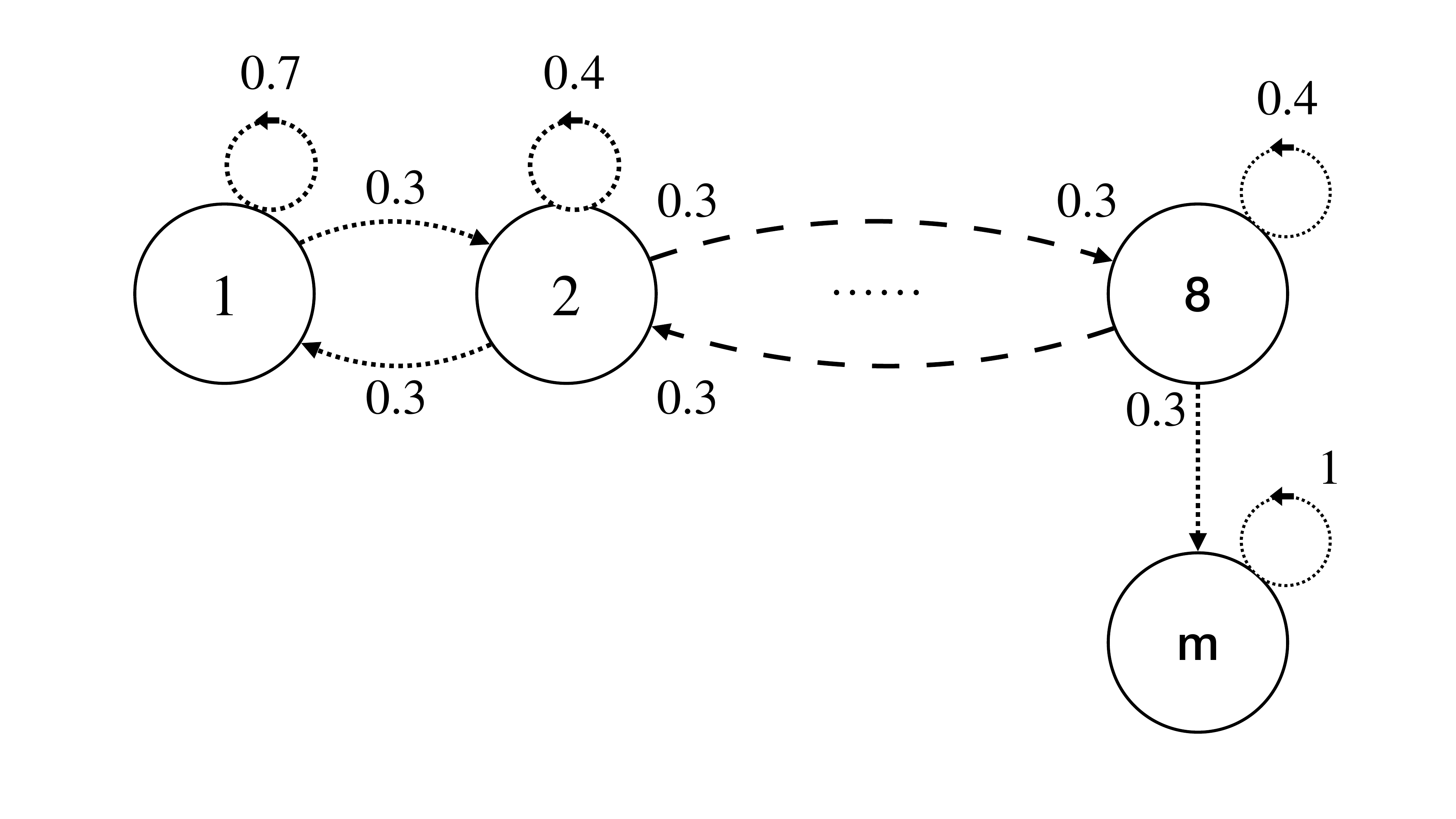}
         \caption{Action {\sf low}.}
           \label{fig:healthcare-1}
    \end{subfigure}
    \begin{subfigure}{0.3\textwidth}
\includegraphics[width=\linewidth]{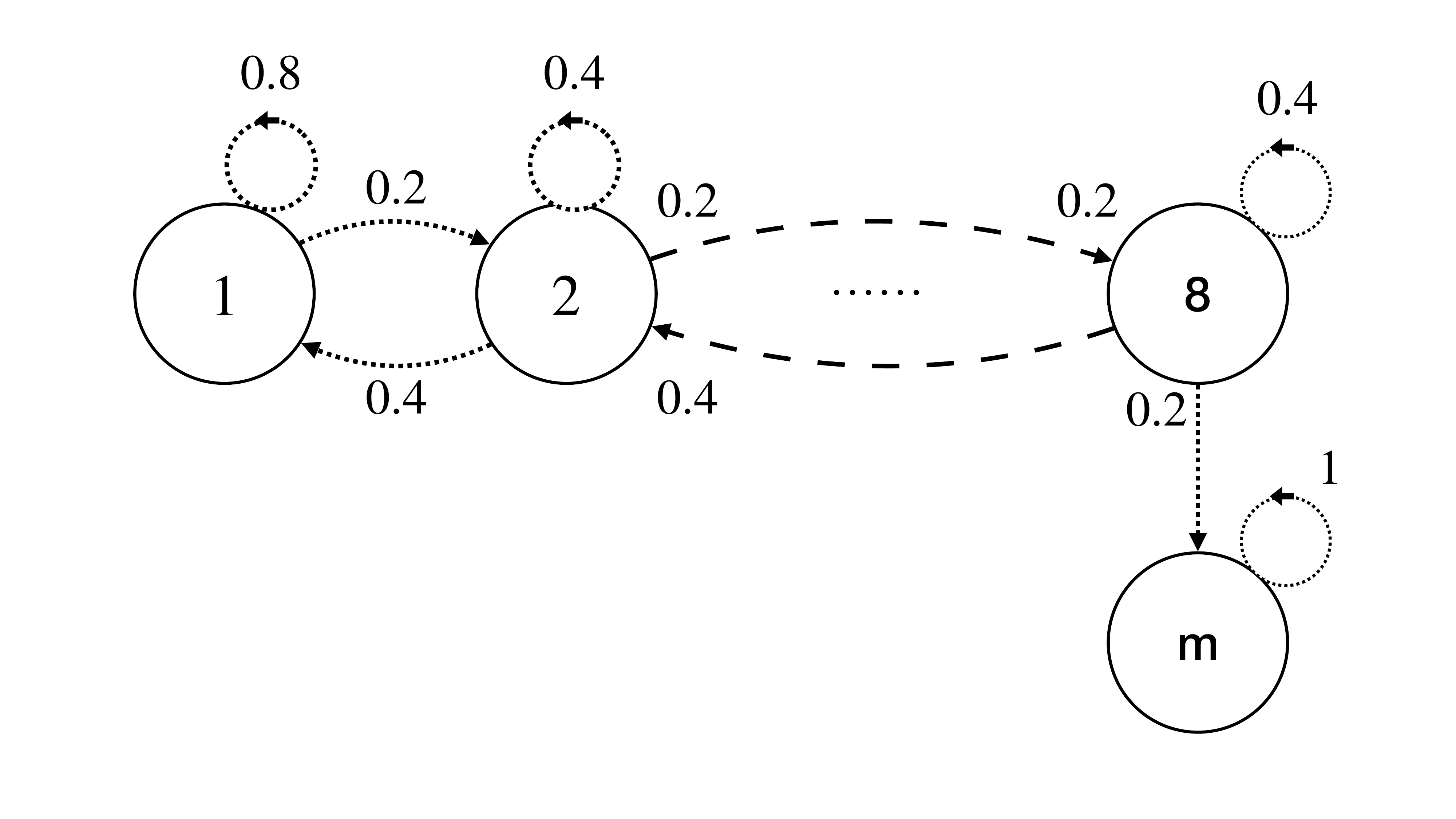}
         \caption{Action {\sf medium}.}
         \label{fig:healthcare-2}
    \end{subfigure}
    \begin{subfigure}{0.3\textwidth}
\includegraphics[width=\linewidth]{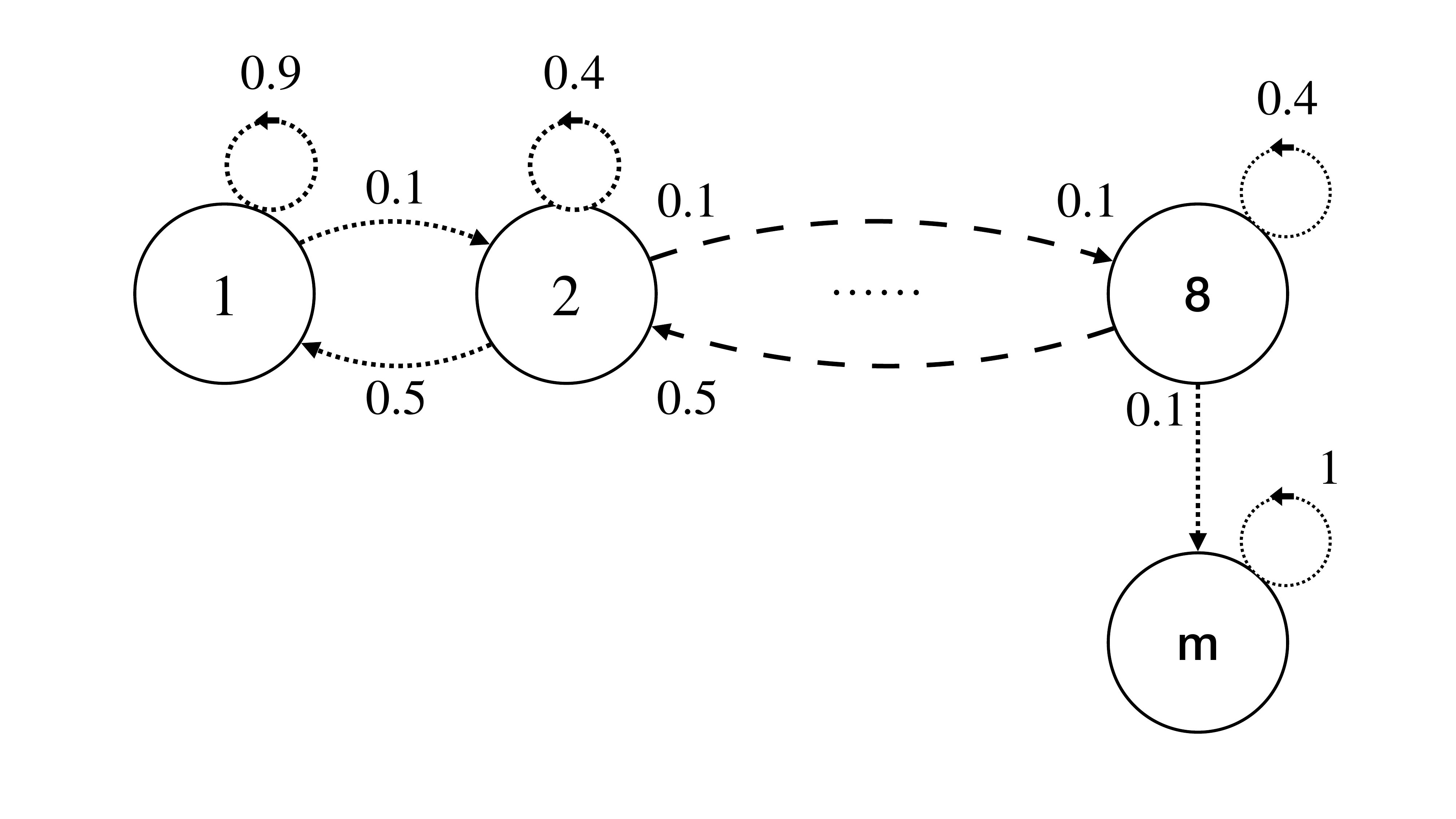}
         \caption{Action {\sf high}.}
         \label{fig:healthcare-3}
    \end{subfigure}
    \end{center}
\caption{Nominal transitions for the MDP instance inspired from healthcare.}    
\label{fig:healthcare-instance}
\end{figure}
\subsubsection{Garnet MDPs} Garnet MDPs are a class of random MDP instances~\citep{archibald1995generation}.In Garnet MDPs, only $N_{\sf b}$ states are reachable from any state-action pair $(s,a)$, and the transitions toward the next states are sampled uniformly at random. The rewards are sampled uniformly at random in $[0,1]$, and we choose $N_{\sf b}=10$ and $A=5$ actions in our numerical experiments.

\subsection{Uncertainty sets and practical implementation}\label{app:simu-uncertainty sets}
\subsubsection{The case of special transitions.}\label{app:simu-special-transitions}
There are absorbing states in the three instances used in our experiments. In non-absorbing states, some transition probabilities toward other states may be $0$, e.g., in the machine replacement instance, transitioning from one of the operative states to a repair state without repairing. Since these transitions may reflect real physical constraints, in the implementation of our algorithms we ensure that the resulting worst-case transition probabilities have the same support as the nominal transitions, e.g., that transitions that are impossible for the nominal transition probabilities remain impossible for the worst-case transitions, and that absorbing states for the nominal transitions remain absorbing for the worst-case transitions.

\subsubsection{Uncertainty based on ellipsoids.}
The ellipsoidal uncertainty set is 
$\U_{sa}^{\sf \ell_2} = \{ \bm{p} \in \Delta(\X) \; | \; \| \bm{p} - \bm{p}_{{\sf nom},sa}\|_{2} \leq \alpha \}, \forall \; (s,a) \in \X \times \A$
where $\alpha>0$ is a scalar. We use a radius $\alpha=0.05$.
To evaluate the Bellman operator $T$ at a vector $\bm{v} \in \R^{\X}$, we would need to solve the following convex program: $\min_{\bm{p} \in \U_{sa}} \bm{p}\tr\bm{v}'$ with $\bm{v}'=(\bm{r}_{sa}+\gamma\bm{v})$, which may considerably increase the computation time of our algorithms, where the Bellman operator needs to be evaluated at every iteration. To ensure a fast implementation, we make the following assumption, building up on the analysis in \cite{grand2021conic,grand2023solving}.
\begin{assumption}\label{ass:assumption-small-radius}
    The radius $\alpha>0$ is such that, for any $(s,a) \in \X \times \A$, we have
    \[ \{ \bm{p} \in \R^{n} \; | \; \bm{p}\tr\bm{e}=1\} \cap \{ \bm{p} \in \R^{\X} \; | \| \bm{p} - \bm{p}_{sa} \|_{2} \leq \alpha \} \subset \Delta(\X).\]
\end{assumption}
Assumption \ref{ass:assumption-small-radius} has the simple interpretation that the orthogonal projection of the $\ell_{2}$-ball $\{ \bm{p} \in \R^{\X} \; | \| \bm{p} - \bm{p}_{sa} \|_{2} \leq \alpha \}$ onto the hyperplane $\{ \bm{p} \in \R^{n} \; | \; \bm{p}\tr\bm{e}=1\}$ is entirely contained into the simplex. Note that from Section \ref{app:simu-special-transitions}, we only consider that transition probabilities in $(0,1)$ are uncertain so that Assumption \ref{ass:assumption-small-radius} holds for a radius $\alpha$ small enough. Under this assumption, we can write $\U_{sa}^{\ell_2}$ as $ \U_{sa}^{\ell_2}  = \bm{p}_{{\sf nom},sa} + \alpha \tilde{B}$ with $\tilde{B}  = \{ \bm{z} \in \R^{\X} \; | \; \bm{z}\tr\bm{e} = 0, \| \bm{z} \|_{2} \leq 1\}. $
Let us now write $\bm{u}_{1},...,\bm{u}_{|\X|-1} \in \R^{|\X|}$ as an orthonormal basis of $\{ \bm{p} \in \R^{n} \; | \; \bm{p}\tr\bm{e}=0\}$, e.g. the one given in closed-form in \cite{egozcue2003isometric}. Then writing $\bm{U} = \left(\bm{u}_{1},...,\bm{u}_{|\X|-1}\right) \in \R^{|\X| \times (|\X|-1)}$, and noting by definition $\bm{U}\tr\bm{U} = \bm{I}_{|\X|-1}$, we have
\[
\tilde{B} = \{ \bm{z} \in \R^{\X} \; | \; \bm{z}\tr\bm{e} = 0, \| \bm{z} \|_{2} \leq 1\} = \{ \bm{Uy} \; | \; \bm{y} \in \R^{|\X|-1}, \|\bm{Uy}\|_{2} \leq 1 \} = \{ \bm{Uy}  \; | \; \bm{y} \in \R^{|\X|-1}, \|\bm{y}\|_{2} \leq 1 \}\]
where the last equality follows from $\bm{U}$ being an orthonormal basis of $\{ \bm{p} \in \R^{n} \; | \; \bm{p}\tr\bm{e}=0\}$.
This shows that under Assumption \ref{ass:assumption-small-radius}, we have the following closed-form update: for $\bm{v} \in \R^{\X}$,
\[\min_{\bm{p} \in \U_{sa}} \bm{p}\tr\bm{v} = \bm{p}_{{\sf nom},sa}\tr\bm{v} + \alpha \min_{\bm{y} \in \R^{|\X|-1}, \| \bm{y} \|_{2} \leq 1} \left(\bm{Uy}\right)\tr\bm{v}  = \bm{p}_{{\sf nom},sa}\tr\bm{v} - \alpha \| \bm{U}\tr\bm{v}\|_{2}.\]
Overall, our analysis shows that in the case of ellipsoidal uncertainty as in $\U^{\sf \ell_{2}}$ and under Assumption \ref{ass:assumption-small-radius}, we can efficiently evaluate $T_{s}(\bm{v})$ as
$T_{s}(\bm{v}) = \max_{a \in \A} \bm{p}_{{\sf nom},sa}\tr\left(\bm{r}_{sa} + \gamma \bm{v}\right) - \alpha \| \bm{U}\tr\left(\bm{r}_{sa} + \gamma \bm{v}\right)\|_{2}, \forall \; s \in \X.$
\subsubsection{Uncertainty based on box inequalities.}
We also consider the following uncertainty set:
$\U_{sa}^{\sf box} = \{ \bm{p} \in \Delta(\X) \; | \; \bm{p}_{sa}^{\sf low} \leq \bm{p} \leq \bm{p}_{sa}^{\sf up} \}, \forall \; (s,a) \in \X \times \A$
with $\bm{p}_{sa}^{\sf low},\bm{p}_{sa}^{\sf up} \in [0,1]^{\X}$ two vectors such that $\bm{p}_{sa}^{\sf low} \leq \bm{p}_{{\sf nom},sa} \leq \bm{p}_{sa}^{\sf up}.$
We build $\bm{p}_{sa}^{\sf low},\bm{p}_{sa}^{\sf up}$ as follows. We choose two scalars $\theta_{\sf up}, \theta_{\sf low}$. The vector $\bm{p}_{sa}^{\sf low}$ is then the convex combination of $\bm{p}_{{\sf nom},sa}$ and $\bm{0}$ with coefficient $\theta_{\sf low}$, and $\bm{p}_{sa}^{\sf up}$ is then the convex combination of $\bm{p}_{{\sf nom},sa}$ and $\bm{e}$ with coefficient $\theta_{\sf up}$: $\bm{p}^{\sf up}_{sa} = (1-\theta_{\sf up}) \bm{p}_{{\sf nom},sa} + \theta_{\sf up}\bm{e},
    \bm{p}^{\sf low}_{sa} = (1-\theta_{\sf low}) \bm{p}_{{\sf nom},sa}.$ In our numerical experiment we use $\theta_{\sf up} = \theta_{\sf low} = 0.05.$
    Note that to evaluate the Bellman operator $T$ at a vector $\bm{v} \in \R^{\X}$, we only need to sort the component of $\bm{v}$ in increasing order and then use the closed-form expression from the method described in proposition 3 in \cite{goh2018data}. 
\
\subsection{Implementation of Algorithm \ref{alg:sarec-limit-discount-values}}\label{app:simu-2pl-PI}
In Algorithm \ref{alg:sarec-limit-discount-values}, at every iteration $t$ we compute $\bm{v}^{\star,\U}_{\gamma_{t}}$ by implementing two-player Strategy Iteration~\cite{hansen2013strategy}, as described in Algorithm \ref{alg:policy-iteration-2pl}.

\begin{algorithm}
\caption{Two-player Strategy Iteration for sa-rectangular robust MDPs}\label{alg:policy-iteration-2pl}
\begin{algorithmic}[1]
\State {\bf Input:} Initial policy $\pi^{\sf start} \in \Pi_{\sf SD}$

\State Initialize $\pi^{0} = \pi^{\sf start} \in \Pi_{\sf SD}$.

\For{$t \in \N$}
    \State Policy Evaluation: $\bm{v}^{t} = \bm{v}^{\pi^{t},\U}_{\gamma}$
    \State Policy Improvement: choose $\pi^{t+1} \in \Pi_{\sf SD}$ such that
    \[ \pi^{t+1}(s) \in \arg \max_{a \in \A} \min_{\bm{p}_{sa} \in \U_{sa}} \bm{p}_{sa} \tr \left(\bm{r}_{sa} + \gamma \bm{v}^{t}\right), \forall \; s \in \X,\]
    with $\pi^{t+1}(s) = \pi^{t}(s)$ if possible.
    \State {\sf stop} if $\pi^{t+1} = \pi^{t}$.
\EndFor
\end{algorithmic}
\end{algorithm}
Recall that in Algorithm \ref{alg:sarec-limit-discount-values}, we compute the discount optimal policies for an increasing sequence of discount factors. Therefore, we warm-start Algorithm \ref{alg:policy-iteration-2pl} at the next iteration with the policy computed at the previous iteration. This considerably reduces the computation time, for two reasons: 
\begin{enumerate}
    \item The discount factors $\gamma_t$ at iteration $t$ and $\gamma_{t+1}$ at iteration $t+1$ are close so that we can expect the optimal policies at iteration $t$ and iteration $t+1$ to be close; 
    \item Since $\U^{\sf box}$ and $\U^{\ell_2}$ are definable sets, Blackwell optimal policies exist. A Blackwell discount factor also exists, and for $\gamma$ large enough, the set of optimal stationary deterministic policies does not change.
\end{enumerate}
In Algorithm \ref{alg:policy-iteration-2pl}, we compute the robust value function of $\pi^t$ at every iteration $t$. We do so by implementing the one-player version of Policy Iteration for the adversarial MDP~\citep{ho2021partial,goh2018data}) as in Algorithm \ref{alg:policy-iteration-1pl}.
\begin{algorithm}
\caption{Policy Iteration for the Adversarial MDPs}\label{alg:policy-iteration-1pl}
\begin{algorithmic}[1]
\State {\bf Input:} Stationary policy $\pi \in \Pi_{\sf S}$, initial transition probabilities $\bm{P}^{\sf start} \in \U$

\State Initialize $\bm{P}^{0} = \bm{P}^{\sf start}$.

\For{$t \in \N$}
    \State Adversarial Policy Evaluation: $\bm{v}^{t} = \bm{v}^{\pi,\bm{P}^{t}}_{\gamma}$
    \State Adversarial Policy Improvement: choose $\bm{P}^{t+1} \in \U$ such that
    \[\bm{p}^{t+1}_{sa} \in \arg \min_{\bm{p}_{sa} \in \U_{sa}} \bm{p}_{sa}\tr\left(\bm{r}_{sa} + \gamma \bm{v}^{t}\right), \forall \; (s,a) \in \X \times \A\]
    with $\bm{P}^{t+1}_{sa} = \bm{P}^{t}_{sa}$ if possible.
    \State {\sf stop} if $\bm{p}^{t+1} = \bm{p}^{t}$.
\EndFor
\end{algorithmic}
\end{algorithm}
We also use warm-starting to accelerate the practical performance of Algorithm \ref{alg:policy-iteration-1pl}.

\end{APPENDICES}

\end{document}